\documentclass[11pt, a4paper]{amsart}
\usepackage{amssymb,amsfonts, amsmath, amsthm}
\usepackage[all,arc]{xy}
\usepackage{enumerate}
\usepackage{mathrsfs}
\usepackage{mathtools}
\usepackage[mathscr]{euscript}
\usepackage{array}
\usepackage{enumitem}
\usepackage{tikz}
\usepackage{tikz-cd}
\usepackage{textcomp}
\usepackage{mathpazo}

\usetikzlibrary{decorations.pathmorphing,shapes}
\usetikzlibrary{patterns}
\usepackage[pagebackref, colorlinks=true, linkcolor=blue, citecolor = red]{hyperref}

\usepackage[nameinlink,capitalise,noabbrev]{cleveref}
\crefname{subsection}{Subsection}{Subsection}
\crefname{intuition}{Intuition}{Intuition}

\usepackage[margin=25mm]{geometry}

\renewcommand*{\backref}[1]{}
\renewcommand*{\backrefalt}[4]{({%
		\ifcase #1 Not cited.%
		\or On p.~#2%
		\else On pp.~#2%
		\fi%
	})}

\DeclareMathAlphabet{\mathbbe}{U}{bbold}{m}{n}

\def\DDelta{{\mathbbe{\Delta}}}
\newcommand{\DD}{\DDelta}

\newcommand{\cC}{\mathcal{C}}
\newcommand{\cD}{\mathcal{D}}
\newcommand{\cE}{\mathcal{E}}
\newcommand{\bN}{\mathbb{N}}
\newcommand{\sN}{\mathscr{N}}
\newcommand{\bR}{\mathbb{R}}
\newcommand{\sS}{\mathscr{S}}

\newcommand{\Hom}{\mathrm{Hom}}

\newcommand{\Map}{\mathrm{Map}}
\newcommand{\Mor}{\mathrm{Mor}}
\newcommand{\map}{\mathrm{map}}
\newcommand{\op}{\mathrm{op}}
\newcommand{\Fun}{\mathrm{Fun}}
\newcommand{\Nat}{\mathrm{Nat}}
\newcommand{\Set}{\mathscr{S}\mathrm{et}}
\newcommand{\Yon}{\mathscr{Y}\mathrm{on}}
\newcommand{\sSet}{s\mathscr{S}\mathrm{et}}
\newcommand{\sSp}{s\mathscr{S}}
\newcommand{\Cat}{\mathscr{C}\mathrm{at}}
\newcommand{\Top}{\mathscr{T}\mathrm{op}}
\newcommand{\id}{\mathrm{id}}
\newcommand{\Obj}{\mathrm{Obj}}
\newcommand{\Sing}{\mathrm{Sing}}
\newcommand{\Sp}{\mathrm{Sp}}
\newcommand{\Kan}{\mathscr{K}\mathrm{an}}
\newcommand{\Ho}{\mathrm{Ho}}
\newcommand{\inccat}{\iota_{\Cat}}
\newcommand{\incspace}{\iota_{\mathscr{S}\mathrm{p}}}

\newcommand{\Comp}{\mathrm{Comp}}
\newcommand{\hoequiv}{\mathrm{hoequiv}}
\newcommand{\hoeqchoice}{\mathrm{hoeqchoice}}
\newcommand{\hoeqcomp}{\mathrm{hoeqcomp}}
\newcommand{\forgcomp}{\mathrm{forgcomp}}
\newcommand{\forgchoice}{\mathrm{forgchoice}}
\newcommand{\Path}{\mathscr{P}\mathrm{ath}}

\newcommand{\LInv}{\mathscr{L}\mathrm{Inv}}
\newcommand{\RInv}{\mathscr{R}\mathrm{Inv}}

\newcommand{\nervediagram}{
 \begin{tikzcd}[row sep=1.2cm, column sep=1.2cm]
   \left \{
  \bullet
\right \}
   \arrow[d, shorten >=1ex,shorten <=1ex] \arrow[r, shorten >=1ex,shorten <=1ex]
   \pgfmatrixnextcell 
   \left \{ 
     \begin{tabular}{c}
      $\bullet$
     \end{tabular}
     \begin{tabular}{c}
      $\rightarrow$ 
     \end{tabular}
      \begin{tabular}{c}
      $\bullet$
     \end{tabular}
     \right \}
   \arrow[d, shorten >=1ex,shorten <=1ex]
   \arrow[l, shift left=1.2] \arrow[l, shift right=1.2] 
   \arrow[r, shift right, shorten >=1ex,shorten <=1ex ] \arrow[r, shift left, shorten >=1ex,shorten <=1ex] 
   \pgfmatrixnextcell 
   \left \{ 
     \begin{tabular}{c}
      $\bullet$
     \end{tabular}
     \begin{tabular}{c}
      $\rightarrow$ 
     \end{tabular}
      \begin{tabular}{c}
      $\bullet$
     \end{tabular}
     \begin{tabular}{c}
      $\rightarrow$ 
     \end{tabular}
      \begin{tabular}{c}
      $\bullet$
     \end{tabular}
     \right \}
   \arrow[d, shorten >=1ex,shorten <=1ex]
   \arrow[l] \arrow[l, shift left=2] \arrow[l, shift right=2] 
   \arrow[r, shorten >=1ex,shorten <=1ex] \arrow[r, shift left=2, shorten >=1ex,shorten <=1ex] \arrow[r, shift right=2, shorten >=1ex,shorten <=1ex]
   \pgfmatrixnextcell 
   \cdots 
   \arrow[l, shift right=1] \arrow[l, shift left=1] \arrow[l, shift right=3] \arrow[l, shift left=3] 
   \\ 
   \left \{ \begin{tabular}{c}
      $\bullet$ \\
      $\downarrow$ \\
      $\bullet$
     \end{tabular}
     \hspace{-0.15in}
     \begin{tabular}{c}
      \strut \\
      $\sim$ \\
      \strut
     \end{tabular}
     \hspace{-0.1in}
     \right \}
   \arrow[d, shift right, shorten >=1ex,shorten <=1ex ] \arrow[d, shift left, shorten >=1ex,shorten <=1ex]
   \arrow[u, shift left=1.2] \arrow[u, shift right=1.2] \arrow[r, shorten >=1ex,shorten <=1ex]
   \pgfmatrixnextcell 
   \left \{ \begin{tabular}{c}
      $\bullet$ \\
      $\downarrow$ \\
      $\bullet$
     \end{tabular}
     \hspace{-0.15in}
     \begin{tabular}{c}
      \strut \\
      $\sim$ \\
      \strut
     \end{tabular}
     \hspace{-0.1in}
     \begin{tabular}{c}
      $\rightarrow$ \\
      \\
      $\rightarrow$
     \end{tabular}
     \begin{tabular}{c}
      $\bullet$ \\
      $\downarrow$ \\
      $\bullet$
     \end{tabular}
     \hspace{-0.15in}
     \begin{tabular}{c}
      \strut \\
      $\sim$ \\
      \strut
     \end{tabular}
     \hspace{-0.1in}
     \right \}
   \arrow[d, shift right, shorten >=1ex,shorten <=1ex ] \arrow[d, shift left, shorten >=1ex,shorten <=1ex]
   \arrow[u, shift left=1.2] \arrow[u, shift right=1.2]
   \arrow[l, shift left=1.2] \arrow[l, shift right=1.2] 
   \arrow[r, shift right, shorten >=1ex,shorten <=1ex ] \arrow[r, shift left, shorten >=1ex,shorten <=1ex] 
   \pgfmatrixnextcell  
   \left \{ \begin{tabular}{c}
      $\bullet$ \\
      $\downarrow$ \\
      $\bullet$
     \end{tabular}
     \hspace{-0.15in}
     \begin{tabular}{c}
      \strut \\
      $\sim$ \\
      \strut
     \end{tabular}
     \hspace{-0.1in}
     \begin{tabular}{c}
      $\rightarrow$ \\
      \\
      $\rightarrow$
     \end{tabular}
     \begin{tabular}{c}
      $\bullet$ \\
      $\downarrow$ \\
      $\bullet$
     \end{tabular}
     \hspace{-0.15in}
     \begin{tabular}{c}
      \strut \\
      $\sim$ \\
      \strut
     \end{tabular}
     \hspace{-0.1in}
     \begin{tabular}{c}
      $\rightarrow$ \\
      \\
      $\rightarrow$
     \end{tabular}
     \begin{tabular}{c}
      $\bullet$ \\
      $\downarrow$ \\
      $\bullet$
     \end{tabular}
     \hspace{-0.15in}
     \begin{tabular}{c}
      \strut \\
      $\sim$ \\
      \strut
     \end{tabular}
     \hspace{-0.1in}
     \right \}
   \arrow[d, shift right, shorten >=1ex,shorten <=1ex ] \arrow[d, shift left, shorten >=1ex,shorten <=1ex]
   \arrow[u, shift left=1.2] \arrow[u, shift right=1.2]
   \arrow[l] \arrow[l, shift left=2] \arrow[l, shift right=2] 
   \arrow[r, shorten >=1ex,shorten <=1ex] \arrow[r, shift left=2, shorten >=1ex,shorten <=1ex] \arrow[r, shift right=2, shorten >=1ex,shorten <=1ex]
   \pgfmatrixnextcell 
   \cdots 
   \arrow[l, shift right=1] \arrow[l, shift left=1] \arrow[l, shift right=3] \arrow[l, shift left=3] 
   \\ 
   \left \{ 
   \begin{tabular}{c}
      $\bullet$ \\
      $\downarrow$ \\
      $\bullet$ \\
      $\downarrow$ \\
      $\bullet$
     \end{tabular}
     \hspace{-0.15in}
     \begin{tabular}{c}
      \strut \\
      $\sim$ \\
      \strut \\
      $\sim$ \\
      \strut 
     \end{tabular}
     \hspace{-0.1in}
     \right \}
   \arrow[d, shorten >=1ex,shorten <=1ex] \arrow[d, shift left=2, shorten >=1ex,shorten <=1ex] \arrow[d, shift right=2, shorten >=1ex,shorten <=1ex]
   \arrow[u] \arrow[u, shift left=2] \arrow[u, shift right=2] 
   \arrow[r, shorten >=1ex,shorten <=1ex]
   \pgfmatrixnextcell 
   \left \{ 
   \begin{tabular}{c}
      $\bullet$ \\
      $\downarrow$ \\
      $\bullet$ \\
      $\downarrow$ \\
      $\bullet$
     \end{tabular}
     \hspace{-0.15in}
     \begin{tabular}{c}
      \strut \\
      $\sim$ \\
      \strut \\
      $\sim$ \\
      \strut 
     \end{tabular}
     \hspace{-0.1in}
     \begin{tabular}{c}
      $\rightarrow$ \\
      \\
      $\rightarrow$ \\
      \\
      $\rightarrow$
     \end{tabular}
     \begin{tabular}{c}
      $\bullet$ \\
      $\downarrow$ \\
      $\bullet$ \\
      $\downarrow$ \\
      $\bullet$
     \end{tabular}
     \hspace{-0.15in}
     \begin{tabular}{c}
      \strut \\
      $\sim$ \\
      \strut \\
      $\sim$ \\
      \strut 
     \end{tabular}
     \hspace{-0.1in}
     \right \}
   \arrow[d, shorten >=1ex,shorten <=1ex] \arrow[d, shift left=2, shorten >=1ex,shorten <=1ex] \arrow[d, shift right=2, shorten >=1ex,shorten <=1ex]
   \arrow[u] \arrow[u, shift left=2] \arrow[u, shift right=2]
   \arrow[l, shift left=1.2] \arrow[l, shift right=1.2] 
   \arrow[r, shift right, shorten >=1ex,shorten <=1ex ] \arrow[r, shift left, shorten >=1ex,shorten <=1ex] 
   \pgfmatrixnextcell 
   \left \{ 
   \begin{tabular}{c}
      $\bullet$ \\
      $\downarrow$ \\
      $\bullet$ \\
      $\downarrow$ \\
      $\bullet$
     \end{tabular}
     \hspace{-0.15in}
     \begin{tabular}{c}
      \strut \\
      $\sim$ \\
      \strut \\
      $\sim$ \\
      \strut 
     \end{tabular}
     \hspace{-0.1in}
     \begin{tabular}{c}
      $\rightarrow$ \\
      \\
      $\rightarrow$ \\
      \\
      $\rightarrow$
     \end{tabular}
     \begin{tabular}{c}
      $\bullet$ \\
      $\downarrow$ \\
      $\bullet$ \\
      $\downarrow$ \\
      $\bullet$
     \end{tabular}
     \hspace{-0.15in}
     \begin{tabular}{c}
      \strut \\
      $\sim$ \\
      \strut \\
      $\sim$ \\
      \strut 
     \end{tabular}
     \hspace{-0.1in}
     \begin{tabular}{c}
      $\rightarrow$ \\
      \\
      $\rightarrow$ \\
      \\
      $\rightarrow$
     \end{tabular}
     \begin{tabular}{c}
      $\bullet$ \\
      $\downarrow$ \\
      $\bullet$ \\
      $\downarrow$ \\
      $\bullet$
     \end{tabular}
     \hspace{-0.15in}
     \begin{tabular}{c}
      \strut \\
      $\sim$ \\
      \strut \\
      $\sim$ \\
      \strut 
     \end{tabular}
     \hspace{-0.1in}
     \right \}
   \arrow[d, shorten >=1ex,shorten <=1ex] \arrow[d, shift left=2, shorten >=1ex,shorten <=1ex] \arrow[d, shift right=2, shorten >=1ex,shorten <=1ex]
   \arrow[u] \arrow[u, shift left=2] \arrow[u, shift right=2]
   \arrow[l] \arrow[l, shift left=2] \arrow[l, shift right=2] 
   \arrow[r, shorten >=1ex,shorten <=1ex] \arrow[r, shift left=2, shorten >=1ex,shorten <=1ex] \arrow[r, shift right=2, shorten >=1ex,shorten <=1ex]
   \pgfmatrixnextcell 
   \cdots 
   \arrow[l, shift right=1] \arrow[l, shift left=1] \arrow[l, shift right=3] \arrow[l, shift left=3] 
   \\
   \ \vdots \ 
   \arrow[u, shift right=1] \arrow[u, shift left=1] \arrow[u, shift right=3] \arrow[u, shift left=3]
   \pgfmatrixnextcell \ \vdots \ 
   \arrow[u, shift right=1] \arrow[u, shift left=1] \arrow[u, shift right=3] \arrow[u, shift left=3]
   \pgfmatrixnextcell \ \vdots \ 
   \arrow[u, shift right=1] \arrow[u, shift left=1] \arrow[u, shift right=3] \arrow[u, shift left=3]
   \pgfmatrixnextcell 
 \end{tikzcd}
}

\newcommand{\simplicialspacediagram}{
   \begin{tikzcd}[row sep=1.2cm, column sep=1.2cm]
   X_{00} 
   \arrow[d, shorten >=1ex,shorten <=1ex] \arrow[r, shorten >=1ex,shorten <=1ex]
   \pgfmatrixnextcell X_{10} 
   \arrow[d, shorten >=1ex,shorten <=1ex]
   \arrow[l, shift left=1.2] \arrow[l, shift right=1.2] 
   \arrow[r, shift right, shorten >=1ex,shorten <=1ex ] \arrow[r, shift left, shorten >=1ex,shorten <=1ex] 
   \pgfmatrixnextcell X_{20} 
   \arrow[d, shorten >=1ex,shorten <=1ex]
   \arrow[l] \arrow[l, shift left=2] \arrow[l, shift right=2] 
   \arrow[r, shorten >=1ex,shorten <=1ex] \arrow[r, shift left=2, shorten >=1ex,shorten <=1ex] \arrow[r, shift right=2, shorten >=1ex,shorten <=1ex]
   \pgfmatrixnextcell \cdots 
   \arrow[l, shift right=1] \arrow[l, shift left=1] \arrow[l, shift right=3] \arrow[l, shift left=3] 
   \\
   X_{01} 
   \arrow[d, shift right, shorten >=1ex,shorten <=1ex ] \arrow[d, shift left, shorten >=1ex,shorten <=1ex]
   \arrow[u, shift left=1.2] \arrow[u, shift right=1.2] \arrow[r, shorten >=1ex,shorten <=1ex]
   \pgfmatrixnextcell X_{11} 
   \arrow[d, shift right, shorten >=1ex,shorten <=1ex ] \arrow[d, shift left, shorten >=1ex,shorten <=1ex]
   \arrow[u, shift left=1.2] \arrow[u, shift right=1.2]
   \arrow[l, shift left=1.2] \arrow[l, shift right=1.2] 
   \arrow[r, shift right, shorten >=1ex,shorten <=1ex ] \arrow[r, shift left, shorten >=1ex,shorten <=1ex] 
   \pgfmatrixnextcell X_{21} 
   \arrow[d, shift right, shorten >=1ex,shorten <=1ex ] \arrow[d, shift left, shorten >=1ex,shorten <=1ex]
   \arrow[u, shift left=1.2] \arrow[u, shift right=1.2]
   \arrow[l] \arrow[l, shift left=2] \arrow[l, shift right=2] 
   \arrow[r, shorten >=1ex,shorten <=1ex] \arrow[r, shift left=2, shorten >=1ex,shorten <=1ex] \arrow[r, shift right=2, shorten >=1ex,shorten <=1ex]
   \pgfmatrixnextcell \cdots 
   \arrow[l, shift right=1] \arrow[l, shift left=1] \arrow[l, shift right=3] \arrow[l, shift left=3] 
   \\
   X_{02} 
   \arrow[d, shorten >=1ex,shorten <=1ex] \arrow[d, shift left=2, shorten >=1ex,shorten <=1ex] \arrow[d, shift right=2, shorten >=1ex,shorten <=1ex]
   \arrow[u] \arrow[u, shift left=2] \arrow[u, shift right=2] 
   \arrow[r, shorten >=1ex,shorten <=1ex]
   \pgfmatrixnextcell X_{12} 
   \arrow[d, shorten >=1ex,shorten <=1ex] \arrow[d, shift left=2, shorten >=1ex,shorten <=1ex] \arrow[d, shift right=2, shorten >=1ex,shorten <=1ex]
   \arrow[u] \arrow[u, shift left=2] \arrow[u, shift right=2]
   \arrow[l, shift left=1.2] \arrow[l, shift right=1.2] 
   \arrow[r, shift right, shorten >=1ex,shorten <=1ex ] \arrow[r, shift left, shorten >=1ex,shorten <=1ex] 
   \pgfmatrixnextcell X_{22} 
   \arrow[d, shorten >=1ex,shorten <=1ex] \arrow[d, shift left=2, shorten >=1ex,shorten <=1ex] \arrow[d, shift right=2, shorten >=1ex,shorten <=1ex]
   \arrow[u] \arrow[u, shift left=2] \arrow[u, shift right=2]
   \arrow[l] \arrow[l, shift left=2] \arrow[l, shift right=2] 
   \arrow[r, shorten >=1ex,shorten <=1ex] \arrow[r, shift left=2, shorten >=1ex,shorten <=1ex] \arrow[r, shift right=2, shorten >=1ex,shorten <=1ex]
   \pgfmatrixnextcell \cdots 
   \arrow[l, shift right=1] \arrow[l, shift left=1] \arrow[l, shift right=3] \arrow[l, shift left=3] 
   \\
   \ \vdots \ 
   \arrow[d, leftrightsquigarrow]
   \arrow[u, shift right=1] \arrow[u, shift left=1] \arrow[u, shift right=3] \arrow[u, shift left=3]
   \pgfmatrixnextcell \ \vdots \ 
   \arrow[d, leftrightsquigarrow]
   \arrow[u, shift right=1] \arrow[u, shift left=1] \arrow[u, shift right=3] \arrow[u, shift left=3]
   \pgfmatrixnextcell \ \vdots \ 
   \arrow[d, leftrightsquigarrow]
   \arrow[u, shift right=1] \arrow[u, shift left=1] \arrow[u, shift right=3] \arrow[u, shift left=3]
   \pgfmatrixnextcell \ \ddots \ 
   \arrow[d, leftrightsquigarrow]
   \\
    X_{0\bullet} \arrow[r, shorten >=1ex,shorten <=1ex]
   \pgfmatrixnextcell X_{1\bullet} 
   \arrow[l, shift left=1.2] \arrow[l, shift right=1.2] 
   \arrow[r, shift right, shorten >=1ex,shorten <=1ex ] \arrow[r, shift left, shorten >=1ex,shorten <=1ex] 
   \pgfmatrixnextcell X_{2\bullet} 
   \arrow[l] \arrow[l, shift left=2] \arrow[l, shift right=2] 
   \arrow[r, shorten >=1ex,shorten <=1ex] \arrow[r, shift left=2, shorten >=1ex,shorten <=1ex] \arrow[r, shift right=2, shorten >=1ex,shorten <=1ex]
   \pgfmatrixnextcell \cdots  
   \arrow[l, shift right=1] \arrow[l, shift left=1] \arrow[l, shift right=3] \arrow[l, shift left=3] 
 \end{tikzcd}
}

\newtheorem{theorem}[equation]{Theorem}
\newtheorem{lemma}[equation]{Lemma}
\newtheorem{proposition}[equation]{Proposition}
\newtheorem{corollary}[equation]{Corollary}

\theoremstyle{definition}
\newtheorem{definition}[equation]{Definition}
\newtheorem{example}[equation]{Example}
\newtheorem{notation}[equation]{Notation}
\newtheorem{construction}[equation]{Construction}

\theoremstyle{remark}
\newtheorem{remark}[equation]{Remark}
\newtheorem{intuition}[equation]{Intuition}

\numberwithin{equation}{section}

\title{Introduction to Complete Segal Spaces}

\author{Nima Rasekh}

\date{June 2026}



\begin{document}

\begin{abstract}
  We introduce $\infty$-categories via complete Segal spaces. We primarily focus on foundational concepts, aiming to provide proper motivation and intuition, requiring only a rudimentary background in category theory. 
\end{abstract}

\maketitle

\tableofcontents

\section{Teaching the Idea of \texorpdfstring{$\infty$}{oo}-Categories via Complete Segal Spaces}\label{sec:intro}
 
\subsection{\texorpdfstring{$\infty$}{oo}-Category Theory: From Papers to Textbooks} 
$\infty$-categories have now become a standard aspect of many areas of mathematics. This started with algebraic topology and homotopy theory. It then generalized to various aspects of (differential or algebraic) geometry, such as derived geometry or spectral geometry. It has also become a standard tool in mathematical physics, and particularly topological field theories, and finally has found its way into the foundations of mathematics, in particular homotopy type theory.

These developments created the need for resources introducing the theory of $\infty$-categories to a wider audience. Unfortunately, in the early days of $\infty$-category theory\footnote{When the author was first learning about these concepts.}, the main sources were the original research papers, several books that combined exposition and research\footnote{Around that time, \cite{lurie2009htt} was probably the most standard introduction to $\infty$-category theory.}, and a few surveys.

Fortunately, in recent years the situation has significantly improved. There is now a wide range of introductions to $\infty$-categories from different angles, as we further discuss at the end (\cref{sec:next}). As one might anticipate, these sources assume reasonable familiarity with category theory, sometimes even $2$-category theory, as well as homotopy theory, and introduce that theory on that basis. Mathematically, this is of course a reasonable and correct approach. Educationally, it can create a challenge conveying intuitions regarding higher categorical notions without delving into too many technical nuances.

\subsection{The Magic of Simplicial Sets} \label{subsec:magic of ss}
One interesting aspect of mathematics education in higher category theory and homotopy theory is the central role of \emph{simplicial sets} in $\infty$-category theory. A lot of modern higher category theory implicitly or explicitly relies on the simplicial machinery. Yet, interestingly enough, most first introductions to algebraic topology do not feature any introduction to simplicial sets \cite{rotman1988algebraictopology,massey1991basicalgebraictopology,munkreskrantzparks2025elementsofalgebraictopology,tomdieck2008algebraictopology}, or include very little \cite{hatcher2002at,may1999algebraictopology}. Similarly, most introductory texts in category theory do not include a comprehensive introduction to simplicial sets \cite{awodey2010categorytheory,leinster2014basiccategorytheory,riehl2016context}, with some exceptions \cite{maclane1998categories}.

It is hence unsurprising that there are specifically-designed introductions to algebraic topology focused on simplicial sets \cite{may1992simplicialobjects,goerssjardine2009simplicialhomotopytheory}. Moreover, most introductions to $\infty$-categories also include some review of simplicial sets. What appears to be missing is a non-technical introduction to simplicial and higher categorical ideas that is accessible to a reader with some background in algebraic topology and category theory.

\subsection{Who is this for?}
This note precisely aims to address this issue. It is designed as a first introduction to the \emph{definition}, \emph{idea} and \emph{intuition} of an $\infty$-category via one particular approach to $\infty$-categories, namely complete Segal spaces. It is aimed at a reader with some rudimentary familiarity with classical category theory and algebraic topology. As a result, this note includes an {\color{blue} Intuition} environment that aims to convey additional intuition to the reader. Ideally, having read this note, the interested reader is motivated and willing to tackle the technical intricacies needed to understand the theory via the many excellent sources suggested in \cref{sec:next}.

This means that this note is \textbf{not} meant as a comprehensive introduction to \emph{$\infty$-category theory}. Only particular aspects of the theory have been covered in a way to pursue the intended mission, namely conveying essential ideas. This in particular means that all proofs in this work are laid out with details and are self-contained. As part of that philosophy, results whose proofs required advanced methods or significant mathematical background, have been stated without proof. In those cases references have been provided for the interested reader.

\subsection{Why Complete Segal Spaces?} \label{subsec:why}
Most modern approaches to $\infty$-category theory use some form of ``model-independence''\footnote{Here the word model should not be interpreted as models of a theory, in the sense of mathematical logic, but rather as concrete implementations of an abstract idea of an $\infty$-category.}, which has rightfully established itself as the most sensible approach to $\infty$-category theory. See \cref{sec:next} for a range of such approaches. However, educationally, it can be instructive to take first steps via a concrete model. This can help connect the axiomatic approaches to concrete examples.

Taking for granted that we want to work with models, it is natural to wonder why we have chosen complete Segal spaces. In particular, some readers might already be familiar with other popular models, such as \emph{quasi-categories} \cite{joyal2008notes,joyal2008theory,lurie2009htt}. Let us hence review some reasons for this choice.
\begin{enumerate}[leftmargin=*]
  \item \textbf{Completeness:} One fascinating aspect that distinguishes categories from $\infty$-categories is the notion of \emph{completeness}. This has become more evident given the importance of the related notion of \emph{univalence} in homotopy type theory \cite{hottbook2013}. Complete Segal spaces are the most natural $\infty$-categorical setting to witness the importance of the completeness condition, as we analyze in \cref{subsec:why css}.
  \item \textbf{Internal $\infty$-Category Theory:} In the same way that categories admit internalization via \emph{internal categories}, there are analogous notions of \emph{internal $\infty$-categories} \cite{rasekh2022cartesian}, with various geometric applications. These internal notions are most naturally expressed in the language of complete Segal spaces.
  \item \textbf{Formalization via Proof Assistants:} Formalization of mathematics enters mathematical statements into computers via suitable programming languages, such as Lean or Rocq, that are then formally verified. As of now the most successful effort formalizing $\infty$-category via proof assistants utilizes a foundation motivated by complete Segal spaces, known as \emph{simplicial homotopy type theory} \cite{riehlshulman2017rezktypes}. This has resulted in a specialized proof assistant, \emph{Rzk}, and formalizations of various $\infty$-categorical results, such as the \emph{$\infty$-categorical Yoneda lemma} \cite{kudasovriehlweinberger2024rzkyoneda}. Understanding the theory of complete Segal spaces is hence crucial for advancing the formalization of $\infty$-category theory via the proof assistant Rzk.
  \item \textbf{Technical Relevance:} While the ``model-independent'' methods have proven to be powerful, there are still results regarding $\infty$-categories that can only be understood via complete Segal spaces. One such example is computing the \emph{filtered colimits} of $\infty$-categories, which is the crucial insight needed to study \emph{filter quotient $\infty$-categories} \cite{rasekh2021filterquotient}.
  \item \textbf{Examples:} Various explicit examples of $\infty$-categories are obtained by constructing explicit complete Segal spaces \cite{lurie2009cobordism,calaquescheimbauer2019cobordism,kock2023wholegrainpetri,mukherjeerasekh2022twisted,rasekh2024model}.
  \item \textbf{$(\infty,n)$-Category Theory:} $\infty$-category theory admits a generalization to a notion of $(\infty,n)$-categories. While there are model-independent approaches to $(\infty,n)$-categories \cite{barwickschommerpries2021unicity}, concrete models still play a far more prominent role in that setting. Many prominent models of $(\infty,n)$-categories, such as \emph{$\Theta_n$-spaces} \cite{rezk2010thetanspaces} or \emph{$n$-fold complete Segal spaces} \cite{barwick2005nfoldsegalspaces}, are direct generalizations of complete Segal spaces. 
\end{enumerate}

\subsection{Background}
As explained, this text is meant to be mostly self-contained. We take for granted basic familiarity with point-set topology \cite[Sections 1--7, 9, 12--19, 23]{munkres2000topology}. It is helpful to have some familiarity with certain aspects of algebraic topology, such as homotopy groups \cite[Section 1.1, Section 4.1]{hatcher2002at}. Moreover, the ideal reader already has some familiarity with category theory, such as \cite[Section 1--3]{maclane1998categories} or \cite[Chapter 1--3]{riehl2016context}, however, all relevant definitions have also been reviewed in \cref{sec:ct}.

\subsection{Notation}
For the reader who has some familiarity with higher categorical literature, some quick comments regarding the notational conventions:
\begin{itemize}[leftmargin=*]
  \item We mostly use the notation introduced by Rezk \cite{rezk2001css}. This in particular includes the notation $F(n)$, $\Delta[l]$ for the generating simplicial spaces (\cref{not:FDelta}).\footnote{More recent papers might denote $F(n)$ as $\Delta[n,0]$ and $\Delta[l]$ as $\Delta[0,l]$.}
  \item In some cases we deviate from Rezk's notation. Concretely, we use $\sSet$ for the category of simplicial sets (\cref{def:sset}) instead of $\sS$. Moreover, we use the (now more common) terminology Rezk nerve $\sN$ (\cref{def:rezk nerve}) instead of \emph{classification diagram} (\cref{rem:classification diagram}).
\end{itemize}

\subsection{Acknowledgements} 
I want to thank my advisor Charles Rezk and PhD committee member Matt Ando for fruitful conversations and suggestions. I thank Matt Feller and Udit Mavinkurve for feedback on an earlier version of this text. Moreover, I thank Matthias Frerichs, Hannes Berkenhagen, and Alessandro Nanto for their detailed feedback, which has significantly improved the exposition.

\section{Simplicial Sets} \label{sec:ss}
The goal of this section is to define the magical simplex category $\DD$ (\cref{subsec:magic of ss}) and the category of simplicial sets, $\sSet$. In the next section we will then see how simplicial sets relate to categories and topological spaces. For more extensive introductions, assuming more category theory, see \cite{riehl2011leisurely,friedman2012elementarysimplicialsets} or \cite[Section I.1]{goerssjardine2009simplicialhomotopytheory}.

\begin{definition}[Simplex Category] \label{def:delta}
  Let $\DD$ be the category defined as follows:
  \begin{itemize}
    \item $\Obj_{\DD} \coloneq \{ [n] \mid n \geq 0 \}$, where $[n] \coloneq \{0 \leq \cdots \leq n\}$ is defined in \cref{ex:ncat}.
    \item For $[m],[n]$ in $\Obj_{\DD}$, $\Hom_{\DD}([m],[n])$ is the set of all functors from $[m]$ to $[n]$ (\cref{ex:ncat}).
    \item The composition and identity are defined as in $\Cat$ (\cref{ex:cat}).
  \end{itemize} 
\end{definition}

\begin{notation} \label{not:delta morphism}
 By definition a morphism $[m] \to [n]$ in $\DD$ is precisely a non-decreasing chain $0 \leq a_0 \leq a_1 \leq \dots \leq a_m \leq n$, where $a_i$ denotes the image of $i \in [m]$. We will hence denote such a morphism by $a_0 a_1 ... a_m$.
\end{notation}

In general $\DD$ has too many morphisms. Indeed, using a basic counting argument from combinatorics we see that $\Hom_{\DD}([m],[n])$ has $\binom{m+n+1}{n}$ elements. However, we can reduce arbitrary morphisms to a smaller set of ``generating morphisms'' via composition.

\begin{notation} \label{not:di si}
  Some morphisms in $\DD$ merit their own notation. 
 \begin{itemize}[leftmargin=*]
  \item For each $n \geq 0$ and $0 \leq i \leq n+1$ there is a unique injective map 
  \[d^i\colon[n] \to [n+1]\]
  whose image does not contain $i$.
  \item For each $n \geq 1$ and $0 \leq i \leq n-1$ there is a unique surjective map 
  \[s^i\colon[n] \to [n-1]\]
  whose pre-image over $i$ is $\{i,i+1\}$ and other pre-images have one element. 
 \end{itemize}
\end{notation}

Why do we specifically single out these morphisms? It turns out these morphisms are enough to get everything else back. The proof is a somewhat tedious combinatorial argument, and hence we refer the reader to \cite[Section VII.5]{maclane1998categories}.
\begin{lemma} \label{lemma:delta unique}
  An arbitrary morphism $\delta\colon [m] \to [n]$ in $\DD$ can be uniquely written as a composition of the form 
  \[\delta = d^{i_1}d^{i_2}\cdots d^{i_k}s^{j_1}\cdots s^{j_l},\]
  where $0 \leq i_k < \dots < i_2 < i_1 \leq n$ and $0 \leq j_1 < j_2 < \dots < j_l < m$.
\end{lemma}

This means we only need to better understand how the morphisms $d^i$ and $s^i$ interact with each other. Here we can use the explicit description given in \cref{not:di si}, to obtain the \emph{cosimplicial identities}.

\begin{lemma}[Cosimplicial identities] \label{lemma:cosimplicial identities}
 The morphisms $d^i$ and $s^i$ satisfy the following identities:
\begin{enumerate}
  \item $d^j d^i = d^i d^{j-1}$ for $i < j$.
  \item $s^j s^i = s^i s^{j+1}$ for $i \leq j$.
  \item $s^j d^i = \begin{cases}
    d^i s^{j-1} & i < j \\
    \id & i = j, j+1 \\
    d^{i-1} s^j & i > j + 1
  \end{cases}$
\end{enumerate}
\end{lemma}

\begin{remark} \label{rem:delta}
  Combining \cref{lemma:delta unique,lemma:cosimplicial identities} the category $\DD$ can be depicted as follows:
\[
    \begin{tikzcd}[row sep=0.5in, column sep=0.5in]
      [0]
      \arrow[r, shift left=1.4, "d^0"] 
      \arrow[r, shift right=1.4, "d^1"'] 
      & 
      {[1]}
      \arrow[l, shorten >=1ex,shorten <=1ex, "s^0" yshift= .2cm]
      \arrow[r] \arrow[r, shift left=2, "d^0"] \arrow[r, shift right=2, "d^2"']   
      & {[2]}
      \arrow[l, shift right, shorten >=1ex,shorten <=1ex ] 
      \arrow[l, shift left, shorten >=1ex,shorten <=1ex] 
      \arrow[r, shift right=1] \arrow[r, shift left=1] \arrow[r, shift right=3] \arrow[r, shift left=3] 
      & 
      \cdots 
      \arrow[l, shorten >=1ex,shorten <=1ex] \arrow[l, shift left=2, shorten >=1ex,shorten <=1ex] 
      \arrow[l, shift right=2, shorten >=1ex,shorten <=1ex]
    \end{tikzcd}.
  \]
\end{remark}

We can now use the category $\DD$ to define simplicial sets. Here we use the category of sets (\cref{ex:set}), and functor categories (\cref{def:functor category}).

\begin{definition}[Simplicial Object]
  Let $\cC$ be a category. A \emph{simplicial object} in $\cC$ is a functor $X\colon \DD^{\op} \to \cC$. 
\end{definition}

\begin{definition}[Simplicial Set] \label{def:sset}
  A \emph{simplicial set} is a simplicial object in $\Set$. The category of simplicial sets, denoted $\sSet$, is the functor category $\Fun(\DD^{\op},\Set)$.
\end{definition}

Following \cref{lemma:delta unique,lemma:cosimplicial identities}, we can now give a more explicit description of what a simplicial set is.

\begin{lemma} \label{lemma:ss diagram}
 Let $\cC$ be a category. A simplicial object $X$ is uniquely determined by the following data:
 \begin{itemize}[leftmargin=*]
  \item A choice of objects $X_0 = X([0]), X_1 = X([1]),  ...$.
  \item For each $n \geq 0$ and $0 \leq i \leq n+1$ morphisms $X(d^i) = d_i\colon X_{n+1} \to X_n$.
  \item For each $n \geq 1$ and $0 \leq i \leq n-1$ there are morphisms $X(s^i) = s_i\colon X_{n-1} \to X_n$.
  \item Satisfying the \emph{simplicial identities}:
  \begin{itemize}[leftmargin=*]
    \item $d_i d_j = d_{j-1} d_i $ for $i < j$.
    \item $s_i s_j = s_{j+1} s_i $ for $i \leq j$.
    \item $d_i s_j = \begin{cases}
    s_{j-1} d_i  & i < j \\
    \id & i = j, j+1 \\
     s_j d_{i-1} & i > j + 1
  \end{cases}$
  \end{itemize}
 \end{itemize}
\end{lemma}

\begin{remark} \label{rem:ss diagram}
  Following this lemma, similar to \cref{rem:delta}, we can depict a simplicial set $X$ as follows:
  \[
  \begin{tikzcd}[column sep=1.2cm]
   X_0 \arrow[r, shorten >=1ex,shorten <=1ex, "s_0"'  yshift= .2cm]
   & X_1 
   \arrow[l, shift left=1.4, "d_1"] \arrow[l, shift right=1.4, "d_0"'] 
   \arrow[r, shift right, shorten >=1ex,shorten <=1ex ] \arrow[r, shift left, shorten >=1ex,shorten <=1ex] 
   & X_2 
   \arrow[l] \arrow[l, shift left=2, "d_2"] \arrow[l, shift right=2, "d_0"'] 
   \arrow[r, shorten >=1ex,shorten <=1ex] \arrow[r, shift left=2, shorten >=1ex,shorten <=1ex] 
   \arrow[r, shift right=2, shorten >=1ex,shorten <=1ex]
   & \cdots 
   \arrow[l, shift right=1] \arrow[l, shift left=1] \arrow[l, shift right=3] \arrow[l, shift left=3] 
 \end{tikzcd}.
  \]
  Notice all arrows are reversed as this functor is mapping out of $\DD^{\op}$.
\end{remark}

\begin{definition}[Degenerate Simplex]
 Let $X$ be a simplicial set. An $n$-simplex $\sigma$ in $X_n$ is called \emph{degenerate} if there exists a simplex $\tau$ in $X_{n-1}$ and $0 \leq i \leq n-1$ such that $\sigma = s_i(\tau)$. Otherwise, $\sigma$ is called \emph{non-degenerate}.  
\end{definition}

\begin{intuition}
  Geometrically, we imagine a simplicial set $X$ as a recipe for constructing a topological space with a notion of ``direction''. Here $X_0$ is the set of vertices, non-degenerate elements in $X_1$ are the set of edges (which is where the direction comes from), non-degenerate elements in $X_2$ are the set of triangles, and so on. The maps $d_0,d_1\colon X_1 \to X_0$ specify the source and target of each edge, while the maps $d_0,d_1,d_2\colon X_2 \to X_1$ specify the edges that form the boundary of each triangle. The degenerate simplices do not manifest geometrically, but do play an important role (see \cref{rem:degenerate in products}).
\end{intuition}

Let us see one simple example of a simplicial set.

\begin{definition}[Constant simplicial set] \label{def:constant simplicial set}
  A simplicial set $X$ is \emph{constant} if all $s_i, d_i$ are bijections.
\end{definition}

Every set $S$ gives rise to a constant simplicial set, as we see in the next example.

\begin{example} \label{ex:constant simplicial set}
  Let $S$ be a set. Define the constant simplicial set $S\colon \DD^{\op} \to \Set$ by $S_n = S$ for all $n \geq 0$, and all face and degeneracy maps are identity maps.
\end{example}

Let us now look at representable functors in the context of simplicial sets. 

\begin{definition}[Standard simplex] \label{def:standard simplex}
  For a given $n \geq 0$, the \emph{standard $n$-simplex} is the representable functor $\Hom_{\DD}(-,[n])$ (\cref{ex:repcontra}), denoted by $\Delta[n]$.
\end{definition}

\begin{remark} \label{rem:delta morphism}
  In \cref{not:delta morphism} we explained how we can depict elements in $\Hom_{\DD}([m],[n]) = \Delta[n]_m$ as a non-decreasing chain $a_0... a_m$. From this perspective, the map $d_i\colon \Delta[n]_{m+1} \to \Delta[n]_m$ is given by deleting the $i$-th element of the chain, while the map $s_i\colon \Delta[n]_{m-1} \to \Delta[n]_m$ is given by repeating the $i$-th element of the chain.
\end{remark}

\begin{remark} \label{rem:non-degenerate}
  As the degeneracy maps $s_i$ add repeating elements (\cref{rem:delta morphism}), an $m$-simplex $[m] \to [n]$ in $\Delta[n]_m$ is non-degenerate if and only if the corresponding sequence is injective, meaning the chain $a_0... a_m$ has no repeated elements.
\end{remark}

\begin{example}[0-Simplex]
  The standard $0$-simplex $\Delta[0]$ is the simplicial set defined by \[\Delta[0]_n \coloneq \Hom_{\DD}([n],[0]) = \{0...0\},\] where in the last step we used \cref{not:delta morphism}. This means we can depict $\Delta[0]$ as follows:
\[
  \begin{tikzcd}[column sep=1.2cm]
   \{0\} \arrow[r, shorten >=1ex,shorten <=1ex, "s_0"'  yshift= .2cm]
   & \{00\} 
   \arrow[l, shift left=1.4, "d_1" near end, "\cong" near start] \arrow[l, shift right=1.4, "d_0"' near end, "\cong"' near start] 
   \arrow[r, shift right, shorten >=1ex,shorten <=1ex ] \arrow[r, shift left, shorten >=1ex,shorten <=1ex] 
   & \{000\} 
   \arrow[l] \arrow[l, shift left=2, "d_2" near end, "\cong" near start] \arrow[l, shift right=2, "d_0"' near end, "\cong"' near start] 
   \arrow[r, shorten >=1ex,shorten <=1ex] \arrow[r, shift left=2, shorten >=1ex,shorten <=1ex] 
   \arrow[r, shift right=2, shorten >=1ex,shorten <=1ex]
   & \cdots 
   \arrow[l, shift right=1] \arrow[l, shift left=1] \arrow[l, shift right=3] \arrow[l, shift left=3] 
 \end{tikzcd}.
  \]
\end{example}

\begin{lemma} \label{lemma:delta 0}
  For every simplicial set $X$, there is a unique map $X \to \Delta[0]$.
\end{lemma}
\begin{proof}
  We can directly see that the unique map on sets $X_n \to \{00...0\}$ is natural, as all $d_i$ and $s_i$ in $\Delta[0]$ are bijections.
\end{proof}

\begin{intuition}
  $\Delta[0]$ has a single non-degenerate simplex, namely $0$ in $\Delta[0]_0$. Hence, geometrically, $\Delta[0]$ is a single point.
  \[ 
  \begin{tikzcd} 
    \overset{0}{\bullet}
  \end{tikzcd}
  \]
\end{intuition}

\begin{example}[1-Simplex]
  The standard $1$-simplex $\Delta[1]$ is the simplicial set defined by 
  \[\Delta[1]_n \coloneq \Hom_{\DD}([n],[1]) = \{0...0, 0...01, 0...011,..., 1...1\},\] 
  where we again used \cref{not:delta morphism}. This means we can depict $\Delta[1]$ as follows:
  \[
  \begin{tikzcd}[column sep=1.2cm]
   \{0,1\} \arrow[r, shorten >=1ex,shorten <=1ex, "s_0"'  yshift= .2cm]
   & \{00,01,11\} 
   \arrow[l, shift left=1.4, "d_1"] \arrow[l, shift right=1.4, "d_0"'] 
   \arrow[r, shift right, shorten >=1ex,shorten <=1ex ] \arrow[r, shift left, shorten >=1ex,shorten <=1ex] 
   & \{000,001,011,111\} 
   \arrow[l] \arrow[l, shift left=2, "d_2" near end, "\cong" near start] \arrow[l, shift right=2, "d_0"' near end, "\cong"' near start] 
   \arrow[r, shorten >=1ex,shorten <=1ex] \arrow[r, shift left=2, shorten >=1ex,shorten <=1ex] 
   \arrow[r, shift right=2, shorten >=1ex,shorten <=1ex]
   & \cdots 
   \arrow[l, shift right=1] \arrow[l, shift left=1] \arrow[l, shift right=3] \arrow[l, shift left=3] 
 \end{tikzcd}.
  \]
\end{example}

\begin{intuition}
  $\Delta[1]$ has three non-degenerate simplices, namely $0$ and $1$ in $\Delta[1]_0$ and $01$ in $\Delta[1]_1$. Hence, geometrically, $\Delta[1]$ is a single edge connecting two vertices.
  \[ 
  \begin{tikzcd} 
    \overset{0}{\bullet} \arrow[r, "01"] & \overset{1}{\bullet}
  \end{tikzcd}
  \]
\end{intuition}

\begin{example}[2-Simplex]
  The standard $2$-simplex $\Delta[2]$ is the simplicial set defined by 
  \[\Delta[2]_n \coloneq \Hom_{\DD}([n],[2]) = \{0...0, 0...01, 0...011,..., 2...2\},\] 
  where we again used \cref{not:delta morphism}. This means we can depict $\Delta[2]$ as follows:
  {\small
  \[
  \begin{tikzcd}[column sep=0.8cm]
   \{0,1,2\} \arrow[r, shorten >=1ex,shorten <=1ex, "s_0"'  yshift= .2cm]
   & \{00,01,11,12,22,02\} 
   \arrow[l, shift left=1.4, "d_1"] \arrow[l, shift right=1.4, "d_0"'] 
   \arrow[r, shift right, shorten >=1ex,shorten <=1ex ] \arrow[r, shift left, shorten >=1ex,shorten <=1ex] 
   & \{000,001,011,111,112,122,222,002,022,012\} 
   \arrow[l] \arrow[l, shift left=2, "d_2"] \arrow[l, shift right=2, "d_0"'] 
   \arrow[r, shorten >=1ex,shorten <=1ex] \arrow[r, shift left=2, shorten >=1ex,shorten <=1ex] 
   \arrow[r, shift right=2, shorten >=1ex,shorten <=1ex]
   & \cdots 
   \arrow[l, shift right=1] \arrow[l, shift left=1] \arrow[l, shift right=3] \arrow[l, shift left=3] 
 \end{tikzcd}.
  \]
  }
\end{example}

\begin{intuition}
  $\Delta[2]$ has seven non-degenerate simplices, namely $0$, $1$, and $2$ in $\Delta[2]_0$, and $01$, $12$, and $02$ in $\Delta[2]_1$, and $012$ in $\Delta[2]_2$. Hence, geometrically, $\Delta[2]$ is a single triangle with three edges and three vertices.
  \[ 
  \begin{tikzcd}[row sep=0.5cm, column sep=0.5cm]
    & \overset{1}{\bullet} \arrow[dr, "12"] & \\
    \overset{0}{\bullet} \arrow[ur, "01"] \arrow[rr, ""{name=U, above}, "02"'] &   & \overset{2}{\bullet}
    \arrow[start anchor={[xshift=-5ex, yshift=4ex]}, to=U, phantom, "012"]
  \end{tikzcd}
  \]
\end{intuition}
     
Having learned about the standard simplices, we now continue our theoretical analysis of simplicial sets, relying on the Yoneda lemma (\cref{thm:yoneda}).

\begin{lemma}[Yoneda lemma for simplicial sets] \label{lemma:delta yoneda}
  For a given $n \geq 0$ and a simplicial set $X$, there is a natural isomorphism $\Hom_{\sSet}(\Delta[n],X) \cong X_n$.
\end{lemma}

We can in particular use the Yoneda lemma to understand the morphisms between standard simplices.

\begin{corollary}
 Let $n, k \geq 0$. There is a natural isomorphism $\Hom_{\sSet}(\Delta[n],\Delta[k]) \cong \Hom_{\DD}([n],[k])$.
\end{corollary}

Later on, beyond the standard simplices, we also need various sub-objects of $\Delta[n]$.

\begin{definition}[Sub-simplicial set] \label{def:sub simplicial set}
  Let $X$ be a simplicial set. A \emph{sub-simplicial set} $Y$ of $X$ is a natural monomorphism $Y \to X$.
\end{definition}

\begin{remark}
  More explicitly, for a simplicial set $X$, a sub-simplicial set $Y$ of $X$ consists of subsets $Y_n \subseteq X_n$ for every $n \geq 0$ such that for every $\delta\colon [m] \to [n]$ in $\DD$, the function $X(\delta)\colon X_n \to X_m$ restricts to a function $Y(\delta)\colon Y_n \to Y_m$.
\end{remark}

One way to create sub-simplicial sets is by picking a subset of the $0$-simplices.

\begin{example} \label{ex:sub simplicial set}
  Let $X$ be a simplicial set and let $S \subseteq X_0$ be a subset of the $0$-simplices. Then the \emph{sub-simplicial set generated by $S$} is the sub-simplicial set $\bar{S}$ in $X$, with $\bar{S}_n \subseteq X_n$ for every $n \geq 0$, given by 
  \[ \bar{S}_n = \{ \sigma \in X_n \mid \forall 0 \leq i \leq n (i^*\sigma \in S)\}.\]
  Note, in particular $\bar{S}_0 = S$. In other words, $\bar{S}$ is the largest sub-simplicial set of $X$ whose $0$-simplices are precisely $S$.  
\end{example}

There are, however, several specific sub-simplicial sets of interest, that we now review.

\begin{definition}[Boundary] \label{def:boundary}
  For every $n \geq 0$, the \emph{boundary} of $\Delta[n]$, denoted $\partial \Delta[n]$, is the sub-simplicial set of $\Delta[n]$ with $\partial \Delta[n]_m \coloneq \{ \delta\colon [m] \to [n] \mid \mathrm{Im}(\delta) \neq [n] \}$. This is indeed a sub-simplicial set as non-surjective maps are closed under precomposition.
\end{definition}

\begin{example}[0-Boundary]
 Every single simplex of $\Delta[0]_n$ is given by a surjective map $[n] \to [0]$. Hence, $\partial \Delta[0]$ is the empty simplicial set, meaning $\partial \Delta[0] = \emptyset$.
\end{example}

\begin{example}[1-Boundary]
  Following \cref{rem:non-degenerate}, $\Delta[1]$ has three non-degenerate simplices. Among those, $0,1$ in $\Delta[1]_0$ are not surjective, while $01$ in $\Delta[1]_1$ is surjective (following our notational convention it corresponds to the identity $[1] \to [1]$). Hence, we can depict $\partial \Delta[1]$ as the simplicial set: 
  \[
  \begin{tikzcd}[column sep=1.2cm]
   \{0,1\} \arrow[r, shorten >=1ex,shorten <=1ex, "s_0"'  yshift= .2cm]
   & \{00,11\} 
   \arrow[l, shift left=1.4, "d_1"] \arrow[l, shift right=1.4, "d_0"'] 
   \arrow[r, shift right, shorten >=1ex,shorten <=1ex ] \arrow[r, shift left, shorten >=1ex,shorten <=1ex] 
   & \{000,111\} 
   \arrow[l] \arrow[l, shift left=2, "d_2"] \arrow[l, shift right=2, "d_0"'] 
   \arrow[r, shorten >=1ex,shorten <=1ex] \arrow[r, shift left=2, shorten >=1ex,shorten <=1ex] 
   \arrow[r, shift right=2, shorten >=1ex,shorten <=1ex]
   & \cdots 
   \arrow[l, shift right=1] \arrow[l, shift left=1] \arrow[l, shift right=3] \arrow[l, shift left=3] 
 \end{tikzcd}.
  \]
  and geometrically, we have the following:
  \[
  \begin{tikzcd} 
    \overset{0}{\bullet} & \overset{1}{\bullet}
  \end{tikzcd}
  \] 
\end{example}

\begin{example}[2-Boundary]
  Following \cref{rem:non-degenerate}, $\Delta[2]$ has seven non-degenerate simplices. Among those, $0,1,2$ in $\Delta[2]_0$ and $01,12,02$ in $\Delta[2]_1$ are not surjective, while $012$ in $\Delta[2]_2$ is the identity map and hence surjective. Hence, we can depict $\partial \Delta[2]$ as the simplicial set: 
  {\small
  \[
  \begin{tikzcd}[column sep=1cm]
   \{0,1,2\} \arrow[r, shorten >=1ex,shorten <=1ex, "s_0"'  yshift= .2cm]
   & \{00,01,11,12,22,02\} 
   \arrow[l, shift left=1.4, "d_1"] \arrow[l, shift right=1.4, "d_0"'] 
   \arrow[r, shift right, shorten >=1ex,shorten <=1ex ] \arrow[r, shift left, shorten >=1ex,shorten <=1ex] 
   & \{000,001,011,111,112,122,222,002,022\} 
   \arrow[l] \arrow[l, shift left=2, "d_2"] \arrow[l, shift right=2, "d_0"'] 
   \arrow[r, shorten >=1ex,shorten <=1ex] \arrow[r, shift left=2, shorten >=1ex,shorten <=1ex] 
   \arrow[r, shift right=2, shorten >=1ex,shorten <=1ex]
   & \cdots 
   \arrow[l, shift right=1] \arrow[l, shift left=1] \arrow[l, shift right=3] \arrow[l, shift left=3] 
 \end{tikzcd}.
  \]
  }
  and geometrically, we have the following:
  \[ 
  \begin{tikzcd}[row sep=0.5cm, column sep=0.5cm]
    & \overset{1}{\bullet} \arrow[dr, "12"] & \\
    \overset{0}{\bullet} \arrow[ur, "01"] \arrow[rr, "02"'] &   & \overset{2}{\bullet}
  \end{tikzcd}
  \]
\end{example}

\begin{definition}[Horn]
  Let $n \geq 1$ and let $0 \leq i \leq n$. The \emph{$i$-th horn} of $\Delta[n]$, denoted $\Lambda[n]_i$, is the sub-simplicial set of $\Delta[n]$ with $(\Lambda[n]_i)_m \coloneq \{ \delta\colon [m] \to [n] \mid \mathrm{Im}(\delta) \cup \{i\} \neq [n] \}$. 
\end{definition}

From the definition, we see that $\Lambda[n]_i$ is always contained in $\partial \Delta[n]$. 

\begin{example}[1-Horns] \label{ex:lambda one horns}
 In the case $n = 1$, we have two horns, $\Lambda[1]_0$ and $\Lambda[1]_1$. $(\Lambda[1]_0)_0 \subseteq \Delta[1]_0 = \{0,1\}$ needs to exclude maps whose complement of the image only contains $0$. Hence, $(\Lambda[1]_0)_0 = \{0\}$. Similarly, $(\Lambda[1]_1)_0 = \{1\}$. Hence $\Lambda[1]_0 \cong \Lambda[1]_1 \cong \Delta[0]$, meaning geometrically we have a single vertex.
\end{example}

\begin{example}[2-Horns] \label{ex:lambda two horns}
  In the case $n=2$, we have three horns $\Lambda[2]_0$, $\Lambda[2]_1$, and $\Lambda[2]_2$. For $\Lambda[2]_0$, we make the following observations:
\begin{itemize}
  \item $(\Lambda[2]_0)_0$ needs to exclude maps whose complement of the image only contains $0$. But every map $[0] \to [2]$ has two elements in the complement of the image and hence $(\Lambda[2]_0)_0 = \{0,1,2\}$.
  \item $(\Lambda[2]_0)_1$ needs to exclude maps whose complement of the image only contains $0$. That is precisely the map $12\colon[1] \to [2]$. Hence, $(\Lambda[2]_0)_1 = \{00,01,11,02,22\}$.
\end{itemize}
 Hence, $\Lambda[2]_0$ is the simplicial set depicted as follows:
  \[
  \begin{tikzcd}[column sep=1.2cm]
   \{0,1,2\} \arrow[r, shorten >=1ex,shorten <=1ex, "s_0"'  yshift= .2cm]
   & \{00,01,11,22,02\} 
   \arrow[l, shift left=1.4, "d_1"] \arrow[l, shift right=1.4, "d_0"'] 
   \arrow[r, shift right, shorten >=1ex,shorten <=1ex ] \arrow[r, shift left, shorten >=1ex,shorten <=1ex] 
   & \{000,001,011,111,222,002,022\} 
   \arrow[l] \arrow[l, shift left=2, "d_2"] \arrow[l, shift right=2, "d_0"'] 
   \arrow[r, shorten >=1ex,shorten <=1ex] \arrow[r, shift left=2, shorten >=1ex,shorten <=1ex] 
   \arrow[r, shift right=2, shorten >=1ex,shorten <=1ex]
   & \cdots 
   \arrow[l, shift right=1] \arrow[l, shift left=1] \arrow[l, shift right=3] \arrow[l, shift left=3] 
 \end{tikzcd}.
  \]
  We can similarly describe $\Lambda[2]_1$ and $\Lambda[2]_2$. Geometrically, these three horns can hence be depicted as follows:
  \[ 
  \begin{tikzcd}[row sep=0.5cm, column sep=0.5cm]
    & \Lambda[2]_0 & \\[-.2cm]
    & \overset{1}{\bullet}  & \\
    \overset{0}{\bullet} \arrow[ur, "01"] \arrow[rr, "02"'] &   & \overset{2}{\bullet}
  \end{tikzcd}
  \qquad
  \begin{tikzcd}[row sep=0.5cm, column sep=0.5cm]
    & \Lambda[2]_1 & \\[-.2cm]
    & \overset{1}{\bullet} \arrow[dr, "12"] & \\
    \overset{0}{\bullet} \arrow[ur, "01"] &   & \overset{2}{\bullet}
  \end{tikzcd}
    \qquad
  \begin{tikzcd}[row sep=0.5cm, column sep=0.5cm]
    & \Lambda[2]_2 & \\[-.2cm]
    & \overset{1}{\bullet} \arrow[dr, "12"] & \\
    \overset{0}{\bullet} \arrow[rr, "02"']  &   & \overset{2}{\bullet}
  \end{tikzcd}
  \]
\end{example}

\begin{definition}[Spine] \label{def:spine}
  For every $n \geq 2$, the \emph{spine} of $\Delta[n]$, denoted $\Sp[n]$, is the sub-simplicial set of $\Delta[n]$ with $\Sp[n]_m \coloneq \{ \delta\colon [m] \to [n] \mid \delta(m) - \delta(0) \leq 1 \}$. This is indeed a sub-simplicial set as such maps are closed under precomposition.
\end{definition}

We can make the following two observations about the spine.

\begin{lemma} \label{lemma:spine zero}
  For every $n \geq 2$, $\Sp[n]_0 = \Delta[n]_0$.
\end{lemma}

\begin{proof}
  Every map $\sigma\colon [0] \to [n]$ necessarily satisfies $\sigma(0) - \sigma(0) = 0 \leq 1$ and hence is in $\Sp[n]_0$. 
\end{proof}

\begin{lemma} \label{lemma:spine bigger 2}
  For every $n,k \geq 2$, every element in $\Sp[n]_k$ is degenerate. 
\end{lemma}

\begin{proof}
  Every map $\sigma\colon [k] \to [n]$ with $k \geq 2$ and $\sigma(k) - \sigma(0) \leq 1$ must be non-injective, and is hence degenerate, by \cref{rem:non-degenerate}.
\end{proof}

\begin{example}[2-Spine]
 For $n = 2$, we can already evaluate $\Sp[2]_0$ (\cref{lemma:spine zero}) and $\Sp[2]_2$ (\cref{lemma:spine bigger 2}), and hence focus on $\Sp[2]_1$. We see that $\Sp[2]_1 = \{00,01,11,12,22\}$ as the map $02\colon [1] \to [2]$ has $\delta(1) - \delta(0) = 2$ and hence does not satisfy the condition. Hence $\Sp[2] \cong \Lambda[2]_1$. 
\end{example}

\begin{example}[3-Spine]
 For $n = 3$, we can already evaluate $\Sp[3]_0$ (\cref{lemma:spine zero}) and $\Sp[3]_2$ (\cref{lemma:spine bigger 2}), and hence focus on $\Sp[3]_1$. Similar to the previous example, we can hence see that $\Sp[3]_1 = \{00,01,11,12,22,23,33\}$. Hence, geometrically, $\Sp[3]$ can be depicted as follows:
  \[ 
  \begin{tikzcd}[row sep=0.5cm, column sep=0.5cm]
    & \overset{3}{\bullet}& \\
    \overset{0}{\bullet} \arrow[rr, "01" near end] &   & \overset{1}{\bullet} \arrow[dl, "12"]  \\[-.2cm]
     & \overset{2}{\bullet} \arrow[uu, "23" near end] & \\
  \end{tikzcd}
  \quad 
  \xhookrightarrow[]{ \text{sits inside}} 
  \quad
  \begin{tikzcd}[row sep=0.5cm, column sep=0.5cm]
    & \overset{3}{\bullet}& \\
    \overset{0}{\bullet} \arrow[ur, "03"] \arrow[dr, "02"] \arrow[rr, "01" near end] &   & \overset{1}{\bullet} \arrow[dl, "12"] \arrow[ul, "13"'] \\[-.2cm]
     & \overset{2}{\bullet} \arrow[uu, "23" near end] & \\
  \end{tikzcd}
  \]
  where the geometric shape of $\Sp[3]$ was chosen in a manner to illustrate how it sits inside $\Delta[3]$ as a sub-simplicial set.
\end{example}

Let us move on to the last topic, namely limits and colimits in $\sSet$. First, let us summarize the main results regarding limits and colimits, as reviewed in \cref{subsec:limits}, in the particular case of simplicial sets.

\begin{definition}[Colimit of simplicial sets] \label{def:colimit sset}
  Let 
  \[X_1 \xleftarrow{f_{1,2}} X_{1,2} \xrightarrow{g_{1,2}} X_2 \xleftarrow{f_{2,3}} X_{2,3} ... \xleftarrow{f_{n-1,n}} X_{n-1,n}\xrightarrow{g_{n-1,n}}  X_n \]
   be a diagram of simplicial sets. The \emph{colimit} of this diagram is the simplicial set denoted 
   \[X_1 \coprod_{X_{1,2}}^{f_{1,2},g_{1,2}} ... \coprod_{X_{n-1,n}}^{f_{n-1,n},g_{n-1,n}} X_n,\] 
   whose $k$-simplices are given by
  \[(X_1 \coprod_{X_{1,2}}^{f_{1,2},g_{1,2}} ... \coprod_{X_{n-1,n}}^{f_{n-1,n},g_{n-1,n}} X_n)_k \coloneq (\coprod_{1 \leq i \leq n}(X_i)_k) / \sim,\]
  where $\coprod$ is the disjoint union of the set of $k$-simplices in $X_i$ and  $\sim$ is the equivalence relation generated by $(f_{i,i+1})_k(x) \sim (g_{i,i+1})_k(x)$ for every $x$ in $(X_{i,i+1})_k$.
\end{definition}

\begin{definition}[Limit of simplicial sets] \label{def:limit sset}
 Let 
  \[X_1 \xrightarrow{f_{1,2}} X_{1,2} \xleftarrow{g_{1,2}} X_2 \xrightarrow{f_{2,3}} X_{2,3} ... \xrightarrow{f_{n-1,n}} X_{n-1,n} \xleftarrow{g_{n-1,n}} X_n \]
  be a diagram of simplicial sets. The \emph{limit} of this diagram is the simplicial set denoted 
  \[X_1 \times_{X_{1,2}}^{f_{1,2},g_{1,2}} \cdots \times_{X_{n-1,n}}^{f_{n-1,n},g_{n-1,n}} X_n,\] 
  whose $k$-simplices are given by
  \begin{multline*}
    (X_1 \times_{X_{1,2}}^{f_{1,2},g_{1,2}} \cdots \times_{X_{n-1,n}}^{f_{n-1,n},g_{n-1,n}} X_n)_k \\
    \coloneq \{(x_1,\dots,x_n)\in \prod_{1 \le i \le n} (X_i)_k \mid (f_{i,i+1})_k(x_i) = (g_{i,i+1})_k(x_{i+1}) \text{ for all } 1 \le i \le n-1\},
  \end{multline*}
 where $\times$ is the Cartesian product of the set of $k$-simplices in $X_i$.
\end{definition}

\begin{remark}
  Following \cref{prop:colimit limit,prop:limit limit}, and further explained in \cref{rem:colimit limit}, maps out of the colimit or into a limit are uniquely given by a compatible tuple of maps out of (or into) the individual simplicial sets.
\end{remark}

We now move on to the explicit computations of limits and colimits of interest in $\sSet$.

\begin{example}
  The simplicial set $\Delta[1] \times \Delta[1]$ is defined as the product of simplicial sets. Hence, following \cref{def:limit sset}, for every $n \geq 0$, we have $(\Delta[1] \times \Delta[1])_n = \Delta[1]_n \times \Delta[1]_n$. This means that for the first three levels we have the following explicit computations:
  {\small
  \begin{itemize}[leftmargin=*]
    \item $(\Delta[1] \times \Delta[1])_0 = \{(0,0), (0,1), (1,0), (1,1)\}$
    \item $(\Delta[1] \times \Delta[1])_1 = \{(00,00), (00,01), (00,11), (01,00), (01,01), (01,11), (11,00), (11,01), (11,11)\}$
    \item $(\Delta[1] \times \Delta[1])_2 = \{(000,000), (000,001), (000,011), (000,111), (001,000), (001,001), (001,011),$ 
    \item[] $(001,111), (011,000), (011,001), (011,011), (011,111), (111,000), (111,001), (111,011), (111,111)\}$
  \end{itemize}
  }
  Interestingly, the non-degenerate simplices of $\Delta[1] \times \Delta[1]$ are given by all $0$-simplices, the five $1$-simplices $(00,01)$, $(11,01)$, $(01,00)$, $(01,11)$, and $(01,01)$, and the two $2$-simplices $(011,001)$ and $(001,011)$. Geometrically, we can hence depict $\Delta[1] \times \Delta[1]$ as follows:
  \[ 
  \begin{tikzcd}[row sep=1cm, column sep=1cm]
    \overset{(0,0)}{\bullet} \arrow[r, "{(00,01)}"] \arrow[d, "{(01,00)}"'] \arrow[dr, "{(01,01)}" description] & \overset{(0,1)}{\bullet} \arrow[d, "{(01,11)}"] \\
    \overset{(1,0)}{\bullet} \arrow[r, "{(11,01)}"'] & \overset{(1,1)}{\bullet}
  \end{tikzcd}
  \]
  where the $2$-simplex $(011,001)$ is the bottom left triangle, and $(001,011)$ is the top right triangle in the square.
\end{example}

\begin{remark}\label{rem:degenerate in products}
  Even though $(00,01)$ is non-degenerate in $\Delta[1] \times \Delta[1]$, $00$ in $\Delta[1]_1$ on its own is not. This shows that a correct notion of product that is compatible with our geometric intuition requires the presence of degenerate simplices.
\end{remark}

\begin{lemma} \label{lemma:boundary colimit}
  There is a bijection of simplicial sets $\partial \Delta[1] \cong \Delta[0] \coprod \Delta[0]$.
\end{lemma}

\begin{proof}
 By \cref{prop:colimit limit}, a map $\Delta[0] \coprod \Delta[0] \to \partial \Delta[1]$ is uniquely given by a pair of maps $\Delta[0] \to \partial \Delta[1]$. By the Yoneda lemma (\cref{lemma:delta yoneda}), each one of these maps corresponds to an element in $\partial \Delta[1]_0 = \{0,1\}$. Hence, we have the map $0 \coprod 1 \colon \Delta[0] \coprod \Delta[0] \to \partial \Delta[1]$, which by construction is an injection. Finally, following \cref{def:colimit sset}, $(\Delta[0] \coprod \Delta[0])_n = \{00...0\} \coprod \{ 00 ... 0\}$ has two elements, and so the map is a bijection.
\end{proof}

\begin{lemma} \label{lemma:horn colimit}
  We have three bijections of simplicial sets.
  \begin{itemize}
    \item $\Lambda[2]_0 \cong \Delta[1] \coprod^{0,0}_{\Delta[0]} \Delta[1]$
    \item $\Lambda[2]_1 \cong \Delta[1] \coprod^{1,0}_{\Delta[0]} \Delta[1]$
    \item $\Lambda[2]_2 \cong \Delta[1] \coprod^{1,1}_{\Delta[0]} \Delta[1]$
  \end{itemize}
\end{lemma}

\begin{proof}
 We focus on the first case. The other cases are analogous. We start by constructing a map of simplicial sets $\Delta[1] \coprod^{0,0}_{\Delta[0]} \Delta[1] \to \Lambda[2]_0$. First observe we have the following commutative diagram
 \[ 
 \begin{tikzcd}
  \Delta[0] \arrow[r, "0"] \arrow[d, "0"'] & \Delta[1] \arrow[d, "01"] \\
  \Delta[1] \arrow[r, "02"] & \Lambda[2]_0
 \end{tikzcd},
 \]
 where we used the fact that $(\Lambda[2]_0)_1 = \{00, 01, 11, 02, 22\}$ and the Yoneda lemma (\cref{lemma:delta yoneda}) to construct the map $\Delta[1] \to \Lambda[2]_0$. Here the commutativity follows from the fact that the image of $0 \in \Delta[0]_0$ in both directions is $0 \in (\Lambda[2]_0)_0$. 

 So, by \cref{prop:colimit limit}, there is a unique map $01 + 02\colon\Delta[1] \coprod_{\Delta[0]} \Delta[1] \to \Lambda[2]_0$. We now directly observe that at level $k$, the map $01 + 02$ is given by 
 \[\{00... 0, ... , 11 ... 1\} \coprod_{00...0} \{00... 0, ... , 11 ...1\} \to \{00... 0, 11...1,00...02, ... , 22...2\}. \]
 Concretely, the map is given by the identity on the first component and maps $00...111$ to $00...222$ on the second component, and is hence a bijection.
\end{proof}

\begin{lemma} \label{lemma:spine colimit}
 Let $n \geq 2$. There is a bijection of simplicial sets 
 \[\Sp[n] \cong \Delta[1] \coprod_{\Delta[0]} ... \coprod_{\Delta[0]} \Delta[1] \]
\end{lemma}

\begin{proof}
 We proceed by induction. The case $n =2$ is already covered by the previous lemma. For the induction step, we repeat the steps from the previous proof, replacing $\Delta[1]$ with $\Sp[n-1]$ and $\Lambda[2]_0$ with $\Sp[n]$.
\end{proof}

\section{Category Theory \texorpdfstring{$\&$}{\&} Homotopy Theory via Simplicial Sets} \label{sec:ct ht via ss}
In the previous section we introduced simplicial sets. We now want to see how we can use simplicial sets to study both categories and topological spaces.

\subsection{From Category Theory to Simplicial Sets} \label{subsec:ct to ss}
We now start relating categories to simplicial sets. Here the key idea is the \emph{nerve} of a category. For more technical introductions see \cite[Section 4]{rezk2022qcats} or \cite[\href{https://kerodon.net/tag/002M}{Subsection 002M}]{lurie2026kerodon}.

By definition of $\DD$ (\cref{def:delta}), there is a functor $\DD \to \Cat$, which is the identity on objects and morphisms. We now have the following definition.

\begin{definition}[Nerve] \label{def:nerve}
  Let $N\colon \Cat \to \sSet$ be the functor given by \cref{constr:restricted yoneda} applied to the functor $\DD \to \Cat$.
\end{definition}

\begin{remark} \label{rem:nerve explicit}
  More explicitly, we can describe the nerve as follows. 
  \begin{itemize}
    \item $n = 0$: $N\cC_0 = \Fun([0],\cC) = \Obj_{\cC}$ is the set of objects of $\cC$.
    \item $n = 1$: $N\cC_1 = \Fun([1],\cC) = \Mor_{\cC}$ is the set of morphisms of $\cC$.
    \item $n \geq 2$: $N\cC_n = \Fun([n],\cC)$ is the set of $n$ composable morphisms in $\cC$. In other words, an element of $N\cC_n$ is a sequence of morphisms
    \[x_0 \xrightarrow{f_1} x_1 \xrightarrow{f_2} x_2 \xrightarrow{f_3} ... \xrightarrow{f_n} x_n\]
    where $x_i \in \Obj_\cC$ and $f_i\colon x_{i-1} \to x_i$.
  \end{itemize}
\end{remark}

\begin{lemma}
 Let $n \geq 0$. Then $N([n]) \cong \Delta[n]$.
\end{lemma}

\begin{proof}
 By definition we have $N([n])_k = \Hom_{\Cat}([k],[n]) = \Hom_\DD([k],[n]) = \Delta[n]_k$. Analogously, $N([n])$ and $\Delta[n]$ have the same face and degeneracy maps. 
\end{proof}

As a functor the nerve construction has very nice injectivity properties.

\begin{proposition} \label{prop:nerve fully faithful}
 The nerve is fully faithful (\cref{def:fully faithful}).
\end{proposition}

\begin{proof}
  Let $\cC$ and $\cD$ be two categories. We need to show that the map
  \[
  N\colon \Hom_{\Cat}(\cC,\cD) \to \Hom_{\sSet}(N\cC,N\cD)
  \]
  is a bijection. Observe that for a given functor $F\colon \cC\to \cD$, the map $N(F)_0\colon \Obj_{\cC} \to \Obj_{\cD}$ is the functor $F$ on objects, and $N(F)_1\colon \Mor_{\cC} \to \Mor_{\cD}$ is the functor $F$ on morphisms. Hence, for two functors $F,G\colon \cC\to \cD$ with $N(F)=N(G)$, we have $F_0=G_0$ and $F_1=G_1$, meaning $F = G$. Hence $N$ is injective.

  It remains to show that $N$ is surjective. Given a map $f\colon N\cC \to N\cD$ of simplicial sets, we will define a functor $F\colon \cC\to \cD$, such that $N(F)=f$. We define $F$ as follows.
  \begin{itemize}[leftmargin=*]
    \item For an object $x\in \Obj_{\cC}=(N\cC)_0$, define $F(x)\coloneq f_0(x)$.
    \item For a morphism $g\colon x\to y$ in $\cC$, define $F(g) \coloneq f_1(g)$.
    \item[] Here the simplicial relations
    \[
      0^*(f_1(g))=f_0(0^*(g))=f_0(x)=F(x),\qquad
      1^*(f_1(g))=f_0(1^*(g))=f_0(y)=F(y),
    \]
    imply that if $g\colon x\to y$, then $F(g)\colon F(x)\to F(y)$.
    \item For a given object $x$ in $\cC$, 
    \[
      F(\id_x) =  f_1(00^*(x)) = 00^*(f_0(x)) = 00^*(F(x)) = \id_{F(x)},
    \]
    hence $F$ preserves identities.
    \item For composable morphisms $x\xrightarrow{g}y\xrightarrow{h}z$ in $\cC$, we have
    {\small
    \[
      F(h\circ g) = f_1(h\circ g) = f_1(02^*(g,h)) = 02^*(f_2(g,h)) = 02^*(f_1(g),f_1(h)) = f_1(h) \circ f_1(g) = F(h)\circ F(g),
    \]
    }
    hence $F$ preserves composition.
  \end{itemize}

  Finally, we confirm that $N(F)=f$. By definition $N(F)_0 = f_0$ and $N(F)_1 = f_1$. More generally let $\sigma\in (N\cC)_n$ be an $n$-simplex, which, by \cref{rem:nerve explicit}, is a string of $n$ composable morphisms
  \[
  \sigma=(x_0\xrightarrow{g_1}x_1\xrightarrow{g_2}\cdots\xrightarrow{g_n}x_n).
  \]
  Then the simplicial identities imply that 
  \[ f_n(\sigma) = (f_1(g_1),..., f_1(g_n)) = (F(g_1),..., F(g_n)) = N(F)_n(\sigma). \]
  This finishes the proof.
\end{proof}

This means we can see categories as special kinds of simplicial sets. But which simplicial sets do we get this way? Here we use the inclusion $i_n\colon\Sp[n] \to \Delta[n]$ (\cref{def:spine}) and the precomposition map (\cref{ex:repcontra}).

 \begin{definition}[Segal condition] \label{def:Segal condition}
  A simplicial set $X$ satisfies the \emph{Segal condition} if the map
  \[i^*_n\colon\Hom_{\sSet}(\Delta[n],X) \to \Hom_{\sSet}(\Sp[n],X) \]
  is a bijection for $n \geq 2$.
 \end{definition}

\begin{lemma} \label{lemma:spine pullback}
 Let $X$ be a simplicial set and let $n \geq 2$. Then there is a bijection of sets
 \[ \Hom_{\sSet}(\Sp[n],X) \cong X_1 \times_{X_0}^{1^*,0^*} X_1 \times_{X_0}^{1^*,0^*} ... \times_{X_0}^{1^*,0^*} X_1. \]
\end{lemma}

\begin{proof}
 We have the following chain of bijections:
 {\small
 \begin{align*}
  \Hom_{\sSet}(\Sp[n],X) & \cong \Hom_{\sSet}(\Delta[1] \coprod^{1,0}_{\Delta[0]} ... \coprod^{1,0}_{\Delta[0]} \Delta[1],X) \hspace{-.8em} & \text{\cref{lemma:spine colimit}}\\
  & \cong \Hom_{\sSet}(\Delta[1],X) \underset{\Hom_{\sSet}(\Delta[0],X)}{\times} ... \underset{\Hom_{\sSet}(\Delta[0],X)}{\times} \Hom_{\sSet}(\Delta[1],X) \hspace{-.8em} & \text{\cref{prop:colimit limit}}\\
  & \cong X_1 \times_{X_0}^{1^*,0^*} ... \times_{X_0}^{1^*,0^*} X_1 \hspace{-.8em} & \text{\cref{lemma:delta yoneda}}
 \end{align*}
 }
\end{proof}

\begin{proposition} \label{prop:Segal condition pullback}
 Let $X$ be a simplicial set. Then the following are equivalent:
 \begin{enumerate}[leftmargin=*]
  \item $X$ satisfies the Segal condition.
  \item For $n \geq 2$, the composition map
  \[X_n \cong \Hom(\Delta[n],X) \xrightarrow{i_n^*} \Hom(\Sp[n],X) \cong X_1 \times_{X_0}^{1^*,0^*} X_1 \times_{X_0}^{1^*,0^*} ... \times_{X_0}^{1^*,0^*} X_1, \]
  is a bijection.
 \end{enumerate}  
\end{proposition}

\begin{proof}
 We have the following diagram of sets 
 \[ 
 \begin{tikzcd}
  \Hom_{\sSet}(\Delta[n],X) \arrow[r] \arrow[d, "\cong"] & \Hom_{\sSet}(\Sp[n],X) \arrow[d, "\cong"] \\
  X_n \arrow[r] & X_1 \times_{X_0}^{1^*,0^*} ... \times_{X_0}^{1^*,0^*} X_1
 \end{tikzcd}
 \]
 Here the vertical maps are bijections, by \cref{lemma:delta yoneda,lemma:spine pullback}. Hence, the top horizontal map is a bijection, which is $(1)$, if and only if the bottom horizontal map is a bijection, which is $(2)$.
\end{proof}

We now have the following important observation.
\begin{proposition} \label{prop:nerve Segal}
 Let $\cC$ be a category. Then the nerve $N\cC$ satisfies the Segal condition.  
\end{proposition}

\begin{proof}
 By definition $N\cC_n = \Hom_{\Cat}([n],\cC)$. Following \cref{lemma:functor as sequence of morphisms}, a functor $[n] \to \cC$ is precisely the data of $n$ composable morphisms $(f_1, ... f_n)$ in $\cC$. By the definition of limit, this is precisely an element in the limit $N\cC_1 \times_{N\cC_0}^{1^*,0^*} N\cC_1 \times_{N\cC_0}^{1^*,0^*} ... \times_{N\cC_0}^{1^*,0^*} N\cC_1$. Hence, $N\cC$ satisfies the Segal condition.
\end{proof}

Before we move on, this observation can help us compute the nerve in many new situations. Here we use the product of categories (\cref{def:product category}).

\begin{proposition} \label{prop:nerve product}
  Let $\cC, \cD$ be two categories. Then $N(\cC \times \cD) \cong N\cC \times N\cD$.
\end{proposition}

\begin{proof}
 Following \cref{rem:nerve explicit,def:product category}, we have bijections
  \[N(\cC \times \cD)_0 \cong \Obj_{\cC \times \cD} \cong \Obj_{\cC} \times \Obj_{\cD} \cong N\cC_0 \times N\cD_0, \]
  \[N(\cC \times \cD)_1 \cong \Mor_{\cC \times \cD} \cong \Mor_{\cC} \times \Mor_{\cD} \cong N\cC_1 \times N\cD_1. \]
 The bijection for higher $n$ then follows from the Segal condition (\cref{prop:Segal condition pullback}).
\end{proof}

What is fascinating is that the other direction also holds.

\begin{theorem} \label{thm:nerve image}
  A simplicial set $X$ is isomorphic to the nerve of a category if and only if $X$ satisfies the Segal condition.
\end{theorem}

\begin{proof}
 We already proved one direction in \cref{prop:nerve Segal}. Let $X$ be a simplicial set that satisfies the Segal condition. Using the Segal conditions, we fix the notation and denote the inverse Segal maps by $\mu_2\colon X_1 \times_{X_0} X_1 \to X_2$ and $\mu_3\colon X_1 \times_{X_0} X_1 \times_{X_0} X_1 \to X_3$. We construct a category $\cC$ such that $N\cC \cong X$ as follows.
\begin{itemize}[leftmargin=*]
  \item $\Obj_{\cC} = X_0$.
  \item For two objects $x,y$ in $\Obj_{\cC}$, the set of morphisms $\Hom_{\cC}(x,y)$ is given by the pullback
  \[
  \begin{tikzcd}
   \Hom_{\cC}(x,y) \arrow[r] \arrow[d] & X_1 \arrow[d, "{(0^*,1^*)}"] \\
   \Delta[0] \arrow[r, "{(x,y)}"] & X_0 \times X_0
  \end{tikzcd}.
  \]
  \item The identity morphism on an object $x$ is given by $00^*(x)$.
  \item For two morphisms $f\colon x \to y$ and $g\colon y \to z$, we define the composition as $g \circ f = 02^*\mu_2(f,g)$.
  \item For three morphisms $f,g,h$ that are composable, we have 
  \[(h \circ g) \circ f = 02^*\mu_2(f,02^*\mu_2(g,h)) = 03^*\mu_3(f,g,h) = 02^*\mu_2(02^*\mu_2(f,g),h) = h \circ (g \circ f).\] 
  Hence, composition is associative.
  \item For a morphism $f\colon x \to y$, we have 
  \[f \circ \id_x = 02^*\mu_2(00^*(x),f) = f \text{ and } \id_y \circ f = 02^*\mu_2(f,00^*(y)) = f.\] 
  Hence, composition is unital.
\end{itemize}
  Finally, we need to show that $N\cC \cong X$. For $n = 0$, we have $N\cC_0 = \Obj_{\cC} = X_0$. For $n = 1$, we have $N\cC_1 = \Mor_{\cC} = X_1$. For $n \geq 2$, we have 
  \[N\cC_n \cong N\cC_1 \times_{N\cC_0}^{1^*,0^*} ... \times_{N\cC_0}^{1^*,0^*} N\cC_1 \cong X_1 \times_{X_0}^{1^*,0^*} ... \times_{X_0}^{1^*,0^*} X_1 \cong X_n. \]
  Hence, $N\cC \cong X$ and we are done.
\end{proof}

\begin{remark} \label{rem:nerve hom set}
  As part of the construction we in particular observed that for two objects $x,y$ in $\cC$, we can reconstruct the hom-set $\Hom_{\cC}(x,y)$ via the nerve of $\cC$ as the pullback (\cref{ex:pullback})
  \[
  \begin{tikzcd}
   \Hom_{\cC}(x,y) \arrow[r] \arrow[d] & N\cC_1 \arrow[d, "{(0^*,1^*)}"] \\
   \Delta[0] \arrow[r, "{(x,y)}"] & N\cC_0 \times N\cC_0
  \end{tikzcd}.
  \]
\end{remark}
Combining these results we can conclude the following.

\begin{definition}
  Let $\sSet^{Seg}$ be the category whose objects are simplicial sets satisfying the Segal condition and whose morphisms are maps of simplicial sets.
\end{definition}

\begin{corollary} \label{cor:nerve equivalence}
  The nerve $N\colon \Cat \to \sSet^{Seg}$ is an equivalence of categories.
\end{corollary}

\begin{proof}
 By \cref{prop:nerve fully faithful}, $N$ is fully faithful and by \cref{thm:nerve image}, $N$ is essentially surjective. Hence, $N$ is an equivalence of categories, by \cref{thm:equivalence of categories}.
\end{proof}

Thus we have shown that, up to equivalence, category theory is simply the study of a particular class of simplicial sets, namely those that satisfy the Segal condition.

\subsection{From Algebraic Topology to Simplicial Sets} \label{subsec:top to ss}
We now want to see how we can use simplicial sets to study the homotopy theory of topological spaces. As we primarily want to understand how these two areas are related, we will review relevant facts from algebraic topology without proof, providing references when needed. Further details can of course be found in any standard textbook on algebraic topology, such as \cite{hatcher2002at}.

Algebraic topology studies topological spaces via algebraic invariants. Prominent examples include the \emph{homotopy groups}, which we review now.

\begin{definition}[Homotopy]
  Let $f,g\colon X \to Y$ be two continuous maps between topological spaces. A \emph{homotopy} from $f$ to $g$ is a continuous map $H\colon X \times [0,1] \to Y$ such that $H|_{X \times \{ 0\} } = f$ and $H|_{X \times \{ 1\} } = g$. We say that $f$ and $g$ are \emph{homotopic} if there exists a homotopy from $f$ to $g$.
\end{definition}

For this next definition recall pointed topological spaces (\cref{ex:pointed top}).
\begin{lemma}[{\cite[Section 1.1]{hatcher2002at}}]
 Let $(X,x), (Y,y)$ be two pointed topological spaces. Then the relation of being homotopic is an equivalence relation on $\Hom_{\Top_*}((X,x),(Y,y))$. 
\end{lemma}

\begin{definition}[Homotopy groups] \label{def:homotopy groups}
  Let $(X,x)$ be a pointed topological space. For $n \geq 0$, the $n$-th homotopy set of $X$ at $x$, denoted $\pi_n(X,x)$, is the set of homotopy classes of maps from the $n$-sphere $S^n$ to $X$ that send the base point of $S^n$ to $x$.
\end{definition}

\begin{remark}
  Notice that $\pi_0(X,x)$, as defined in \cref{def:homotopy groups}, is independent of the choice of base point $x$, and is in fact the set of path-connected components of $X$. Hence, we will often write $\pi_0(X)$ instead of $\pi_0(X,x)$.
\end{remark}

\begin{remark}
  It is a classical result that $\pi_n(X,x)$ is a group for $n \geq 1$, and is abelian for $n \geq 2$ (\cite[Section 4.1]{hatcher2002at}). Hence, following common convention, we will also adopt the term \emph{homotopy group}. However, the group structure has no relevance for our purposes.
\end{remark}

We can now define new invariants of topological spaces via homotopy groups.

\begin{definition}[Weak homotopy equivalence] \label{def:weak homotopy equivalence}
  A continuous map $f\colon X \to Y$ between topological spaces is a \emph{weak homotopy equivalence} if it induces a bijection $\pi_0(f)\colon \pi_0(X) \to \pi_0(Y)$ and, for every base point $x \in X$,  a bijection $\pi_n(f)\colon \pi_n(X,x) \to \pi_n(Y,f(x))$ for all $n \geq 1$.
\end{definition}

One major question in algebraic topology is to determine whether two topological spaces are \emph{weakly homotopy equivalent}, via some continuous map that is a weak homotopy equivalence. However, determining whether a given map is a weak homotopy equivalence is often impossibly difficult. One major breakthrough in this direction is the introduction of \emph{CW-complexes}, and notably the \emph{CW-approximation theorem}, as proven in \cite[Proposition 4.13]{hatcher2002at}.

\begin{definition}[CW-complex]
  A topological space $X$ is a \emph{CW-complex} if it is of the form $X = \bigcup_{n \geq 0} X_n$, where $X_0$ is a discrete topological space and, for $n \geq 1$, $X_{n}$ is obtained from $X_{n-1}$ by attaching $n$-cells $\partial D^{n} \to D^{n}$ along continuous maps $\partial D^n \to X_{n-1}$. 
\end{definition}
  
\begin{theorem}[CW-approximation theorem]
  For every topological space $X$, there exists a CW-complex $\tilde{X}$ and a weak homotopy equivalence $\tilde{X} \to X$.
\end{theorem}

This means that, if we only care about topological spaces up to weak homotopy equivalence, we can restrict our attention to CW-complexes. Fortunately, the inductive nature of CW-complexes allows us to effectively determine weak homotopy equivalences between them.

\begin{definition}[Homotopy equivalence]
  A continuous map $f\colon X \to Y$ between topological spaces is a \emph{homotopy equivalence} if there exists a continuous map $g\colon Y \to X$ such that both compositions $f \circ g$ and $g \circ f$ are homotopic to the identity.
\end{definition}

We now have the following first lemma relating homotopy equivalences and weak homotopy equivalences (\cref{def:weak homotopy equivalence}). Concretely, we can rely on standard properties of homotopy groups, such as $\pi_n$ being a functor and identifying homotopic maps, to show that homotopy equivalences are weak homotopy equivalences.

\begin{lemma} \label{lemma:homotopy implies weak homotopy}
  Every homotopy equivalence is a weak homotopy equivalence.
\end{lemma}

\begin{proof}
 If $f\colon X \to Y$ is a homotopy equivalence, then there exists $g\colon Y \to X$ such that $f \circ g$ and $g \circ f$ are homotopic to the identity. Hence, for every base point $y \in Y$, we have 
 \[\pi_n(f,g(y)) \circ \pi_n(g,y) = \pi_n(f \circ g,y) = \pi_n(\id_Y,y),\] 
 where for the last equality we have identified $\pi_n(Y,y)$ with $\pi_n(Y,f \circ g(y))$, as $y$ and $f \circ g(y)$ are path-connected. Similarly, for every $x \in X$, $\pi_n(g,f(x)) \circ \pi_n(f,x) = \pi_n(\id_X,x)$, which means that $\pi_n(f,x)\colon \pi_n(X,x) \to \pi_n(Y,f(x))$ is a bijection for all $n \geq 1$.

 We can similarly observe that $\pi_0(f)$ is a bijection on path-components. Hence, we are done.
\end{proof}

The converse does not hold in general. See \cite[Exercise 1.3.7]{hatcher2002at} for a space that is weakly equivalent, but not homotopy equivalent, to a point. However, for CW-complexes, we have the following major result \cite[Theorem 4.5]{hatcher2002at}.

\begin{theorem}[Whitehead's theorem]
  A weak homotopy equivalence between CW-complexes is a homotopy equivalence.
\end{theorem}

This state of affairs naturally leads us to the following question: How can we effectively determine homotopy equivalences of CW-complexes? Here simplicial sets can provide a powerful approach.

\begin{definition}[Topological n-simplex]
  For $n \geq 0$, let $|\Delta[n]|$ denote the topological $n$-simplex, which is the topological subspace of $\bR^{n+1}$ given by 
  \[|\Delta[n]| \coloneq \{(t_i)_{0 \leq i \leq n} \in \bR^{n+1} \mid \sum_{i = 0}^n t_i = 1, t_i \geq 0 \text{ for all } i \}. \] 
\end{definition}

\begin{lemma} \label{lemma:restricted geometric realization}
  The assignment $[n] \mapsto |\Delta[n]|$ defines a functor $|\Delta[\bullet]|\colon \DD \to \Top$.
\end{lemma}

\begin{proof}
 We construct the functor on the morphisms $d^i,s^i$ and then show the cosimplicial identities are satisfied (\cref{lemma:cosimplicial identities}). Let $d^i\colon |\Delta[n]| \to |\Delta[n+1]|, s^i\colon |\Delta[n]| \to |\Delta[n-1]|$ be defined as follows:
 \[
 d^i(t_0,...,t_n) = (t_0,...,t_{i-1},0,t_i,...,t_n)
 \]
 \[
 s^i(t_0,...,t_n) = (t_0,...,t_{i-1},t_i + t_{i+1},t_{i+2},...,t_n).
 \]
 Checking the cosimplicial identities is very straightforward. For example, 
 \[s^i \circ d^i(t_0,...,t_n) = s^i(t_0,...,t_{i-1},0,t_i,...,t_n) = (t_0,...,t_{i-1},t_i + 0,t_{i+1},...,t_n) = (t_0,...,t_n),\]
 and we leave the rest to the reader.
\end{proof}
 
\begin{definition}
 Let $\Sing\colon \Top \to \sSet$ be the functor defined by applying \cref{constr:restricted yoneda} to $|\Delta[\bullet]|$.
\end{definition}

\begin{remark} \label{rem:sing explicit}
  More explicitly, for a given topological space $X$, we have the equality $\Sing(X)_n = \Hom_{\Top}(|\Delta[n]|,X)$.
\end{remark}

Having defined $\Sing$, we now want to understand how well it recovers the homotopy theory of topological spaces. For that we first develop a homotopy theory of simplicial sets, and then we show that $\Sing$ is compatible with it.

\begin{definition}
  Two maps $f,g\colon X \to Y$ between simplicial sets are \emph{homotopic} if there exists a map $H\colon X \times \Delta[1] \to Y$ such that $H|_{X \times \{ 0\} } = f$ and $H|_{X \times \{ 1\} } = g$. 
\end{definition}

\begin{definition}[Homotopy equivalence]
  A map $f\colon X \to Y$ between simplicial sets is a \emph{homotopy equivalence} if there exists a map $g\colon Y \to X$ such that both compositions $f \circ g$ and $g \circ f$ are homotopic to the identity.
\end{definition}

One early observation is that $\Sing$ is compatible with homotopies, which in particular will follow from explicit understanding of the functor $\Sing$ \cite[Chapter I]{goerssjardine2009simplicialhomotopytheory}.

\begin{proposition} \label{prop:sing homotopy}
  The functor $\Sing$ sends homotopic maps to homotopic maps, and preserves homotopy equivalences. 
\end{proposition}

Thus we can study homotopy theory of topological spaces via simplicial sets. But do we really need all simplicial sets? Or can we restrict to a smaller class? Here we make the following observations.

\begin{definition}[Kan complex] \label{def:kan complex}
 A simplicial set $K$ is a \emph{Kan complex} if for every $n \geq 1$ and $0 \leq k \leq n$, every map $\Lambda[n]_k \to K$ extends to a map $\Delta[n] \to K$
 \[
 \begin{tikzcd}
    \Lambda[n]_k \arrow[d] \arrow[r] & K \\
    \Delta[n] \arrow[ur, dashed]
   \end{tikzcd},
 \]
 or, equivalently, the map 
 \[\Hom(\Delta[n], K) \to \Hom(\Lambda[n]_k, K) \]
 is surjective.
\end{definition}

\begin{proposition}[{\cite[Proposition 2.3]{goerssjardine2009simplicialhomotopytheory}}] \label{prop:sing kan}
  For every topological space $X$, the simplicial set $\Sing(X)$ is a Kan complex.
\end{proposition}

\begin{definition}[Category of Kan complexes]
  Let $\Kan$ be the category whose objects are Kan complexes and whose morphisms are maps of simplicial sets.
\end{definition}

\cref{prop:sing kan} shows that we have a functor $\Sing\colon \Top \to \Kan$. Unlike the nerve functor, this functor is not an equivalence of categories. However, if we only consider homotopically relevant data, then things change. Let us see an example of that.

\begin{definition}[Set of path-components] \label{def:pizero sset}
  Let $S$ be a simplicial set. Let $\sim$ be the relation on $S_0$ given by $x \sim y$ if there exists an $f$ in $S_1$, meaning a map $\Delta[1] \to S$, such that $0^*(f) = x$ and $1^*(f) = y$. Let $\pi_0(S) = S_0/\sim$, called the set of \emph{path-components of $S$}, be the quotient of $S_0$ with respect to the equivalence relation generated by $\sim$.
\end{definition}

For Kan complexes, \cref{def:pizero sset} behaves as one might expect, demonstrating the homotopical nature of Kan complexes.

\begin{lemma} \label{lemma:homotopy equivalence relation}
 Let $K$ be a Kan complex. Then $\sim$ is already an equivalence relation on $K_0$.
\end{lemma}

\begin{proof}
  We need to show the relation is reflexive, symmetric, and transitive. Reflexivity is witnessed by the constant map $\Delta[1] \to K$. Now, let $H\colon\Delta[1] \to K$ be a homotopy from $f$ to $g$. Let $\sigma\colon \Lambda[2]_0 \to K$ be the map sending $01$ to $H$ and $02$ to $\id_{f}$, where we used the description of $\Lambda[2]_0$ as a colimit (\cref{lemma:horn colimit}). As $K$ is a Kan complex, there is a lift $H'\colon \Delta[2] \to K$. By the simplicial identities, $12^*H'\colon\Delta[1] \to K$ gives us the desired homotopy from $g$ to $f$.

  Similarly, if $H_1$ is a homotopy from $f$ to $g$ and $H_2$ is a homotopy from $g$ to $h$, we can construct a homotopy from $f$ to $h$ by first constructing the map $H_1 + H_2\colon \Lambda[2]_1 \to K$ sending $01$ to $H_1$ and $12$ to $H_2$, then lifting it to a map $H_{12}\colon \Delta[2] \to K$, and finally restricting to $02^*H_{12}\colon \Delta[1] \to K$.
\end{proof}

\begin{lemma} \label{lemma:sing pizero}
 Let $X$ be a topological space. Then there is a bijection of sets $\pi_0(\Sing(X)) \cong \pi_0(X)$.
\end{lemma}

\begin{proof}
 By \cref{rem:sing explicit}, we have $\Sing(X)_0= \Hom(|\Delta[0]|,X)$ is the set of points of $X$, and $\Sing(X)_1 = \Hom(|\Delta[1]|,X)$ is precisely the set of paths in $X$. Hence, by \cref{def:pizero sset}, $\pi_0(\Sing(X))$ are the path-components of $X$, which is precisely $\pi_0(X)$.
\end{proof}

We can use this result to obtain the classical ``path-component'' decomposition. Here we consider $\pi_0(X)$ as a constant simplicial set, as defined in \cref{ex:constant simplicial set}.

\begin{lemma}
 Let $S$ be a simplicial set. There is a map of simplicial sets $\tau_0\colon S \to \pi_0(S)$.
\end{lemma}

\begin{proof}
 Let $(\tau_0)_n\colon S_n \to \pi_0(S)$ send a $\sigma\colon \Delta[n] \to S$ to the class $[0^*\sigma]$ in $\pi_0(S)$. As all $d_i$ and $s_i$ in the simplicial set $\pi_0(S)$ are identities (\cref{ex:constant simplicial set}), we only need to verify that for every $\sigma$ in $S_n$, $(\tau_0)_{n-1}(d_i\sigma) = (\tau_0)_{n}(\sigma)$ and $(\tau_0)_{n+1}(s_i\sigma) = (\tau_0)_{n}(\sigma)$. From \cref{rem:delta morphism}, we have
 \begin{itemize}[leftmargin=*]
  \item $(\tau_0)_{n+1}(s_i\sigma) = [0^*s_i\sigma] = [0^*\sigma] = (\tau_0)_{n}(\sigma)$
  \item $(\tau_0)_{n-1}(d_i\sigma) = [0^*d_i\sigma] = [0^*\sigma] = (\tau_0)_{n}(\sigma)$ if $i > 0$
  \item $(\tau_0)_{n-1}(d_0\sigma) = [0^*d_0\sigma] = [1^*\sigma] = [0^*\sigma] = (\tau_0)_{n}(\sigma)$, where we used the fact that $0^*\sigma$ and $1^*\sigma$ are connected by the path $01^*\sigma\colon \Delta[1] \to S$.
 \end{itemize}
 This proves naturality and shows we have a map of simplicial sets $\tau_0\colon S \to \pi_0(S)$.
\end{proof}

\begin{definition}[Path-component of a simplicial set]
 Let $S$ be a simplicial set. For $x$ in $S_0$, the path-component of $x$, denoted $S_x$, is the sub-simplicial set of $S$ given by the pullback (\cref{ex:pullback})
 \[
 \begin{tikzcd}
 S_x \arrow[r, "\iota_x"] \arrow[d] & S \arrow[d, "\tau_0"] \\
 \Delta[0] \arrow[r, "{[x]}"] & \pi_0(S)
 \end{tikzcd}.
 \]  
\end{definition}

\begin{remark}
  More explicitly $(S_x)_n$ is the pre-image of the class $[x]$ under the map $(\tau_0)_n\colon S_n \to \pi_0(S)$, meaning $S_x$ is indeed a sub-simplicial set of $S$.
\end{remark}

\begin{lemma} \label{lemma:path component decomposition}
 Let $S$ be a simplicial set, then $\coprod_{[x] \in \pi_0(S)}\iota_x \colon \coprod_{[x] \in \pi_0(S)} S_x \to S$ is an isomorphism of simplicial sets. Hence, every simplicial set decomposes as a disjoint union of its path-components.
\end{lemma}

\begin{proof}
 By definition each $S_x \hookrightarrow S$ is a sub-simplicial set of $S$. Moreover, for $[x] \neq [y]$ in $\pi_0(S)$, the intersection of $S_x$ and $S_y$ is empty. This means the map $\coprod_{[x] \in \pi_0(S)}\iota_x$ is injective. Finally, every simplex $\sigma$ in $S_n$ is contained in some $(S_{\sigma(0)})_n$, so $\coprod_{[x] \in \pi_0(S)} \iota_x$ is also surjective. 
\end{proof}

These results provide some first evidence that simplicial sets have a functioning homotopy theory and that $\Sing$ does indeed match the homotopy theory of topological spaces with Kan complexes. Ideally there should be a stronger result providing a precise match. While true, unfortunately, the general result is more involved, and we will only state it without proof. 

\begin{definition}[Homotopy category of CW-complexes]
  Let $\Ho(\Top^{\text{CW}})$ be the category with objects CW-complexes and morphisms homotopy classes of continuous maps.
\end{definition}

\begin{definition}[Homotopy category of Kan complexes]
  Let $\Ho(\Kan)$ be the category with objects Kan complexes and morphisms homotopy classes of maps.
\end{definition}

\cref{prop:sing homotopy} shows that there is an induced functor $\Sing\colon\Ho(\Top^{\text{CW}}) \to \Ho(\Kan)$. We now have the following major result.

\begin{theorem}[{\cite[Theorem 11.4]{goerssjardine2009simplicialhomotopytheory}}] \label{thm:sing equivalence}
 The functor $\Sing\colon \Ho(\Top^{\text{CW}}) \to \Ho(\Kan)$ is an equivalence of categories.
\end{theorem}

The proof of this statement is a central result of simplicial homotopy theory. Proving it goes far beyond the scope of this introduction, and we refer the interested reader to \cite{goerssjardine2009simplicialhomotopytheory} for details.

\begin{remark} \label{rem:two out of three}
  One particular implication of these results is that homotopy equivalences of Kan complexes satisfy the famous $2$-out-of-$3$ property, which states that if $f\colon X \to Y$ and $g\colon Y \to Z$ are maps of Kan complexes such that two of the three maps $f,g,g \circ f$ are homotopy equivalences, then so is the third. For a proof of this, we refer the reader to \cite[Proposition 9.2 CM2]{goerssjardine2009simplicialhomotopytheory}.
\end{remark}

\section{Category Theory \texorpdfstring{$\&$}{\&} Homotopy Theory: Two Paths towards Simplicial Spaces} \label{sec:twopaths}
In the previous section we studied two different aspects of the category of simplicial sets. We showed that one subcategory, namely $\sSet^{Seg}$ is equivalent to $\Cat$ and another subcategory, namely $\Kan$ has an equivalent homotopy theory to that of $\Top$. A higher category should combine the categorical aspects of $\Cat$ and the homotopical aspects of $\Top$. In particular both categories and topological spaces should be examples of higher categories. However, these two subcategories only intersect very trivially. 

 \begin{lemma} \label{lemma:groupoid kan complex}
  Let $\cC$ be a category. Then $N\cC$ is a Kan complex if and only if $\cC$ is a groupoid.
 \end{lemma}

\begin{proof}
  Let $\cC$ be a category such that $N\cC$ is a Kan complex. Let $f\colon x \to y$ be a morphism in $\cC$. We need to show that $f$ is an isomorphism. We now have the following bijections:
  \begin{align*}
    \Hom(\Lambda[2]_0, N\cC) & \cong \Hom(\Delta[1] \coprod^{0,0}_{\Delta[0]} \Delta[1],N\cC) & \text{\cref{lemma:horn colimit}}\\ 
    &  \cong \Hom(\Delta[1], N\cC) \times^{0^*,0^*}_{\Hom(\Delta[0], N\cC)} \Hom(\Delta[1], N\cC) & \text{\cref{prop:colimit limit}}\\ 
      & \cong N\cC_1 \times^{0^*,0^*}_{N\cC_0} N\cC_1 & \text{\cref{lemma:delta yoneda}}\\
      & \cong \Mor_\cC \times^{0^*,0^*}_{\Obj_\cC} \Mor_\cC & \text{\cref{rem:nerve explicit}} 
  \end{align*}
  So, the surjectivity of the Kan lifting condition in \cref{def:kan complex} means that for every pair of morphisms $(f,g)$ with the same source, there exists a pair $(f,h)$, such that $h\circ f = g$. Applying this to the pair $(f,\id_x)$, we get a left inverse $h$ of $f$. 

  Now, repeating the same argument with $\Lambda[2]_2$, we similarly deduce the existence of a right inverse of $f$. Hence, $f$ is an isomorphism and $\cC$ is a groupoid.

  The other direction is much more involved and so we refer the reader to \cite[\href{https://kerodon.net/tag/0037}{Proposition 0037}]{lurie2026kerodon}. 
\end{proof}  

So, how can we combine two different types of objects which intersect almost trivially in a larger setting? Let us try to motivate the solution via a very simple example.

\begin{example} \label{ex:combine functions}
  Let $P$ be the set of polynomial functions $p(x)\colon \bR \to \bR$ and $E$ the set of exponential functions $a \cdot b^x\colon \bR \to \bR$. Then $P, E$ are two subsets of functions from $\bR$ to $\bR$. How can we directly obtain a set of functions that includes $P,E$? Here is an easy trick, namely adding a new variable:
  \[ PE \coloneq \{ f\colon \bR \times \bR \to \bR \mid \forall x_0, y_0 \in \bR,  (f(x,y_0) \in P \text{ and } f(x_0,y) \in E)\}. \]
  Thus $PE$ includes both $P$ and $E$, as functions in $PE$ which only include the variable $x$, or $y$, respectively. 
\end{example}

This example suggests a solution to our predicament.
 
\begin{definition}[Simplicial space]
  A \emph{simplicial space} is a functor $X\colon \DD^{\op} \to \sSet$. The category of simplicial spaces is the functor category $\sSp \coloneq \Fun(\DD^{\op},\sSet)$.
\end{definition}

\begin{remark} \label{rem:simplicial space depiction}
  Following \cref{prop:functor category as currying}, we have an isomorphism of categories  
  \[\Fun(\DD^{\op},\sSet) = \Fun(\DD^{\op},\Fun(\DD^{\op},\Set)) \cong \Fun(\DD^{\op} \times \DD^{\op} , \Set).\]
  Thus we have two equivalent ways to depict a simplicial space $X$:
  \[
   \simplicialspacediagram
 \]
 Here $X_{nm}$ are sets, whereas $X_{n\bullet}$ are simplicial sets.  
 \end{remark}

By analogy with \cref{ex:combine functions} we need two ways to include simplicial sets into simplicial spaces, giving us a categorical and spatial direction.

\begin{definition}[Categorical inclusion] \label{def:inccat}
 The \emph{categorical inclusion}, denoted  $\inccat\colon \sSet \to \sSp$, is the functor $\pi_1^*\colon \Fun(\DD^{\op},\Set) \to \Fun(\DD^{\op} \times \DD^{\op},\Set)$ induced by the projection $\pi_1\colon \DD \times \DD \to \DD$ onto the first factor.
\end{definition}

\begin{intuition} \label{int:inccat}
 Given a simplicial set $S$, the diagram in \cref{rem:simplicial space depiction} of $\inccat S$ is the simplicial space with $(\inccat S)_{nk} = S_n$ and all vertical maps identities. 
\end{intuition}

\begin{definition}[Spatial inclusion] \label{def:incspace}
  The \emph{spatial inclusion}, denoted $\incspace\colon \sSet \to \sSp$, is the functor $\pi_2^*\colon \Fun(\DD^{\op},\Set) \to \Fun(\DD^{\op} \times \DD^{\op},\Set)$ induced by the projection $\pi_2\colon \DD \times \DD \to \DD$ onto the second factor.
\end{definition}

\begin{intuition} \label{int:incspace}
 Given a simplicial set $S$, the diagram in \cref{rem:simplicial space depiction} of $\incspace S$ is the simplicial space with $(\incspace S)_{nk} = S_k$ and all horizontal maps identities. 
\end{intuition}

\begin{remark}
  Based on \cref{int:inccat,int:incspace}, $\inccat$ is also called the \emph{horizontal embedding} (as only the horizontal direction varies) and $\incspace$ the \emph{vertical embedding} (as only the vertical direction varies).
\end{remark}

\begin{remark}
 For now the names categorical and spatial inclusion are just suggestive, but we will see in the coming sections (\cref{sec:reedy,sec:segal spaces}) that the naming conventions are indeed justified.  
\end{remark}

\begin{notation} \label{not:FDelta}
  Recall $\Delta[n]$ (\cref{def:standard simplex}), $\partial \Delta[n]$ (\cref{def:boundary}), and $\Sp[n]$ (\cref{def:spine}). We adopt the notational conventions:
  \begin{itemize}[leftmargin=*]
    \item $F(n) \coloneq \inccat(\Delta[n])$, $\partial F(n) \coloneq \inccat(\partial \Delta[n])$, $G(n) \coloneq \inccat(\Sp[n])$.
    \item $\Delta[l] \coloneq \incspace(\Delta[l])$, $\partial \Delta[l] \coloneq \incspace(\partial \Delta[l])$ (here we are identifying $\Delta[l]$ with its image under $\incspace$ to simplify notation).
    \item More generally, for a simplicial set $S$, we will again denote the simplicial space $\incspace S$ by $S$, and note this is consistent with the previous item.
  \end{itemize} 
\end{notation} 

We saw before that we can use the Yoneda lemma to characterize maps out of $\Delta[n]$. Now via the Yoneda lemma (\cref{thm:yoneda}), we have the following analogous result.

\begin{lemma} \label{lemma:FDelta yoneda}
  Let $X$ be a simplicial space and let $n,l \geq 0$. Then there are natural bijections 
  \[
   \Hom_{\sSp}(F(n) \times \Delta[l], X) \cong X_{nl}.
  \]
\end{lemma}

Given that we have more simplicial information now, we need to generalize our hom sets as well. For this next definition we recall the Yoneda embedding (\cref{def:yoneda embedding}), opposite category functor (\cref{lemma:opposite functor}), and the product functor (\cref{rem:products functor}).

\begin{definition}[Mapping space] \label{def:mapping space}
  Let $X,Y$ be two simplicial spaces. The \emph{mapping simplicial set} $\Map(X,Y)$ is defined as the simplicial set 
  \[ \DD^{\op} \xrightarrow{ \Yon^{\op} } \sSet^{\op} \xrightarrow{\incspace^{\op}}  \sSp^{\op} \xrightarrow{ (X \times -)^{\op}} \sSp^{\op} \xrightarrow{ \Hom(-,Y)} \Set\]
\end{definition}

\begin{remark} \label{rem:mapping space explicit}
  Very explicitly we can characterize the $n$-simplices of $\Map(X,Y)$ as $\Hom(X \times \Delta[n], Y)$, and the face and degeneracy maps are induced by the corresponding maps of $\Delta[n]$. In particular $\Map(X,Y)_0 = \Hom(X,Y)$, so the vertices of $\Map(X,Y)$ are precisely the maps from $X$ to $Y$.
\end{remark}

\begin{remark}
 By \cref{not:FDelta}, if $X,Y$ are simplicial sets, then $\Map(X,Y) = \Map(\incspace X, \incspace Y)$. 
\end{remark}

We now make several important observations about $\Map$ that mirror the behavior of $\Hom$. First we have a generalized form of the Yoneda lemma.

\begin{lemma} \label{lemma:map yoneda}
  Let $X$ be a simplicial space and let $n \geq 0$. Then there is a natural isomorphism of simplicial sets
  \[\Map(F(n), X) \cong X_n.\]
\end{lemma}

\begin{proof}
 For every $k \geq 0$, via \cref{lemma:FDelta yoneda} we have the following chain of bijections
 \[ \Map(F(n), X)_k = \Hom(F(n) \times \Delta[k], X) \cong X_{nk} = (X_n)_k,\]
 giving us the desired isomorphism between $\Map(F(n), X)$ and $X_n$.
\end{proof}

Beyond maps out of $F(n)$, which we can evaluate via the Yoneda lemma, there are two other cases we can easily compute.

\begin{example} \label{ex:initial}
  Let $X$ be a simplicial space, and $\emptyset$ denote the empty simplicial space (all levels are the empty set). Then $\Map(\emptyset, X) \cong \Delta[0]$.
\end{example}

\begin{example} \label{ex:final}
 Let $X$ be a simplicial space. Then $\Map(X, F(0)) \cong \Map(X,\Delta[0]) \cong \Delta[0]$.
\end{example}

Following \cref{rem:mapping space explicit}, the mapping space is suitably functorial and concretely we have the following observation.

\begin{lemma}
  Let $f\colon X \to Y$ be a map of simplicial spaces and let $Z$ be a simplicial space. The collection of maps $(f^*)_n\colon \Map(Y,Z)_n \to \Map(X,Z)_n$, given by 
  \[ (\sigma\colon Y \times \Delta[n] \to Z) \mapsto (\sigma \circ (f \times \id_{\Delta[n]})\colon X \times \Delta[n] \to Z),\]
  defines a map of simplicial sets $f^*\colon \Map(Y,Z) \to \Map(X,Z)$. 
\end{lemma}

\begin{proof}
 For a given $\delta\colon [n] \to [m]$, we need to show that the following diagram commutes
  \[
  \begin{tikzcd}
    \Map(Y,Z)_m \arrow[r, "{(f^*)_m}"] \arrow[d, "{\delta^*}"] & \Map(X,Z)_m \arrow[d, "{\delta^*}"] \\
    \Map(Y,Z)_n \arrow[r, "{(f^*)_n}"] & \Map(X,Z)_n
  \end{tikzcd}.
  \]  
  For a given $\sigma \colon Y \times \Delta[m] \to Z$, this unfolds to the evident equality
  \[(\sigma \circ (f \times \id_{\Delta[m]})) \circ (\id_X \times \delta) = \sigma \circ (\id_Y \times \delta) \circ (f \times \id_{\Delta[n]}),\]
  and hence we are done.
\end{proof}

Similarly, we have the following analogous lemma.

\begin{lemma}
  Let $f\colon X \to Y$ be a map of simplicial spaces and let $Z$ be a simplicial space. The collection of maps $(f_*)_n\colon \Map(Z,X)_n \to \Map(Z,Y)_n$, given by 
  \[ (\sigma\colon Z \times \Delta[n] \to X) \mapsto (f \circ \sigma \colon Z \times \Delta[n] \to Y),\]
  defines a map of simplicial sets $f_*\colon \Map(Z,X) \to \Map(Z,Y)$. 
\end{lemma}

We now move on to colimits and limits in $\sSp$. Similar to the situation for simplicial sets (\cref{def:colimit sset,def:limit sset}), the colimits and limits in $\sSp$ are computed individually for each $(k,n)$-simplex, and so we will not repeat those formulas here. However, we do have the following generalization of \cref{prop:colimit limit,prop:limit limit} for mapping spaces.

\begin{proposition} \label{prop:colimit limit map}
  Let 
  \[X_1 \xleftarrow{F_{1,2}} X_{1,2} \xrightarrow{G_{1,2}} X_{2} \xleftarrow{F_{2,3}} ... \xrightarrow{G_{n-1,n}} X_n\] 
  be a diagram of simplicial spaces, and let $C$ in $\sSp$ denote its colimit. Moreover, for a simplicial space $Y$, let $L$ denote the limit of the diagram of simplicial sets 
  \[\Map(X_1,Y) \xrightarrow{F_{1,2}^*} ... \xleftarrow{G_{n-1,n}^*} \Map(X_n,Y).\] 
  Then $L \cong \Map(C,Y)$.
\end{proposition}

\begin{proof}
  We want to show that $L_k \cong \Map(C,Y)_k = \Hom(C \times \Delta[k],Y)$. Here, we make the following elegant observation that follows from direct computation. $C \times \Delta[k]$ is the colimit of the diagram 
  \[X_1 \times \Delta[k] \xleftarrow{F_{1,2} \times \id_{\Delta[k]}} X_{1,2} \times \Delta[k] \xrightarrow{G_{1,2} \times \id_{\Delta[k]}} X_{2} \times \Delta[k] \xleftarrow{F_{2,3} \times \id_{\Delta[k]}} ... \xrightarrow{G_{n-1,n} \times \id_{\Delta[k]}} X_n \times \Delta[k]\] 
  So, by \cref{prop:colimit limit},
  $\Hom(C \times \Delta[k],Y)$ is bijective to the limit of the diagram of sets
  \[\Hom(X_1 \times \Delta[k],Y) \xrightarrow{F_{1,2}^*} ... \xleftarrow{G_{n-1,n}^*} \Hom(X_n \times \Delta[k],Y).\]
  This is, by definition, equal to the diagram of sets
  \[\Map(X_1,Y)_k \xrightarrow{F_{1,2}^*} ... \xleftarrow{G_{n-1,n}^*} \Map(X_n,Y)_k,\]
  which is how we defined $L_k$.
\end{proof}

Using the same argument, we also get the following generalization of \cref{prop:limit limit}.
\begin{proposition} \label{prop:limit limit map}
 Let 
  \[X_1 \xrightarrow{F_{1,2}} X_{1,2} \xleftarrow{G_{1,2}} X_{2} \xrightarrow{F_{2,3}} ... \xleftarrow{G_{n-1,n}} X_n\] 
  be a diagram of simplicial spaces, and let $L$ in $\sSp$ denote its limit. Moreover, for a simplicial space $Y$, let $L'$ denote the limit of the diagram of simplicial sets 
  \[\Map(Y,X_1) \xrightarrow{(F_{1,2})_*} ... \xleftarrow{(G_{n-1,n})_*} \Map(Y,X_n).\] 
  Then $L' \cong \Map(Y,L)$.
\end{proposition}

Let us see some examples.

\begin{example} \label{ex:map coproduct}
  Let $X$ be a simplicial space. Then, combining \cref{prop:colimit limit map} and \cref{lemma:map yoneda}, we get 
  \[\Map(F(0) \coprod ... \coprod F(0),X) \cong \Map(F(0),X) \times ... \times \Map(F(0),X) \cong X_0 \times ... \times X_0.\]
\end{example}

\begin{example} \label{ex:map boundary}
  In particular, analogous to \cref{lemma:boundary colimit}, we have $\partial F(1) \cong F(0) \coprod F(0)$, and so for any simplicial space $X$, we have $\Map(\partial F(1),X) \cong X_0 \times X_0$.
\end{example}

A general simplicial space is evidently not a suitable model for a higher category, as there is not even a notion of composition. We need to add some conditions.

\section{Interlude: Kan Fibrations as Dependent Spaces} \label{sec:kanfibrations}
In \cref{subsec:top to ss} we discussed Kan complexes as an analogue of topological spaces in the realm of simplicial sets. However, in the set-theoretic context we also use ``dependent sets'', otherwise known as functions. For example, we have the following result, the proof of which we leave as an exercise.

\begin{lemma}
  Assume we have the following commutative diagram of sets
  \[
  \begin{tikzcd}
   Y \arrow[rr, "f"] \arrow[dr, "p"'] & & Z \arrow[dl, "q"] \\ 
   & X 
  \end{tikzcd}.
  \] 
  Then $f$ is a bijection if and only if for every $x \in X$ the induced map on pre-images $f_x\colon p^{-1}(x) \to q^{-1}(x)$ is a bijection.
\end{lemma}

Intuitively this result is stating that a ``global bijection'' between the sets $Y$ and $Z$ that sits over $X$ is determined by ``local bijections'' on the fibers. Ideally, in the context of spaces we should have similar statements, where bijection is replaced by homotopy equivalence. However, unfortunately, this will fail without some assumption. Here we use the groupoid $I(1)$, reviewed in \cref{ex:ione}.

\begin{lemma} \label{lemma:i contractible}
 The map $N(0)\colon N[0] \to N(I(1))$ is a homotopy equivalence of Kan complexes.
\end{lemma}

\begin{proof}
  Following \cref{lemma:groupoid kan complex}, $\Delta[0] = N[0]$ and $N(I(1))$ are both Kan complexes. Let $t\colon I(1) \to [0]$ be the unique functor. By uniqueness, we must have $t \circ 0 = \id_{[0]}$. Moreover, let $H\colon [1] \times I(1) \to I(1)$ be the functor defined as follows:
    \[H(0,0) = 0, H(1,0) = 0, H(0,1) = 0, H(1,1) = 1,\]
  and for morphisms the choices are unique. Now, by \cref{prop:nerve product}, $N([1] \times I(1)) \cong N([1]) \times N(I(1)) = \Delta[1] \times N(I(1))$. Then $NH \colon \Delta[1] \times N(I(1)) \to N(I(1))$ is a homotopy between $0 \circ t$ and $\id_{N(I(1))}$, as $NH|_{\{ 0\} \times N(I(1))} = N(0 \circ t)$ and $NH|_{\{ 1\} \times N(I(1))} = N(\id_{I(1)})$. 
\end{proof}

\begin{example} \label{ex:preimage fail}
   The functor $0\colon [0] \to I(1)$ induces a commutative diagram of Kan complexes
  \[ 
   \begin{tikzcd}
    \Delta[0] \arrow[rr, "N(0)"] \arrow[dr, "N(0)"'] & & N(I(1)) \arrow[dl, "\id"] \\
    & N(I(1))
   \end{tikzcd}.
  \]
  By \cref{lemma:i contractible}, the top map is a homotopy equivalence. However, the pre-image over $1$ is the map $\emptyset \to \Delta[0]$, which is not a homotopy equivalence, and in fact admits no map in the other direction.
\end{example}

This is not the only situation where a suitable notion of dependency is needed. Here is another simple example, which we will leave as an exercise.

\begin{lemma} \label{lemma:pullback invariance}
  Assume we have the following commutative diagram of sets
  \[ 
  \begin{tikzcd}
   A_1 \arrow[r] \arrow[d, "\cong", "f_1"'] & A_{1,2} \arrow[d, "\cong", "f_{1,2}"'] & ... \arrow[r] \arrow[l] & A_{n-1,n} \arrow[d, "\cong", "f_{n-1,n}"'] & A_n \arrow[d, "\cong", "f_n"'] \arrow[l] \\ 
   B_1 \arrow[r] & B_{1,2} & ... \arrow[r] \arrow[l] & B_{n-1,n} & B_n \arrow[l]
  \end{tikzcd},
  \]
  such that all vertical maps are bijections. Then the induced map (\cref{lemma:limit map}) is a bijection.
\end{lemma}

Again, we would like to have a similar statement in the context of Kan complexes, where bijections are replaced by homotopy equivalences. However, once more this will fail without some assumption.

\begin{example} \label{ex:pullback invariance fail}
  Take the diagram of Kan complexes, where, by \cref{lemma:i contractible}, all vertical maps are homotopy equivalences,
  \[ 
  \begin{tikzcd}
   \Delta[0] \arrow[r, "0"] \arrow[d, "="] & N(I(1)) \arrow[d, "\simeq"] & \Delta[0] \arrow[d, "="] \arrow[l, "1"'] \\ 
   \Delta[0] \arrow[r] & \Delta[0] & \Delta[0] \arrow[l]
  \end{tikzcd}
  \]
  For the top pullback, we have:
  \[(\Delta[0] \times_{N(I(1))}^{0,1} \Delta[0])_0 = \{(0,0) \in \Delta[0]_0^2 \mid 0 = 1 \} = \emptyset.\]
  So, the induced map on pullbacks is the map $\emptyset \to \Delta[0]$, which again is not a homotopy equivalence.
\end{example}

We hence need a proper notion of ``dependent spaces'' or ``fibrations''. Fortunately, there is no need to reinvent the wheel, and we can use the well-established notion of Kan fibrations.

\begin{definition}[Kan fibration]
  A map of simplicial sets $Y \to X$ is a \emph{Kan fibration} if for every $n \geq 1, 0 \leq i \leq n$ the following diagram lifts
  \[
    \begin{tikzcd}[row sep=1.2cm, column sep=1.2cm]
     \Lambda[n]_i  \arrow[r] \arrow[d] & Y \arrow[d, twoheadrightarrow] \\
     \Delta[n] \arrow[ur, dashed] \arrow[r] & X
    \end{tikzcd}
  \]
\end{definition}

\begin{notation}
  In diagrams we denote Kan fibrations by $\twoheadrightarrow$.
\end{notation}

\begin{remark} \label{rem:kan fibration kan complex}
  Kan fibrations are indeed relative versions of Kan complexes, as Kan fibrations $Y \to \Delta[0]$ precisely correspond to a choice of Kan complex $Y$. Here we used the fact that the map to $\Delta[0]$ is unique (\cref{lemma:delta 0}).
\end{remark}

Kan fibrations are determined via lifting conditions, and so compose.

\begin{lemma} \label{lemma:kan fibration compose}
  The composition of two Kan fibrations is a Kan fibration.
\end{lemma}

\begin{proof}
 Let $p\colon Z \to Y$ and $q\colon Y \to X$ be two Kan fibrations. We want to show that $q \circ p\colon Z \to X$ is a Kan fibration. Let $n \geq 1, 0 \leq i \leq n$, and consider the following diagram:
  \[
    \begin{tikzcd}
     \Lambda[n]_i  \arrow[r] \arrow[dd, hookrightarrow] & Z \arrow[d, "p", twoheadrightarrow] \\
      & Y \arrow[d, "q", twoheadrightarrow] \\
     \Delta[n] \arrow[ur, dashed] \arrow[uur, dashed] \arrow[r] & X
    \end{tikzcd}.
  \]
  As both $p$ and $q$ are Kan fibrations, we can construct the lift in two steps, first lifting along $q$ and then along $p$.
\end{proof}

\begin{lemma} \label{lemma:kan fibration pullback}
  Let the following diagram be a pullback square of simplicial sets (\cref{def:pullback square})
  \[
  \begin{tikzcd}
    W \arrow[r] \arrow[d, "q"'] \arrow[dr, phantom, "\ulcorner", very near start] & Y \arrow[d, "p", twoheadrightarrow] \\
    Z \arrow[r] & X
  \end{tikzcd}.
  \]
  If $p$ is a Kan fibration, then $q$ is also a Kan fibration.
\end{lemma}

\begin{proof}
  Let $n \geq 1, 0 \leq i \leq n$, and consider the following diagram:
  \[
    \begin{tikzcd}
     \Lambda[n]_i  \arrow[r, "h"] \arrow[d, hookrightarrow] & W \arrow[d, "q"]  \arrow[r, "k"] & Y \arrow[d, "p", twoheadrightarrow] \\
     \Delta[n] \arrow[ur, dashed]  \arrow[r, "f"'] & Z \arrow[r, "g"'] & X
    \end{tikzcd}.
  \]
  By \cref{prop:limit limit}, such a lift $\Delta[n] \to W$ corresponds to a commutative diagram of the form 
  \[
    \begin{tikzcd}
      \Lambda[n]_i \arrow[d] \arrow[rr, "kh"] & & Y \arrow[d, "p", twoheadrightarrow] \\ 
       \Delta[n] \arrow[r, "f"'] \arrow[urr, dashed] & Z \arrow[r, "g"'] & X
    \end{tikzcd}.
  \]
  However, this lift exists as $p$ is a Kan fibration, and hence we are done.
\end{proof}

\begin{remark}
  The astute reader might have noticed that the proofs of \cref{lemma:kan fibration compose,lemma:kan fibration pullback} only depend on the fact that Kan fibrations are characterized via a lifting condition. This is indeed the case, and so the same proofs work for any class of maps that is defined analogously. See \cite[Section 11.1]{riehl2014categoricalhomotopytheory} for a more abstract approach.
\end{remark}

The pullback property also has some additional useful implications about Kan fibrations.

\begin{lemma} \label{lemma:kan fibration factorization}
 Let $p\colon Y \to X$ be a Kan fibration. Assume $Y \xrightarrow{ \ q \ } A \xrightarrow{ \ i \ } X$ is a factorization of $p$, such that $i$ is an inclusion of simplicial sets. Then $q$ is a Kan fibration.
\end{lemma}

\begin{proof}
  The inclusion assumption implies that the following diagram is a pullback square
  \[
    \begin{tikzcd}
      Y \arrow[r, equal] \arrow[d, "q"', twoheadrightarrow] \arrow[dr, phantom, "\ulcorner", very near start] & Y \arrow[d, "p", twoheadrightarrow] \\
      A \arrow[r, "i", hookrightarrow] & X
    \end{tikzcd},
  \]
  and so the result follows from \cref{lemma:kan fibration pullback} and that $p$ is a Kan fibration.
\end{proof}

We now have the following powerful technical result about pullbacks and Kan fibrations, analogous to \cref{lemma:pullback invariance}. Unfortunately, a proof requires developing some theory and goes beyond our aims. See \cite[Proposition 13.3.9]{hirschhorn2003modelcategories} for a proof of the case $n=2$.

\begin{lemma} \label{lemma:homotopy pullback}
 Assume we have the following diagram of Kan complexes,
   \[ 
  \begin{tikzcd}
   A_1 \arrow[r, twoheadrightarrow] \arrow[d, "\simeq", "f_1"'] & A_{1,2} \arrow[d, "\simeq", "f_{1,2}"'] & ... \arrow[r, twoheadrightarrow] \arrow[l, twoheadrightarrow] & A_{n-1,n} \arrow[d, "\simeq", "f_{n-1,n}"'] & A_n \arrow[d, "\simeq", "f_n"'] \arrow[l, twoheadrightarrow] \\ 
   B_1 \arrow[r, twoheadrightarrow] & B_{1,2} & ... \arrow[r, twoheadrightarrow] \arrow[l, twoheadrightarrow] & B_{n-1,n} & B_n \arrow[l, twoheadrightarrow]
  \end{tikzcd},
  \]
 where the horizontal maps are Kan fibrations and the vertical maps are homotopy equivalences. Then the induced map on limits is also a homotopy equivalence.
\end{lemma}

\begin{remark} \label{rem:homotopy pullback}
 In fact \cite[Proposition 13.3.9]{hirschhorn2003modelcategories} states that if $n =2$ in \cref{lemma:homotopy pullback}, then it suffices for $A_1 \to A_{1,2}$ and $B_1 \to B_{1,2}$ to be Kan fibrations for the result to hold.
\end{remark}

\begin{lemma}[{\cite[\href{https://kerodon.net/tag/00XC}{Proposition 00XC}]{lurie2026kerodon}}] \label{lemma:homotopy pullback preimage}
  Assume we have the following commutative diagram of Kan complexes
  \[
  \begin{tikzcd}
   Y \arrow[rr, "f"] \arrow[dr, "p"', twoheadrightarrow] & & Z \arrow[dl, "q", twoheadrightarrow] \\ 
   & X 
  \end{tikzcd},
  \] 
  where the maps $p$ and $q$ are Kan fibrations. Then $f$ is a homotopy equivalence if and only if for every $x$ in $X$ the induced map on pre-images (which we also call fibers) $f_x\colon p^{-1}(x) \to q^{-1}(x)$ is a homotopy equivalence.
\end{lemma}

We now study how mapping spaces (\cref{def:mapping space}) interact with Kan fibrations. This is in fact one of the most important and powerful results about Kan fibrations. Unfortunately, the proof is very technical, so we refer the reader to \cite[Corollary 5.3]{goerssjardine2009simplicialhomotopytheory} as a suitable source.

\begin{proposition} \label{prop:kan fibration mapping space}
  Let $X$ be a Kan complex, and $i\colon A \hookrightarrow B$ be an inclusion of simplicial sets (a sub-simplicial set). Then the induced map of simplicial sets
  \[i^*\colon \Map(B,X) \twoheadrightarrow \Map(A,X)\]
  is a Kan fibration.
\end{proposition}

There are particularly relevant examples of this result.

\begin{example} \label{ex:map kan complex}
  Let $X$ be a Kan complex, and $A$ be a simplicial set. Then the map $\emptyset \to A$ induces a Kan fibration $\Map(A,X) \twoheadrightarrow \Map(\emptyset,X) = \Delta[0]$ (\cref{ex:initial}), which means that $\Map(A,X)$ is a Kan complex (\cref{rem:kan fibration kan complex}).
\end{example}

We can use this example to get some new equivalences of Kan complexes.

\begin{example} \label{ex:map kan complex equivalence}
  Let $K$ be a Kan complex. Then the maps $\Delta[0] \xrightarrow{ 0 } \Delta[n] \to \Delta[0]$ induce homotopy equivalences of Kan complexes $\Map(\Delta[0], K) \simeq \Map(\Delta[n], K)$.
\end{example}

For the next example, we introduce new notions.

\begin{definition}[Path space] \label{def:path space}
  Let $K$ be a Kan complex and let $x,y$ be two points in $K$, corresponding to two maps $x,y \colon \Delta[0] \to K$. We define the \emph{path space} between $x$ and $y$ as the following pullback
  \[
  \begin{tikzcd}
   \Path_K(x,y) \arrow[dr, phantom, "\ulcorner", very near start] \arrow[r] \arrow[d, twoheadrightarrow] & \Map(\Delta[1],K) \arrow[d, "{(0^* , 1^*)}", twoheadrightarrow] \\
   \Delta[0] \arrow[r, "{(x , y)}"'] & \Map(\partial \Delta[1],K)
  \end{tikzcd},
  \]
  where for the bottom map we use the explicit characterization of $\Map(\partial \Delta[1],K)_0 \cong K_0 \times K_0$, relying on the computation $\partial \Delta[1] \cong \Delta[0] \coprod \Delta[0]$ (\cref{lemma:boundary colimit}) and \cref{prop:colimit limit map}.
 \end{definition}

 \begin{lemma}
  Let $K$ be a Kan complex and let $x,y$ be two points in $K$. Then the path space $\Path_K(x,y)$ is a Kan complex.
 \end{lemma}

 \begin{proof}
  Applying \cref{prop:kan fibration mapping space} to the inclusion $\partial \Delta[1] \hookrightarrow \Delta[1]$, it follows that $\Map(\Delta[1],K) \twoheadrightarrow \Map(\partial \Delta[1],K)$ is a Kan fibration. Thus, by \cref{lemma:kan fibration pullback}, the pullback $\Path_K(x,y) \to \Delta[0]$ is a Kan fibration, which, by \cref{rem:kan fibration kan complex}, means that $\Path_K(x,y)$ is a Kan complex.
 \end{proof}
 
Let us see one source of many Kan fibrations.

\begin{definition}[Inclusion of path-components] \label{def:inclusion path components}
  A map $i\colon A \to B$ of Kan complexes is an \emph{inclusion of path-components} if the following diagram is a pullback
  \[
  \begin{tikzcd}
    A \arrow[r, "\tau_0"] \arrow[d, hookrightarrow, "i"] \arrow[dr, phantom, "\ulcorner", very near start] & \pi_0(A) \arrow[d, hookrightarrow, "\pi_0(i)"] \\
    B \arrow[r, "\tau_0"] & \pi_0(B)
  \end{tikzcd}.
  \]
\end{definition}

\begin{intuition} \label{intuition:inclusion path components}
 Recall, from \cref{lemma:path component decomposition}, that we have a path-component decomposition $B \cong \coprod_{b \in \pi_0(B)}B_b$, for every simplicial set $B$. From this perspective we can very intuitively say that an inclusion of path-components is simply a choice of subset $A' \subseteq \pi_0(B)$ and $A = \coprod_{b\in A'} B_b$, at which point $i\colon A \to B$ is simply $\coprod_{b \in A'} B_b \hookrightarrow \coprod_{b \in \pi_0(B)} B_b$, meaning it is indeed an inclusion of path-components.
\end{intuition}

We now have the following lemmas and results.

\begin{lemma} \label{lemma:delta lambda path component}
 Let $n \geq 1, 0 \leq l \leq n$. Then $\pi_0(\Delta[n]) \cong \pi_0(\Lambda[n]_l) \cong \{[0]\}$.
\end{lemma}

\begin{proof}
 Any two vertices in $\Lambda[n]_l$ are connected by a chain of $1$-simplices, indeed in the cases $n \geq 3$, any two vertices are directly connected, and in the case of $n = 2$, there is always a chain of length 2 (\cref{ex:lambda two horns}). We have a similar argument for $\Delta[n]$. Hence, $\pi_0(\Delta[n]) \cong \pi_0(\Lambda[n]_l) \cong \{[0]\}$.
\end{proof}

We now have the following result about maps into constant simplicial sets, as defined in \cref{def:constant simplicial set}.

\begin{lemma} \label{lemma:constant boundary sset}
 Let $S$ be a constant simplicial set. Then for $n \geq 1, 0 \leq l \leq n$, every map $\Delta[n] \to S$ and $\Lambda[n]_l \to S$ is constant.
\end{lemma}

\begin{proof}
  We prove the case for $\Lambda[n]_l$ and the case for $\Delta[n]$ is similar (and even easier). A map $f\colon \Lambda[n]_l \to S$ induces a map on path-components $\pi_0(\Lambda[n]_l) \to \pi_0(S)$ (\cref{def:pizero sset}). As $S$ is constant, $\pi_0(S) = S_0$. Hence, by \cref{lemma:delta lambda path component}, $\pi_0(f)\colon\pi_0(\Lambda[n]_l) \to S_0$ is necessarily constant with image $s$ in $S_0$. 

  Now, let $\sigma$ in $(\Lambda[n]_l)_k$. Then, by the previous paragraph $f_0(0^*\sigma) = s$. Now, by naturality $0^*f_k(\sigma) = f_0(0^*\sigma) = s$, but $0^*\colon S_k \to S_0$ is a bijection, so $f_k(\sigma) = s$ and we are done.
\end{proof}

\begin{lemma} \label{lemma:constant kan fibration}
  Let $S,T$ be constant simplicial sets. Then any map $f\colon S \to T$ is a Kan fibration.
\end{lemma}

\begin{proof}
 We need to show that for every $n \geq 1, 0 \leq k \leq n$, the following diagram admits a lift
  \[
    \begin{tikzcd}[row sep=1.2cm, column sep=1.2cm]
     \Lambda[n]_k  \arrow[r] \arrow[d] & S \arrow[d, "f"] \\ 
     \Delta[n]  \arrow[r] & T
    \end{tikzcd}.
  \]
  By \cref{lemma:constant boundary sset}, this diagram will factor as follows
  \[
    \begin{tikzcd}[row sep=1.2cm, column sep=1.2cm]
     \Lambda[n]_k  \arrow[r] \arrow[d] & \Delta[0] \arrow[d, equal] \arrow[r, "s"] & S \arrow[d, "f"] \\ 
     \Delta[n] \arrow[r] & \Delta[0] \arrow[r, "t"] & T
    \end{tikzcd}.
  \]
  Now, the desired lift is given by the map $\Delta[n] \to \Delta[0] \xrightarrow{s} S$, and so we are done.
\end{proof}

\begin{intuition}
  A Kan fibration is seen as a ``dependent space'', i.e.~a space that varies over another space. So, if the spatial structure is trivial, then intuitively one would expect the Kan fibration condition to follow automatically. This is the content of \cref{lemma:constant kan fibration}.
\end{intuition}

\begin{proposition} \label{prop:inclusion path component kan fibration}
  Let $i\colon A \hookrightarrow B$ be an inclusion of path-components. Then $i$ is a Kan fibration.
\end{proposition}

\begin{proof}
 By \cref{lemma:constant kan fibration}, the map $\pi_0(i)\colon \pi_0(A) \hookrightarrow \pi_0(B)$ is a Kan fibration. By \cref{def:inclusion path components} $i$ is a pullback of $\pi_0(i)$. So, by \cref{lemma:kan fibration pullback}, the map $i\colon A \to B$ is a Kan fibration. 
\end{proof}

\begin{intuition}
  Following \cref{intuition:inclusion path components}, an inclusion of path-components $A \to B$ can be seen as a map of the form $\coprod_{i \in I_0} B_i \to \coprod_{i \in I} B_i$. As $\Lambda[n]_k, \Delta[n]$ are connected, for any commutative diagram of the form
  \[  
  \begin{tikzcd}
   \Lambda[n]_k \arrow[r] \arrow[d, hookrightarrow] & A \arrow[d, hookrightarrow, "i"] \\
   \Delta[n] \arrow[r] & B
  \end{tikzcd},
  \]
  the map $\Lambda[n]_k \to A$ needs to factor through a specific path-component $B_i$ of $A$. Moreover, by commutativity, the image of $\Delta[n] \to B$ also needs to factor through $B_i$. Corestricting the right-hand side to $B_i$, we hence only need to find a lift for the $\id_{B_i}\colon B_i \to B_i$, which of course always exists.
\end{intuition}

\section{Reedy Fibrancy} \label{sec:reedy}
As explained in \cref{sec:twopaths}, our simplicial spaces should behave like ``spaces'' (i.e.~Kan complexes) in what we already called the spatial (or ``vertical'') direction. Naively, we might hence expect to impose a level-wise Kan complex condition. However, later on we aim to impose additional conditions. For example, as we showed in \cref{prop:nerve Segal}, the nerve $N\cC$ of a category $\cC$ satisfies the Segal condition, which, as we saw in \cref{prop:Segal condition pullback}, is in particular a statement about the pullback $N\cC_1 \times_{N\cC_0} N\cC_1$.  We want an analogous statement in the context of simplicial spaces, meaning we want to analyze pullbacks of Kan complexes $X_1 \times_{X_0} X_1$. However, Kan complexes operate in the realm of homotopy theory and so all statements about them, including pullbacks, should be invariant under homotopy equivalences. Here we can benefit from the insight of \cref{sec:kanfibrations}. 

As we saw in \cref{lemma:homotopy pullback}, if the maps $0^*,1^*\colon X_1 \to X_0$ are Kan fibrations, then the pullback behaves correctly. This is one of the situations where the various spaces constructing a simplicial space need to be dependent (in the sense of \cref{sec:kanfibrations}). For that reason we now introduce the Reedy fibrancy condition.

\begin{definition}[Reedy fibrancy] \label{def:reedy}
  A simplicial space $X$ is \emph{Reedy fibrant} if for every $n \geq 0, k \geq 1, 0 \leq i \leq k$, the following diagram lifts 
  \[ 
  \begin{tikzcd}[row sep=1.2cm, column sep=1.2cm]
   \partial F(n) \times \Delta[k] \coprod_{\partial F(n) \times \Lambda[k]_i} F(n) \times \Lambda[k]_i \arrow[r] \arrow[d] & X \\
   F(n) \times \Delta[k] \arrow[ur, dashed]
  \end{tikzcd} 
  \]
\end{definition}

The Reedy fibrancy condition has many powerful implications and equivalent characterizations. However, proving that in full generality requires some more advanced machinery. We refer the reader to \cite[Section 2.5]{rezk2001css} for a review and to \cite[Proposition 9.3.1, Theorem 15.3.4]{hirschhorn2003modelcategories} for more detailed proofs.
\begin{lemma}\label{lemma:reedy}
 Let $X$ be a simplicial space. Then the following are equivalent:
  \begin{enumerate}
   \item $X$ is Reedy fibrant.
   \item For every $n \geq 0$, the induced map of mapping spaces
   \[\Map(F(n), X) \twoheadrightarrow \Map(\partial F(n), X)\]
   is a Kan fibration.
   \item For every inclusion of simplicial spaces $A \hookrightarrow B$, the induced map of mapping spaces
   \[\Map(B, X) \twoheadrightarrow \Map(A, X)\]
   is a Kan fibration.
  \end{enumerate}
\end{lemma}

We can use the Reedy fibrancy condition to get several important dependent Kan complexes, in the sense of \cref{sec:kanfibrations}. Let us see various manifestations of this dependency condition.

\begin{example} \label{ex:reedy fibrant zero}
  The map $\emptyset \to F(n)$ is an inclusion of simplicial spaces. So, for a given Reedy fibrant simplicial space $X$, the induced map $\Map(F(n), X) \twoheadrightarrow \Map(\emptyset, X)$ is a Kan fibration. By \cref{lemma:map yoneda,ex:initial}, this is isomorphic to the map $X_n \twoheadrightarrow \Delta[0]$, which, by \cref{rem:kan fibration kan complex}, means that $X_n$ is a Kan complex.
\end{example}

\begin{example} \label{ex:reedy fibrant one}
  The map $\partial F(1) \to F(1)$ is an inclusion of simplicial spaces. So, for a given Reedy fibrant simplicial space $X$, the induced map $\Map(F(1), X) \twoheadrightarrow \Map(\partial F(1), X)$ is a Kan fibration. By \cref{lemma:map yoneda,ex:map boundary}, this is isomorphic to the map $X_1 \twoheadrightarrow X_0 \times X_0$. Now, by \cref{ex:reedy fibrant zero}, $X_0 \to \Delta[0]$ is a Kan fibration, and so by \cref{lemma:kan fibration pullback}, both maps $\pi_1, \pi_2\colon X_0 \times X_0 \to X_0$ are Kan fibrations. As Kan fibrations compose (\cref{lemma:kan fibration compose}), this means $\pi_1 \circ (0^*,1^*) = 0^*$ and $\pi_2 \circ (0^*,1^*) = 1^*$ are Kan fibrations as well.
\end{example}

\begin{example} \label{ex:reedy fibrant n}
 Let $n \geq 2$ and let $0 + 1 + ... + n\colon \coprod_{0 \leq i \leq n} F(0) \to F(n)$ be the map that sends the $i$-th copy of $F(0)$ to the vertex $i$ in $F(n)_0 = \{0, ..., n\}$. This is an inclusion of simplicial spaces. So, for a given Reedy fibrant simplicial space $X$, the induced map $\Map(F(n), X) \twoheadrightarrow \Map(\coprod_{0 \leq i \leq n} F(0), X)$ is a Kan fibration. Again, by \cref{lemma:map yoneda,ex:map coproduct}, it hence follows that $X_n \to (X_0)^{n+1}$ is a Kan fibration. 
\end{example}

These observations show how Reedy fibrancy helps realize our goal of having a collection of Kan complexes with a well-behaved notion of pullback. Let us end with some examples.

\begin{lemma} \label{lemma:fn discrete}
  Let $A$ be a simplicial set. Then $\Map(F(n), \inccat(A))$ is a constant simplicial set on the set $\Hom(F(n), \inccat(A))$.
\end{lemma}

\begin{proof}
 We have 
 {\small
 \[\Map(F(n), \inccat(A))_l = \Hom(F(n) \times \Delta[l], \inccat(A)) \cong \inccat(A)_{nl} = \inccat(A)_{n0} \cong  \Hom(F(n), \inccat(A)),\]
 }
 finishing the proof.
\end{proof}

\begin{lemma} \label{lemma:boundary fn discrete}
 Let $A$ be a simplicial set. Then $\Map(\partial F(n), \inccat(A))$ is a constant simplicial set on the set $\Hom(\partial F(n), \inccat(A))$.
\end{lemma}

\begin{proof}
  By \cref{lemma:map yoneda}, we have the following commutative diagram
  \[
  \begin{tikzcd}
    \Map(\partial F(n), \inccat(A))_0 \arrow[r, hookrightarrow] \arrow[d, "\cong"] & \Map(\partial F(n), \inccat(A))_l \arrow[d, "\cong"] \\
    \Hom_{\sSp}(\partial F(n), \inccat(A)) \arrow[r, hookrightarrow, "\pi_1^*"] & \Hom_{\sSp}(\partial F(n) \times \Delta[l]  , \inccat(A))
  \end{tikzcd},
  \]
  where the vertical maps are bijections and the horizontal maps are injections. We want to show the top map is a surjection. Because of the bijections this is equivalent to the bottom map being a surjection. Unwinding definitions, this means every map $\sigma\colon \partial F(n) \times \Delta[l] \to \inccat(A)$ factors through $\pi_1\colon \partial F(n) \times \Delta[l] \to \partial F(n)$.

  Fix a $k \geq 0$. Then $\partial F(n)_k = \Hom([k],[n])^{\text{non-surj}}$ is a constant simplicial set (\cref{not:FDelta,def:boundary}). Moreover, building on \cref{not:FDelta}, $\Delta[l]_k = (\incspace \Delta[l])_k = \Delta[l]$ is a connected simplicial set (\cref{lemma:delta lambda path component}). So, the decomposition into path-components (\cref{lemma:path component decomposition}) of $\partial F(n)_k \times \Delta[l]$ is precisely $\coprod_{\Hom([k],[n])^{\text{non-surj}}} \Delta[l]$. Similarly, $\partial F(n)_k$ admits a decomposition into path-components given by $\coprod_{\Hom([k],[n])^{\text{non-surj}}} \Delta[0]$.

  With these observations at hand, we can now restate our goal as trying to show that for every map $\sigma\colon \partial F(n) \times \Delta[l] \to \inccat(A)$, the map $\sigma_k \colon \coprod_{\Hom([k],[n])^{\text{non-surj}}} \Delta[l] \to \inccat(A)$ admits a factorization of the form
  \[\coprod_{\Hom([k],[n])^{\text{non-surj}}} \Delta[l] \to \coprod_{\Hom([k],[n])^{\text{non-surj}}} \Delta[0] \to A_k.\]
  This is equivalent to showing that for every non-surjective map $f\colon [k] \to [n]$ the restriction of $\sigma$ on $\{f\} \times \Delta[l] \cong \Delta[l]$ factors through $\{f\} \times \Delta[0]$. This was proven in \cref{lemma:constant boundary sset}, and hence we are done.
\end{proof}

\begin{proposition} \label{prop:inccat reedy fibrant}
  Let $A$ be a simplicial set. Then $\inccat(A)$ is Reedy fibrant.
\end{proposition}

\begin{proof}
  By \cref{lemma:reedy}, it suffices to show that for every $n \geq 0$, the map 
  \[\Map(F(n), \inccat(A)) \twoheadrightarrow \Map(\partial F(n), \inccat(A))\] 
  is a Kan fibration. By \cref{lemma:fn discrete,lemma:boundary fn discrete} both domain and codomain are constant simplicial sets, and so the desired result follows from \cref{lemma:constant kan fibration}.
\end{proof}

\begin{example}
  Combining \cref{not:FDelta} and \cref{prop:inccat reedy fibrant} it follows that the simplicial spaces $F(n), \partial F(n)$ are Reedy fibrant for every $n \geq 0$, and  $G(n)$ is Reedy fibrant for $n \geq 2$.
\end{example}

\begin{example} \label{ex:delta not reedy fibrant}
  The simplicial space $\Delta[1]$, which is our notation for $\incspace(\Delta[1])$, is not Reedy fibrant. Indeed, following \cref{ex:reedy fibrant zero}, the simplicial set $\incspace(\Delta[1])_0 = \Delta[1] = N([1])$ would need to be a Kan complex, which, following \cref{lemma:groupoid kan complex}, would require $[1]$ to be a groupoid, which is not the case (\cref{ex:groupoids}).
\end{example}

\begin{example} \label{ex:ione not reedy fibrant}
  Recall from \cref{lemma:groupoid kan complex} that $N(I(1))$ is in fact a Kan complex. Yet, the simplicial space $N(I(1))$, which is our notation for $\incspace(N(I(1)))$, is not Reedy fibrant. Indeed, following \cref{ex:reedy fibrant one}, the diagonal map $\Delta\colon N(I(1)) \to N(I(1)) \times N(I(1))$ needs to be a Kan fibration. However, recalling that $\Lambda[1]_0 \cong \Delta[0]$, by \cref{ex:lambda one horns}, the diagram
  \[
  \begin{tikzcd}
   \Lambda[1]_0 \arrow[r, "0"] \arrow[d] & N(I(1)) \arrow[d, "\Delta"] \\
   \Delta[1] \arrow[r, "{(01,00)}"] & N(I(1)) \times N(I(1))
  \end{tikzcd}
  \]
  does not admit a lift. Indeed, at the level of $0$-simplices, such a lift would imply that $(1,0)$ in $(N(I(1)) \times N(I(1)))_0$ is in the image of the diagonal map, which is evidently not the case.
\end{example}

\section{Interlude: Uniqueness vs. Contractibility} \label{sec:contractibility}
Many properties in category theory (composition, inverses, ... ) involve uniqueness. Hence, before we proceed with our study of higher categories it is instructive to review the concept of ``uniqueness'' in a homotopical context. 

Recall that in set theory, for a given set $X$ we can define formulas $\varphi(x)$ on $X$. An easy example is the formula $n > 3$ on the set $\bN$. Given such a formula $\varphi$ on $X$, we can construct the subset of $X$ consisting of all elements that satisfy $\varphi$, i.e.~the subset $\{x \in X \mid \varphi(x)\}$. In our example this would simply be $\{n \in \bN \mid n > 3\} = \{4,5,...\}$. We now have the following classical definition.

\begin{definition}
 Let $X$ be a set and let $\varphi$ be a property on $X$. We say $x_0 \in X$ satisfies $\varphi$ uniquely if $x_0$ is the only element in $\{x \in X \mid \varphi(x)\}$.
\end{definition}

These definitions are used all throughout mathematics and category theory. 

\begin{example}
  Let $\cC$ be a category and let $f\colon X \to Y$, $g\colon Y \to Z$ be two morphisms. Then the set $\{h \colon X \to Z \mid h \text{ is the composition of } f \text{ and } g \}$ has a unique element, namely the composite $g \circ f$.
\end{example}

In the context of homotopy theory, these notions of uniqueness need to be adapted. Indeed, we are replacing sets with Kan complexes and bijections with homotopy equivalences. We hence need a homotopical analogue of uniqueness.

\begin{definition}[Contractibility]
  A Kan complex $K$ is called \emph{contractible} if the unique map $K \to \Delta[0]$ (\cref{lemma:delta 0}) is a homotopy equivalence.
\end{definition}

\begin{remark}
  Contractible Kan complexes are precisely the homotopical analogue of singleton sets.
\end{remark}

Previously we considered subsets associated to a given formula $\varphi$. This can be directly generalized to the homotopical setting.

\begin{definition}[Sub-Kan complex]
  Let $K$ be a Kan complex. A \emph{sub-Kan complex} $L \subseteq K$ is a sub-simplicial set that is also a Kan complex.
\end{definition}

\begin{definition}[Homotopically unique] \label{def:homotopicallyunique}
  Let $K$ be a Kan complex. A \emph{homotopical property} on $K$ is a sub-Kan complex $L \subseteq K$. It is \emph{homotopically unique} if $L$ is contractible.
\end{definition}

Homotopical uniqueness and contractibility are sufficient when we have a property on a single Kan complex and we want to analyze whether the property is homotopically unique. However, as we discussed in \cref{sec:kanfibrations}, we are often in a dependent situation, where we have a family of objects or properties that vary over some base. In the context of sets, for example, we have the following results.

\begin{lemma} \label{lemma:dependent uniqueness}
 Let $f\colon Y \to Z$ be a map of sets. Then $f$ is a bijection if and only if for every $z \in Z$ the pre-image $f^{-1}(z)$ is a singleton. 
\end{lemma}

\begin{intuition}
  We can think of the set $Y = \coprod_{z \in Z} f^{-1}(z)$ as a disjoint union of the pre-images. From this perspective, we can think of bijections as a ``global form of uniqueness over $Z$'' and \cref{lemma:dependent uniqueness} as the statement that this global form of uniqueness is equivalent to a ``local form of uniqueness'' for every $z$ in $Z$. 
\end{intuition}

We want an analogous result for Kan complexes, where bijections are replaced by homotopy equivalences and singletons are replaced by contractible Kan complexes. However, the following example shows that this will fail without some assumption.

\begin{example}
  Let $0\colon\Delta[0] \to N(I(1))$ be the homotopy equivalence defined in \cref{ex:preimage fail}. Then the pre-image of $1$ is the empty set, which is not contractible. Hence, being a global equivalence does not imply local homotopical uniqueness. 
\end{example}

Building on \cref{sec:kanfibrations}, we want to characterize dependent homotopical uniqueness as a suitable condition on Kan fibrations, as they model dependent structures.

\begin{definition}[Trivial fibration]
   A map of simplicial sets $f\colon Y \to Z$ is a \emph{trivial fibration} if it is a Kan fibration and a homotopy equivalence.
\end{definition}

\begin{remark} \label{rem:trivial fibration over point}
  This definition is indeed a simultaneous generalization of contractibility and dependence. In particular, a trivial Kan fibration $K \to \Delta[0]$ is precisely the data of a contractible Kan complex $K$ (\cref{rem:kan fibration kan complex}).
\end{remark}

We now have the following very powerful result about trivial fibrations, demonstrating its suitability for dependent homotopical uniqueness. For the proof, we refer the reader to \cite[Theorem 7.10]{goerssjardine2009simplicialhomotopytheory}.

\begin{proposition} \label{prop:trivial fib}
  Let $f\colon Y \to X$ be a map of simplicial sets. Then the following are equivalent.
  \begin{itemize}[leftmargin=*]
    \item $f$ is a trivial fibration.
    \item $f$ is a Kan fibration and for every $x$ in $X$, the pre-image $f^{-1}(x)$ is a contractible Kan complex.
    \item For every $n \geq 0$ and commutative diagram
    \[  
    \begin{tikzcd}
     \partial \Delta[n] \arrow[r] \arrow[d] & Y \arrow[d, "f", twoheadrightarrow, "\simeq"'] \\
     \Delta[n] \arrow[r] \arrow[ur, dashed] & X
    \end{tikzcd}
    \]
   there exists a lift.
  \end{itemize}
\end{proposition}

Before we apply dependent homotopical uniqueness in the next section, let us record how the third condition of \cref{prop:trivial fib} and a proof analogous to \cref{lemma:kan fibration pullback} give us the following result.

\begin{lemma} \label{lemma:trivial fibration pullback}
  Let the following diagram be a pullback square of simplicial sets
  \[
  \begin{tikzcd}
    W \arrow[r] \arrow[d, "q"', twoheadrightarrow, "\simeq"] & Y \arrow[d, "p", twoheadrightarrow, "\simeq"'] \\
    Z \arrow[r] & X
  \end{tikzcd}.
  \]
  If $p$ is a trivial Kan fibration, then $q$ is also a trivial Kan fibration.
\end{lemma}

\begin{remark} \label{rem:hott}
  Categories are usually defined in a set-theoretic context, and the notion of uniqueness in that foundation matches the intended categorical applications. As we saw in this section, homotopical uniqueness is more intricate in the set-theoretic context, and requires the contractibility of Kan complexes. This raises the question whether it is possible to adjust the underlying foundation so that uniqueness in that foundation matches homotopical uniqueness. This has been successfully achieved via the new foundation of \emph{homotopy type theory} \cite{hottbook2013}. See \cref{sec:next} for suggestions on how to actually develop higher category theory in such a setting. 
\end{remark}

\section{Segal Spaces} \label{sec:segal spaces} 
Let us recall the diagram of a simplicial space, as depicted in \cref{rem:simplicial space depiction}. In \cref{sec:reedy}, we introduced Reedy fibrancy to ensure that ``vertically'' our simplicial space behaves like a space. In this section we want to ensure that ``horizontally'' our simplicial space behaves like a category. We recall from \cref{thm:nerve image} that the key property was the Segal condition. Hence, we simply introduce the homotopical analogue to the Segal condition. We then observe how the resulting objects, Segal spaces, already exhibit many categorical features. Here, we fundamentally rely on the material in \cite[Section 5]{rezk2001css}, even though some exposition differs.

Recall from \cref{not:FDelta} the simplicial space $G(n) = \inccat(\Sp[n])$ that represents the spine of $F(n) = \inccat(\Delta[n])$. We can combine \cref{lemma:spine colimit,prop:colimit limit map} and repeat the proof from \cref{lemma:spine pullback} to get the following result.

\begin{lemma}
  Let $X$ be a simplicial space and $n \geq 2$. Then there is a natural isomorphism of simplicial sets
  \[\Map(G(n), X) \cong X_1 \times_{X_0}^{1^*,0^*} ... \times_{X_0}^{1^*,0^*} X_1,\]
  where there are $n$ factors of $X_1$ and the pullbacks are taken over the maps $0^*, 1^*\colon X_1 \to X_0$.
\end{lemma}

Previously, the Segal condition on a simplicial set gave us strict compositions. In the context of higher categories, we should generalize this condition.

\begin{definition}[Segal space] \label{def:segal space}
  A \emph{Segal space} is a Reedy fibrant simplicial space $T$ such that for every $n \geq 2$ the Kan fibration
  \[i^*_n\colon T_n \cong \Map(F(n),T) \twoheadrightarrow \Map(G(n),T) \cong T_1 \times_{T_0}^{1^*,0^*} ... \times_{T_0}^{1^*,0^*} T_1\]
  is a trivial Kan fibration.
\end{definition}

\begin{remark}
  Here the Reedy fibrancy of $T$ ensures that maps from $G(n)$ to $T$ correspond homotopy invariantly to compatible tuples in $T_1$, as explained in \cref{lemma:homotopy pullback}.
\end{remark}

\begin{remark} \label{rem:segal space lifting}
  Let us note here that we can characterize Segal spaces also with lifting conditions, analogous to \cref{def:reedy}. Namely, a Reedy fibrant simplicial space $T$ is a Segal space, if for all $n \geq 2, k \geq 0$, the following diagram admits a lift
  \[
  \begin{tikzcd}[row sep=1.2cm, column sep=1.2cm]
   G(n) \times \Delta[k] \coprod_{G(n) \times \partial \Delta[k]} F(n) \times \partial \Delta[k] \arrow[r] \arrow[d] & T \\
   F(n) \times \Delta[k] \arrow[ur, dashed]
  \end{tikzcd}.
  \]
  The proof follows from understanding how lifting conditions interact with products and hom-sets. See \cite[Section 21]{rezk2022qcats}, and particularly \cite[Proposition 21.5]{rezk2022qcats}, for an elegant introduction to this concept and a proof.
\end{remark}

We now commence with the category theory of Segal spaces.

 \begin{definition}[Object of a Segal space] \label{def:objects of segal space}
  Let $T$ be a Segal space. We define the \emph{objects of $T$} as the set of vertices in $T_0$. Thus, $\Obj_T = T_{00}$.
 \end{definition}

 \begin{remark}
  This definition is consistent with our characterization of nerves of categories, where $N\cC_0$ is precisely the set of objects (\cref{rem:nerve explicit}).
 \end{remark}

 In a category $\cC$, given two objects $x, y$, we have a set of morphisms $\Hom_\cC(x,y)$. In the higher categorical realm, we expect more structure.

 \begin{definition}[Morphism of a Segal space] \label{def:morphisms of segal space}
  Let $T$ be a Segal space. A \emph{morphism} in $T$ is a point in $T_1$, so an element in $T_{10}$.
 \end{definition}

For each morphism $f$, $0^*f$ is the \emph{source} and $1^*f$ is the \emph{target}. For two objects $x, y$ in $T$, we commonly use the notation $f\colon x \to y$, for a morphism $f$ from $x$ to $y$. 

 \begin{definition}[Mapping space of a Segal space] \label{def:mapping space Segal}
  For two objects $x,y$ in $T$, we define the \emph{mapping space} by the following pullback: 
  \[
  \begin{tikzcd}
    \map_{T} (x,y) \arrow[r] \arrow[d, twoheadrightarrow] & T_1 \arrow[d, "{(0^*,1^*)}", twoheadrightarrow] \\
    \Delta[0] \arrow[r, "{(x,y)}"] & T_0 \times T_0
  \end{tikzcd}
  \]
  or in other words the fiber of $(0^*,1^*)$ over the point $(x,y)$.
 \end{definition}

Notice we can actually call it a space, as we have a Kan complex.

\begin{lemma} \label{lemma:mapping space kan complex}
  Let $T$ be a Segal space and let $x,y$ be two objects. Then the mapping space $\map_T(x,y)$ is a Kan complex.
\end{lemma}

\begin{proof}
 By \cref{ex:reedy fibrant one}, the map $(0^*,1^*)$ is a Kan fibration. So, by \cref{lemma:kan fibration pullback}, the pullback along $(x,y)$, namely $\map_T(x,y) \to \Delta[0]$, is also a Kan fibration, which, by \cref{rem:kan fibration kan complex}, means that $\map_T(x,y)$ is a Kan complex.
\end{proof}

\begin{remark}
 Recall from \cref{rem:nerve hom set} that for a category $\cC$, we recover the hom set via an analogous pullback. This provides further evidence that this is a good way to define the mapping space in a Segal space.
\end{remark}

Up until now we have not used the Segal condition in any way. However, we want to use it to define composition.
 
\begin{intuition}[Composition via the Segal condition for n=2] \label{int:segal space n2}
  The first condition states that the map $T_2 \overset{\simeq}{\twoheadrightarrow} T_1 \times_{T_0}^{1^*,0^*} T_1$ is a trivial Kan fibration. Note that $T_2$ is the space of $2$-cells. Concretely, we depict a $2$-cell $\sigma$ as follows:
   \[
    \begin{tikzcd}[row sep=0.5cm, column sep=0.5cm]
     & y \arrow[dr, "g"]& \\
     x \arrow[rr, ""{name=U, above}, "h"'] \arrow[ur, "f"] & & z
     \arrow[start anchor={[xshift=-5ex, yshift=4ex]}, to=U, phantom, "\sigma"]
    \end{tikzcd}
   \]
   Similarly, we think of $T_1 \times_{T_0}^{1^*,0^*} T_1$ as the space of \emph{two composable arrows}, which we depict as:
    \[
    \begin{tikzcd}[row sep=0.5cm, column sep=0.5cm]
     & y \arrow[dr, "g"]& \\
     x \arrow[ur, "f"] & & z
    \end{tikzcd}
   \]
   The Segal condition states that every such diagram can be filled out to a complete $2$-cell: 
    
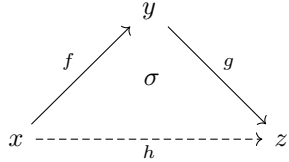
\[
    \begin{tikzcd}[row sep=0.5cm, column sep=0.5cm]
     & y \arrow[dr, "g"]& \\
     x \arrow[rr, ""{name=U, above}, "h"', dashed] \arrow[ur, "f"] & & z
     \arrow[start anchor={[xshift=-5ex, yshift=4ex]}, to=U, phantom, "\sigma"]
    \end{tikzcd}
   \]
   From this point of view, we think of $h$ as the composition of $f$ and $g$, and we think of $\sigma$ as a witness of that composition.
   Thus we often depict $h$ as $g \circ f$.
\end{intuition}

Based on this intuition we now have the following definition.

\begin{definition}[Composition in a Segal space]
  Let $T$ be a Segal space, let $x,y,z$ be objects, and let $f\colon x \to y$ and $g\colon y \to z$ be two morphisms. A \emph{composition} is a choice of lift
  \[ 
   \begin{tikzcd}[column sep=1.5cm]
    & G(2) \arrow[r, "{(f,g)}"] \arrow[d] & T \\ 
    F(1) \arrow[r, "02"'] \arrow[urr, dashed, "g \circ f" near start] & F(2) \arrow[ur, dashed]
   \end{tikzcd}
  \]
\end{definition}

Now if we have three morphisms $f\colon x \to y, g\colon y \to z, h\colon z \to w$, we might hope that the two compositions $h \circ (g \circ f)$ and $(h \circ g) \circ f$ coincide. Let us observe an intuitive argument why this should be the case.

\begin{intuition}[Associativity via the Segal condition for n=3] \label{int:segal space n3}
   The second condition states that the map $T_3 \overset{\simeq}{\twoheadrightarrow} T_1 \times^{1^*,0^*}_{T_0} T_1 \times^{1^*,0^*}_{T_0} T_1$ is a trivial Kan fibration. Here $T_3$ is the space of $3$-cells, which we can depict as a tetrahedron.
   
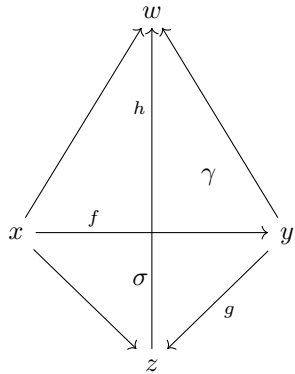
\[
    \begin{tikzcd}[row sep=0.5cm, column sep=0.5cm]
     & w & \\
     & & \\
     x \arrow[uur] \arrow[dr] \arrow[rr, ""{name=U, above}, "f" near start] & & y \arrow[dl, "g"] \arrow[uul] \\
     & z \arrow[uuu, "h" near end]
     \arrow[start anchor={[xshift=-1ex, yshift=1ex]}, to=U, phantom, "\sigma"]
     \arrow[start anchor={[xshift=5ex, yshift=10ex]}, to=U, phantom, "\gamma"]
    \end{tikzcd}
   \]
   where all $2$-cells and the middle $3$-cell are filled out. 
   On the other hand $T_1 \times_{T_0}^{1^*,0^*} T_1 \times_{T_0}^{1^*,0^*} T_1$ can be depicted as three composable arrows.
   
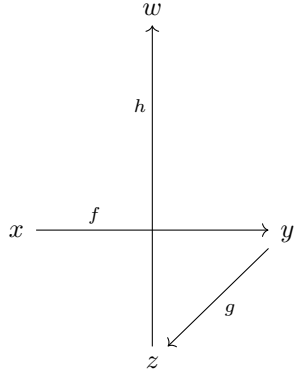
\[
    \begin{tikzcd}[row sep=0.5cm, column sep=0.5cm]
     & w & \\
     & & \\
     x \arrow[rr, ""{name=U, above}, "f" near start] & & y \arrow[dl, "g"] \\
     & z \arrow[uuu, "h" near end]
    \end{tikzcd}
   \]
   The Segal condition states that every such diagram can be filled out to a complete $3$-cell: 
   \[
    \begin{tikzcd}[row sep=1cm, column sep=1cm]
     & w & \\
     & & \\
     x \arrow[uur, dashed, "h \circ (g \circ f)" near start, "(h \circ g) \circ f" near end] 
     \arrow[dr, dashed, "g \circ f"'] \arrow[rr, ""{name=U, above}, "f" near start] & & 
     y \arrow[dl, "g"] \arrow[uul, dashed, "h \circ g"'] \\
     & z \arrow[uuu, "h" near end]
    \end{tikzcd}
   \]
   From the point of view of \cref{int:segal space n2}, the $1$-cell from $x$ to $w$ is both a choice of composition $(h \circ g) \circ f$ and $h \circ (g \circ f)$. 
\end{intuition}
   
\cref{int:segal space n2} and \cref{int:segal space n3} are first indications of how a Segal space has categorical features. We now want to make these intuitive ideas into a precise argument. This requires generalized mapping spaces, the space of compositions, and homotopic morphisms.

\begin{definition}[Homotopic morphisms in a Segal space] \label{def:homotopic morphisms}
 Let $T$ be a Segal space, let $x,y$ be two objects and let $f,g$ in $\map_T(x,y)$ be two morphisms. We say $f$ and $g$ are \emph{homotopic} if there is a map $\Delta[1] \to \map_T(x,y)$ sending $0$ to $f$ and $1$ to $g$. 
\end{definition}

\begin{notation}
  If $f,g\colon x \to y$ are two morphisms in a Segal space $T$ that are homotopic, we write $f \simeq g$.
\end{notation}

Our notation is justified by the following basic observation, which we proved in \cref{lemma:homotopy equivalence relation}.

\begin{lemma} \label{lemma:homotopy equivalence relation segal}
  Let $T$ be a Segal space and let $x,y$ be two objects. Then the relation $\simeq$ on $\map_T(x,y)$ is an equivalence relation.
\end{lemma}

\begin{definition}[Generalized mapping space of a Segal space] \label{def:space of compositions}
 Let $T$ be a Segal space and let $x_0,x_1, ..., x_n$ be objects in $T$. The \emph{generalized mapping space} $\map_{T}(x_0,...,x_n)$ is defined as the following pullback:
 \[ 
 \begin{tikzcd}
    \map_{T}(x_0,...,x_n) \arrow[r] \arrow[d, twoheadrightarrow] & T_n \arrow[d, "{(0^*,1^*,...,n^*)}", twoheadrightarrow] \\
    \Delta[0] \arrow[r, "{(x_0,x_1, ..., x_n)}"] & (T_0)^{n+1}.
 \end{tikzcd}
 \]
\end{definition}

Using \cref{ex:reedy fibrant n}, we have the following result, whose statement and proof is analogous to \cref{lemma:mapping space kan complex}.

\begin{lemma} \label{lemma:generalized mapping space kan complex}
  Let $T$ be a Segal space and let $x_0,x_1, ..., x_n$ be objects in $T$. Then $\map_T(x_0,...,x_n)$ is a Kan complex.
\end{lemma}

Again, up to here we have not used the Segal condition. However, using the general Segal condition we have the following result about generalized mapping spaces and mapping spaces.

\begin{proposition} \label{prop:generalized mapping space mapping space}
  Let $T$ be a Segal space and let $x_0,x_1, ..., x_n$ be objects in $T$, where $n \geq 2$. Then the trivial Kan fibration given in \cref{def:segal space} induces a trivial Kan fibration
  \[\map_T(x_0,...,x_n) \overset{ \ \ \simeq \ \ }{\twoheadrightarrow} \map_T(x_0,x_1) \times  ... \times \map_T(x_{n-1},x_n).\]
\end{proposition}

\begin{proof}
  We have the following commutative diagram
  \[
  \begin{tikzcd}[column sep=1.5cm]
    \Delta[0] \arrow[r, "{(x_0,x_1, ..., x_n)}"]  \arrow[d, equal] & (T_0)^{n+1} \arrow[d, equal] & T_n \arrow[d, "\simeq", twoheadrightarrow]  \arrow[l, twoheadrightarrow] \\
    \Delta[0] \arrow[r, "{(x_0,x_1, ..., x_n)}"] & (T_0)^{n+1} & T_1 \times_{T_0} ... \times_{T_0} T_1 \arrow[l, twoheadrightarrow]
  \end{tikzcd}.
  \]
  Here all vertical maps are homotopy equivalences. So, by \cref{rem:homotopy pullback}, the induced map on pullbacks will be a homotopy equivalence as well. By \cref{def:space of compositions}, the pullback of the top row is precisely $\map_T(x_0,...,x_n)$. Thus we are left with computing the pullback of the bottom row.

  Let us denote the pullback of the bottom diagram by $P_{n+1}$. We proceed by induction on $n$. If $n =1$, then by \cref{def:mapping space Segal}, $P_{1+1} = P_2 = \map_T(x_0,x_1)$. Now assuming the result holds for $n$, we have the following commutative diagram
  \[ 
  \begin{tikzcd}
    \Delta[0] \arrow[r, "{(x_0, ..., x_{n-1})}"] \arrow[d, equal] & T_0^n \arrow[d, "\pi_{n}"] & T_1 \times^{1^*,0^*}_{T_0} ... \times^{1^*,0^*}_{T_0} T_1  \arrow[d, "1^*\pi_{n-1}"]  \arrow[l, twoheadrightarrow, "{(0^*, ..., (n-1)^*)}"'] &[-.7cm] P_{n} \arrow[d] \\
    \Delta[0] \arrow[r, "{x_{n-1}}"] & T_0 & T_0 \arrow[l, equal] \arrow[r, rightsquigarrow] & \Delta[0]\\
    \Delta[0] \arrow[r, "{(x_{n-1},x_n)}"] \arrow[u, equal] & T_0 \times T_0 \arrow[u, "\pi_1"'] \arrow[d, rightsquigarrow] & T_1 \arrow[u, "0^*"'] \arrow[l, twoheadrightarrow, "{(0^*,1^*)}"'] & \map_T(x_{n-1},x_n) \arrow[u] \\ 
    \Delta[0] \arrow[r, "{(x_0,... ,x_n)}"] & T_0^{n+1} & T_1 \times_{T_0} ... \times_{T_0} T_1 \arrow[l, twoheadrightarrow, "{(0^*,...,n^*)}"'] \\
  \end{tikzcd}.
  \]
  Taking the pullback of each column gives us the diagram below it, which is the desired diagram, whereas taking the pullback of each row gives us the diagram to the right. By pullback commutativity (\cref{lemma:pullback commutativity}), the pullbacks of these two resulting diagrams coincide, meaning we have equivalences 
  \[P_{n+1} \overset{\simeq}{\twoheadrightarrow} P_n \times \map_T(x_{n-1},x_n) \overset{\simeq}{\twoheadrightarrow} ... \overset{\simeq}{\twoheadrightarrow} \map_T(x_0,x_1) \times ... \times \map_T(x_{n-1},x_n),\]
  finishing the result.
\end{proof}

We can now use these results to define and study the space of compositions.

\begin{definition}[Space of compositions in a Segal space]
 Let $T$ be a Segal space, let $x,y,z$ be three objects, and let $f\colon x \to y, g\colon y \to z$ be two morphisms. The \emph{space of compositions} of $f$ and $g$ is defined by the following pullback:
 \[ 
 \begin{tikzcd}
   \Comp_T(f,g) \arrow[r] \arrow[d, twoheadrightarrow] & \map_T(x,y,z) \arrow[d, twoheadrightarrow, "{(01^*,12^*)}"] \\
   \Delta[0] \arrow[r, "{(f,g)}"] & \map_T(x,y) \times \map_T(y,z).
 \end{tikzcd}
 \]  
\end{definition}

\begin{proposition} \label{prop:composition contractible}
  Let $T$ be a Segal space, let $x,y,z$ be three objects, and let $f\colon x \to y, g\colon y \to z$ be two morphisms. Then the space of compositions $\Comp_T(f,g)$ is a contractible Kan complex.
\end{proposition}

\begin{proof}
 By \cref{prop:generalized mapping space mapping space}, the map $\map_T(x,y,z) \to \map_T(x,y) \times \map_T(y,z)$ is a trivial Kan fibration. This means $\Comp_T(f,g) \to \Delta[0]$ is a trivial Kan fibration (\cref{lemma:trivial fibration pullback}), meaning it is a contractible Kan complex (\cref{rem:trivial fibration over point}).
\end{proof}
 
\begin{intuition}
 A point in $\Comp_T(f,g)$ is precisely a choice of composition of $f$ and $g$, i.e.~a diagram
 \[ 
  \begin{tikzcd}[row sep=0.5cm, column sep=0.5cm]
    & Y \arrow[dr, "g"]& \\
    X \arrow[rr, ""{name=U, above}, "h"'] \arrow[ur, "f"] & & Z
    \arrow[start anchor={[xshift=-5ex, yshift=4ex]}, to=U, phantom, "\sigma"]
  \end{tikzcd}.
 \]  
 From this perspective, \cref{prop:composition contractible} is telling us that the composition of $f$ and $g$ is \emph{homotopically} unique, in the sense of \cref{sec:contractibility}.
\end{intuition}

We can use this contractibility result to show any two choices of compositions in the target mapping space are homotopic.

\begin{lemma} \label{lemma:homotopic compositions}
Let $T$ be a Segal space, let $x,y,z$ be three objects, and let $f\colon x \to y, g\colon y \to z$ be two morphisms. 
Then any two choices of compositions of $f$ and $g$ are homotopic in $\map_T(x,z)$.
\end{lemma}

\begin{proof}
Let $h, h' \colon x \to z$ be two choices of compositions of $f$ and $g$. Let $\sigma, \sigma'\colon \Delta[0] \to \Comp_T(f,g)$ be witnesses for the compositions $h,h'$, respectively.
We have the following diagram 
\[ 
\begin{tikzcd}
  \partial \Delta[1] \arrow[rrr, "h + h'", bend left = 10] \arrow[r, " \sigma + \sigma'"'] \arrow[d] & \Comp_T(f,g) \arrow[d, twoheadrightarrow, "\simeq"] \arrow[r] & \map_T(x,y,z) \arrow[r, "02^*"'] \arrow[d, twoheadrightarrow, "\simeq"]  & \map_T(x,z) \\
 \Delta[1] \arrow[r] \arrow[ur, "H" description, dashed] & \Delta[0] \arrow[r, "{(f,g)}"] & \map_T(x,y) \times \map_T(y,z)
\end{tikzcd},
\]
which lifts by \cref{prop:trivial fib}. The composition map $\Delta[1] \to \map_T(x,z)$ gives us the desired homotopy between $h$ and $h'$.
\end{proof}

In fact we can use this proof method to show a related, and very important, result.

\begin{lemma} \label{lemma:homotopic compositions well defined}
 Let $T$ be a Segal space, let $x,y,z$ be three objects, and let $f,f' \colon x \to y, g, g' \colon y \to z$ be four morphisms. Moreover, let $H_f$ be a homotopy from $f$ to $f'$, and let $H_g$ be a homotopy from $g$ to $g'$. Then any chosen composition of \(f\) and \(g\) is homotopic to any chosen composition of \(f'\) and \(g'\).
\end{lemma}

\begin{proof}
  Let $h, h' \colon x \to z$ be choices of compositions of $f,g$ and $f',g'$, respectively. Moreover, let $\sigma, \sigma'\colon \Delta[0] \to \map_T(x,y,z)$ be witnesses for $h,h'$, respectively. We have the following diagram
 \[ 
 \begin{tikzcd}[column sep = 2cm]
  \partial \Delta[1] \arrow[rr, "h + h'", bend left = 10] \arrow[r, " \sigma + \sigma'"'] \arrow[d] & \map_T(x,y,z) \arrow[r, "02^*"'] \arrow[d, twoheadrightarrow, "\simeq"]  & \map_T(x,z) \\
 \Delta[1] \arrow[ur, "H" description, dashed] \arrow[r, "{(H_f,H_g)}"'] & \map_T(x,y) \times \map_T(y,z)
\end{tikzcd},
\]
which lifts by \cref{prop:trivial fib}. The composition map $\Delta[1] \to \map_T(x,z)$ gives us the desired homotopy between the two compositions.
\end{proof}

We can use these observations to construct ``pre-composition'' and ``post-composition'' maps.

\begin{lemma} \label{lemma:pre post composition}
  Let $T$ be a Segal space, let $x,y,z$ be three objects, and let $f\colon x \to y$ be a morphism. Then there are maps of Kan complexes
  \begin{itemize}
    \item $f^*\colon\map_T(y,z) \to \map_T(x,z)$ sending $g$ to $g \circ f$
    \item $f_*\colon\map_T(z,x) \to \map_T(z,y)$ sending $g$ to $f \circ g$
  \end{itemize}
\end{lemma}

\begin{proof}
 We consider the first case and the second one is analogous. We have the following diagram
\[
 \begin{tikzcd}
  \map_T(y,z) \arrow[r, "\{f\} \times \id"] & \map_T(x,y) \times \map_T(y,z) \arrow[r, dashed, bend left= 10] & \map_T(x,y,z) \arrow[l, "\simeq", twoheadrightarrow] \arrow[r] & \map_T(x,z)
 \end{tikzcd}.
 \]
By \cref{prop:generalized mapping space mapping space}, there is a homotopy inverse map $\map_T(x,y) \times \map_T(y,z) \to  \map_T(x,y,z)$. Composing all these morphisms gives us a map $\map_T(y,z) \to \map_T(x,z)$, which by definition pre-composes with $f$. 
\end{proof}

\begin{remark}
 We should note here the maps $f^*,f_*$ defined in \cref{lemma:pre post composition} depend on a concrete choice of homotopy inverse for the given trivial Kan fibration, which generally cannot be chosen in a canonical or natural manner.
\end{remark}

\begin{lemma} \label{lemma:homotopic morphisms pre post composition}
  Let $T$ be a Segal space, let $x,y,z$ be three objects, and let $f,g\colon x \to y$ be two morphisms. If $f \simeq g$, then $f^* \simeq g^*$ and $f_* \simeq g_*$.
\end{lemma}

\begin{proof}
 We consider the case of pre-composition and the other case is analogous. Let $H\colon \Delta[1] \to \map_T(x,y)$ be a homotopy from $f$ to $g$. We have the following diagram
\[
 \begin{tikzcd}
  \Delta[1] \times \map_T(y,z) \arrow[r, "H \times \id"] & \map_T(x,y) \times \map_T(y,z) \arrow[r, dashed, bend left= 10, start anchor={[yshift=-0.15cm, xshift=0.15cm]}] & \map_T(x,y,z) \arrow[l, "\simeq", twoheadrightarrow] \arrow[r] & \map_T(x,z)
 \end{tikzcd}.
 \]
 Again, choosing a homotopy inverse for the middle map and composing gives us a map of simplicial sets $\Delta[1] \times \map_T(y,z) \to \map_T(x,z)$, which by definition is a homotopy between $f^*$ and $g^*$.
\end{proof}

We can now move on to further aspects of the category theory of Segal spaces, namely identity morphisms and associativity.

\begin{definition}[Identity map in a Segal space] \label{def:identity map}
  Let $T$ be a Segal space and let $x$ be an object in $T$. The \emph{identity map} of $x$ is the morphism $\id_x = 00^*(x)$ in $T_1$.
\end{definition}

\begin{remark}
  The simplicial identities $0^*00^* = 1^*00^* = \id$ ensure that $\id_x$ indeed has source and target $x$.
\end{remark}

We now have the following result. 

 \begin{proposition} \label{prop:associativity and identity of segal space}
  Let $T$ be a Segal space, let $x,y,z,w$ be four objects. Let $f\colon x \to y$, $g\colon y \to z$ and $h\colon z \to w$. Then, $h \circ (g \circ f)$ and $(h \circ g ) \circ f$ are homotopic, and both $\id_y \circ f, f \circ \id_x$ are homotopic to $f$. This means composition is homotopically associative and unital.
 \end{proposition}

 \begin{proof}
  We have the following commutative diagram
  \[
   \begin{tikzcd}[row sep =0.5cm, column sep=0.3cm]
    & \map_T(x,y) \times \map_T(y,z) \times \map_T(z,w) & \\
    & \map_T(x,y,z,w) \arrow[dl, "013^*"', twoheadrightarrow] \arrow[dr, "023^*", twoheadrightarrow]  \arrow[u, twoheadrightarrow, "\simeq"] & \\
    \map_T(x,y,w) \arrow[dr, "02^*"', twoheadrightarrow] & & \map_T(x,z,w) \arrow[dl, "02^*", twoheadrightarrow] \\
    & \map_T(x,w) &
   \end{tikzcd}
  \]
  Let us take a point $(f,g,h)$ in $\map_T(x,y) \times \map_T(y,z) \times \map_T(z,w)$. Via the equivalence, we can lift it to an element $\sigma$ in $\map_T(x,y,z,w)$. Going through the left-hand map gives us $(h \circ g) \circ f$, but the right-hand map gives us $h \circ (g \circ f)$. This proves associativity.
   
  For the identity relation, let $f\colon x \to y$. This gives us a $2$-cell $001^*(f)$ in $\map_T(x,x,y)$, which satisfies $(01^*001^*(f),02^*001^*(f),12^*001^*(f)) = (\id_x,f,f)$, as we can evidently compute. This proves that composition of $f$ and $\id_x$ is $f$. The other case follows similarly via $011^*(f)$.
 \end{proof}

Identity and associativity help us better understand the properties of the pre-composition and the post-composition maps on mapping spaces.

\begin{lemma}\label{lemma:pre post composition id}
  Let $T$ be a Segal space and let $x,y$ be two objects. Then $\id_x^* \simeq \id_{\map_T(x,y)}\simeq (\id_y)_*$.
\end{lemma}

\begin{proof}
 We only prove the first equivalence and the second one is analogous. For a morphism $f$, $\id_x^*(f) = f \circ \id_x \simeq f = \id_{\map_T(x,y)}(f)$.
\end{proof}

\begin{lemma} \label{lemma:pre post composition compose}
  Let $T$ be a Segal space, let $x,y,z$ be three objects, and let $f\colon x \to y$, $g\colon y \to z$ be two morphisms. Then $f^* \circ g^* \simeq (g \circ f)^*$ and $g_* \circ f_* \simeq (g \circ f)_*$.
\end{lemma}

\begin{proof}
 We prove the case of pre-composition and the other case is analogous. For a given object $w$ and morphism $h\colon z \to w$, we have 
 \[f^* \circ g^*(h) \simeq f^*(h \circ g) \simeq (h \circ g) \circ f \simeq h \circ (g \circ f)  \simeq (g \circ f)^*(h),\]
 where we used the definition of $f^*,g^*$ and associativity (\cref{prop:associativity and identity of segal space}).
\end{proof}

Let us move on to functors of Segal spaces. 

\begin{definition}[Functor of Segal spaces] \label{def:functor of segal spaces}
  Let $T$ and $U$ be two Segal spaces. A \emph{functor} $F\colon T \to U$ is a map of simplicial spaces.
\end{definition}

On the surface, this seems to differ from the classical definition of functors between categories. However, we have the following result.

\begin{proposition} \label{prop:functoriality}
  Let $F\colon T \to U$ be a functor between Segal spaces. Then $F$ induces the following:
  \begin{itemize}[leftmargin=*]
    \item A function $\Obj_T \to \Obj_U$ between the objects of $T$ and $U$.
    \item For every two objects $x,y$ in $T$, a map of mapping spaces $\map_T(x,y) \to \map_U(Fx,Fy)$.
    \item For every three objects $x,y,z$ in $T$ and morphisms $f\colon x \to y$, $g\colon y \to z$, a homotopy between $F(g \circ f)$ and $Fg \circ Ff$.
    \item For every object $x$ in $T$, a homotopy between $F(\id_x)$ and $\id_{Fx}$.
  \end{itemize}
\end{proposition}

\begin{proof}
 First, $F$ induces a function $F_{00}\colon \Obj_T = T_{00} \to U_{00} = \Obj_U$. Next, we have the commutative diagram 
 \[ 
  \begin{tikzcd}
     {\Delta[0]} \arrow[r, "{(x,y)}"] \arrow[d, equal] & T_0 \times T_0 \arrow[d, "F_0 \times F_0"] &  T_1 \arrow[d, "F_1"] \arrow[l,twoheadrightarrow]\\
    {\Delta[0]} \arrow[r, "{(Fx,Fy)}"] & U_0 \times U_0 & U_1 \arrow[l,twoheadrightarrow] 
  \end{tikzcd}.
 \]
 As the mapping spaces are just the pullbacks of these diagrams (\cref{def:mapping space Segal}), the existence of a map $F_1\colon\map_T(x,y) \to \map_U(Fx,Fy)$ follows from \cref{lemma:limit map}. Next, let $\sigma$ be a point in $\Comp(f,g)$, meaning $01^*\sigma = f$, $12^*\sigma = g$. Then, by the simplicial identities $F_2\sigma$ is a point in $\Comp(F_1f,F_1g)$, as we have $01^*\circ F_2\sigma = F_1 \circ 01^*\sigma = F_1f$, $12^*\circ F_2\sigma = F_1 \circ 12^*\sigma = F_1g$. This gives us the desired homotopy between $F_1(02^*\sigma) = F_1(g \circ f)$ and $F_1g \circ F_1f$, as any two possible compositions of $F_1f$ and $F_1g$ are homotopic (\cref{lemma:homotopic compositions}). Finally, for a given object $x$ in $T$, $F_1(\id_x) = F_1(00^*(x)) = 00^*(F_0x) = \id_{F_0x}$, where we again used the simplicial identities.
\end{proof}

With this result at hand, we can also define fully faithful functors of Segal spaces.

\begin{definition}[Fully faithful functor of Segal spaces]
  Let $F\colon T \to U$ be a functor between Segal spaces. We say $F$ is \emph{fully faithful} if for every two objects $x,y$ in $T$, the induced map of mapping spaces $\map_T(x,y) \to \map_U(Fx,Fy)$ is a homotopy equivalence of Kan complexes.
\end{definition}

Having defined functors between Segal spaces, we can also consider homotopies and equivalences.  

\begin{definition}[Homotopy of functors of Segal spaces]
  Let $F,G\colon T \to U$ be two functors between Segal spaces. A \emph{homotopy} $H \colon F \to G$ is a map of simplicial spaces $H\colon T \times \Delta[1] \to U$ such that $H|_{T \times 0^*} = F$ and $H|_{T \times 1^*} = G$.
\end{definition}

\begin{definition}[Equivalence of Segal spaces]
  Let $F\colon T \to U$ be a functor between Segal spaces. We say $F$ is an \emph{equivalence} if there exists a functor $G\colon U \to T$, such that $G \circ F$ is homotopic to $\id_T$ and $F \circ G$ is homotopic to $\id_U$.
\end{definition}

\begin{remark} \label{rem:equivalence to homotopy equivalence}
 Recall from \cref{def:incspace}, $\Delta[1]_n = \Delta[1]$. Thus, a homotopy $H\colon T \times \Delta[1] \to U$ gives us a homotopy $H_n \colon T_n \times \Delta[1] \to U_n$ for every $n \geq 0$. This means every equivalence $F\colon T \to U$ induces a homotopy equivalence of spaces $F_n \colon T_n \to U_n$ for every $n \geq 0$.
\end{remark}

We end this section with the observation that every Segal space has an associated homotopy category.

\begin{construction} \label{constr:homotopy category of segal space}
  Let $T$ be a Segal space. We now construct a directed graph with a composition operation and identities. This construction is crucial when we look at the \emph{homotopy category} (\cref{thm:homotopy category of segal space}). Let $\Ho(T)$ consist of the following data:
  \begin{itemize}[leftmargin=*]
    \item \textbf{Objects:} The objects of $\Ho(T)$ are the objects of $T$.
    \item \textbf{Morphisms:} For two objects $x,y$ in $\Ho(T)$ we define the set of morphisms as
    \[\Hom_{\Ho(T)}(x,y) \coloneq \map_{T}(x,y) / \sim,\]
    where $\sim$ is the homotopy relation defined in \cref{def:homotopic morphisms,lemma:homotopy equivalence relation}.
    \item \textbf{Composition:} For three objects $x,y,z$ in $\Ho(T)$, the composition map 
    \[\map_{T}(x,y) \times \map_{T}(y,z) \to \map_{T}(x,z)\]
    gives us a composition map
    \[ \Hom_{\Ho(T)}(x,y) \times \Hom_{\Ho(T)}(y,z) \to \Hom_{\Ho(T)}(x,z),\]
    because homotopy respects composition, by \cref{lemma:homotopic compositions,lemma:homotopic compositions well defined}.
  \end{itemize}
\end{construction}

This construction gives us the following theorem.

\begin{theorem}[Homotopy category of a Segal space] \label{thm:homotopy category of segal space}
  Let $T$ be a Segal space. Then $\Ho(T)$ is a category, called the \emph{homotopy category} of $T$.
\end{theorem}

\begin{proof}
  We only need to verify that the composition operation is well-defined, associative and unital. However, this is the statement of \cref{lemma:homotopic compositions,lemma:homotopic compositions well defined,prop:associativity and identity of segal space}. 
\end{proof}

This theorem allows us to construct an ordinary category out of every Segal space, confirming the connection between Segal spaces and categories. This connection also extends to functors.

\begin{theorem}
  Let $F\colon T \to U$ be a functor between Segal spaces. Then $F$ induces a functor between homotopy categories $\Ho(F)\colon \Ho(T) \to \Ho(U)$.
\end{theorem}

\begin{proof}
  By \cref{prop:functoriality}, $F$ induces a function between the objects of $\Ho(T)$ and $\Ho(U)$, and a map between the morphism sets. Moreover, the homotopies ensure that composition and identities are preserved. Thus, $\Ho(F)$ is indeed a functor.
\end{proof}
 
\section{Homotopy Equivalences in Segal Spaces} \label{sec:homotopy equiv in Segal Spaces}
In the previous section we saw many categorical aspects of Segal spaces, such as composition, identities, and functors. In this section we will focus on the homotopical aspects of Segal spaces, in particular homotopy equivalences, which are part of the theory that distinguishes Segal spaces from ordinary categories. In particular we will present several equivalent characterizations of homotopy equivalences, and observe a variety of homotopical uniqueness results (in the sense of \cref{sec:contractibility}) regarding homotopy equivalences. 

\begin{remark}
  All the definitions and results presented in this section are well-known to experts. However, there does not appear to be a systematic presentation of these results in the literature, beyond the first couple of steps that can already be found in \cite{rezk2001css}. As a result this section is somewhat more demanding than the rest of this note.
\end{remark}
 
\subsection{Homotopy Equivalences: First Steps}
We commence with the definition and first properties of homotopy equivalences in Segal spaces.

\begin{definition}[Homotopy equivalence in a Segal space] \label{def:homotopy equivalence}
  Let $T$ be a Segal space and let $x,y$ be objects in $T$. A morphism $f$ in $\map_{T}(x,y)$ is a homotopy equivalence if there are maps $g,h$ in $\map_{T}(y,x)$ such that $f \circ g$ is homotopic to $\id_y$ and $h \circ f$ is homotopic to $\id_x$.
\end{definition}
 
\begin{remark}
 We can interpret the condition in \cref{def:homotopy equivalence} as saying that $\id_y$ is a choice of composition of $f$ and $g$, and, similarly, $\id_x$ is a choice of composition of $h$ and $f$. 
\end{remark}

\begin{intuition} \label{int:homotopy equiv in ss}
  Although it might appear that the definition only involves the existence of three maps, in reality the nature of a  Segal space demands that there are several other important pieces of information. In particular, the fact that the compositions $f \circ g$ and $h \circ f$ correspond to the identity involves $2$-cells that witness these compositions. The information can be captured in a diagram of the following form 
  \[
   \begin{tikzcd}[row sep =0.5in, column sep=0.5in]
    x  \arrow[r, "\id_x"] \arrow[dr, "f"] & x \\
    y \arrow[u, "g"] \arrow[r, "\id_y"'] & y \arrow[u, "h"'] 
   \end{tikzcd}
  \]
\end{intuition}

Note that, similar to the case for classical categories, inverses are homotopically unique, consistent with our philosophy (\cref{sec:contractibility}).

\begin{lemma} \label{lemma:inverses homotopic}
 Let $T$ be a Segal space, let $x,y$ be objects in $T$, and let $f$ in $\map_{T}(x,y)$ be a homotopy equivalence with inverses $g$ and $h$. Then $g$ and $h$ are homotopic.
\end{lemma}

\begin{proof}
  We have the following chain of homotopies
  \begin{align*}
    g & \simeq \id_x \circ g  & \text{\cref{prop:associativity and identity of segal space}} \\  
      & \simeq (h \circ f) \circ g & \text{Definition of } h  \\  
      & \simeq h \circ (f \circ g) & \text{\cref{prop:associativity and identity of segal space}} \\  
      & \simeq  h \circ \id_y & \text{Definition of } g  \\  
      & \simeq h & \text{\cref{prop:associativity and identity of segal space}}
  \end{align*}  \qedhere
\end{proof}

Given their importance we want various alternative ways to characterize homotopy equivalences. First we observe homotopy equivalences via isomorphism in the homotopy category.

\begin{proposition} \label{prop:homotopy equiv iff iso in homotopy category}
  Let $T$ be a Segal space and let $f\colon x \to y$ be a morphism in $T$. Then $f$ is a homotopy equivalence if and only if $[f]$ is an isomorphism in $\Ho(T)$.
\end{proposition}

\begin{proof}
 If $f$ is a homotopy equivalence, then there are maps $g,h$ in $\map_{T}(y,x)$ such that $f \circ g$ is homotopic to $\id_y$ and $h \circ f$ is homotopic to $\id_x$. This means $[f] \circ [g] = [\id_y]$ and $[h] \circ [f] = [\id_x]$, so $[f]$ is an isomorphism in $\Ho(T)$.
 
 On the other hand, if $[f]$ is an isomorphism in $\Ho(T)$, then there is a morphism $[g]$ in $\Ho(T)$ such that $[f] \circ [g] = [\id_y]$ and $[g] \circ [f] = [\id_x]$. Following \cref{thm:homotopy category of segal space}, $[g] \circ [f] = [g\circ f] = [\id_x]$, meaning $g \circ f$ is homotopic to $\id_x$. Similarly, $f\circ g$ is homotopic to $\id_y$. Thus, $f$ is a homotopy equivalence.
\end{proof}

\begin{intuition}
  Fundamentally, the homotopy category $\Ho(T)$ only remembers the objects and homotopy classes of morphisms in $T$. Thus \cref{prop:homotopy equiv iff iso in homotopy category} can be understood as saying that ``being a homotopy equivalence'' does not contain any non-trivial homotopical information. See \cref{intuition:hoeqchoice is hoequiv} for a more detailed discussion of this point.
\end{intuition}

Let us move on to another characterization of homotopy equivalences. In classical category theory, we can identify isomorphisms as those morphisms whose pre-composition or post-composition induces bijections on hom sets (\cref{cor:iso iff natural iso}). Let us generalize that.

\begin{proposition} \label{prop:pre post composition homotopy equivalence}
  Let $T$ be a Segal space, let $x,y$ be two objects in $T$, and let $f\colon x \to y$ be a morphism. Then the following are equivalent:
\begin{enumerate}[leftmargin=*]
  \item $f$ is a homotopy equivalence.
  \item For every object $z$, the maps $f_*\colon\map_T(z,x) \to \map_T(z,y)$ are homotopy equivalences of Kan complexes.
  \item For every object $z$, the maps $f^*\colon\map_T(y,z) \to \map_T(x,z)$ are homotopy equivalences of Kan complexes.
\end{enumerate}
\end{proposition}

\begin{proof}
  We prove $(1) \Leftrightarrow (2)$ and the case $(1) \Leftrightarrow (3)$ is analogous. Let us assume $f$ is a homotopy equivalence. Let us denote the right inverse of $f$ by $g$, and note there is a homotopy $f \circ g \simeq \id_y$. Then, we have 
  \begin{align*}
    f_* \circ g_* & \simeq (f \circ g)_* & \text{\cref{lemma:pre post composition compose}} \\  
      & \simeq (\id_y)_* & \text{\cref{lemma:homotopic morphisms pre post composition}}  \\  
      & \simeq \id & \text{\cref{lemma:pre post composition id}}.
   \end{align*}
   Using a left inverse $h$ of $f$ similarly gives us $h_* \circ f_* \simeq \id$. Thus, $f_*$ is a homotopy equivalence with right inverse $g_*$ and left inverse $h_*$.

   On the other hand, let us assume $f_*$ is a homotopy equivalence. Then, as every homotopy equivalence is a weak equivalence (\cref{lemma:homotopy implies weak homotopy}), for every object $z$, we get a bijection of hom sets $\Ho ([f]_*)\colon \Hom_{\Ho (T)}(z,x) \to \Hom_{\Ho (T)}(z,y)$. Thus, by \cref{cor:iso iff natural iso}, $[f]$ is an isomorphism in $\Ho(T)$, and by \cref{prop:homotopy equiv iff iso in homotopy category}, $f$ is a homotopy equivalence.
\end{proof}

We can use equivalences in Segal spaces to generalize various common categorical notions.
 
 \begin{definition}[Segal space groupoid] \label{def:segal space groupoid}
  We say a Segal space $T$ is a \emph{Segal space groupoid} if every morphism is a homotopy equivalence.
 \end{definition}

\begin{definition}[Essentially surjective functor of Segal spaces] \label{def:essentially surjective segal}
  Let $F\colon T \to U$ be a functor between Segal spaces. We say $F$ is \emph{essentially surjective} if for every object $y$ in $U$, there exists an object $x$ in $T$ such that there is a homotopy equivalence from $Fx$ to $y$ in $U$.
\end{definition}

\subsection{Homotopy Equivalences: Structure vs. Property}
We now move on to more intrinsic characterizations of homotopy equivalences, via various spaces of equivalences. We construct several spaces of equivalences, using well-chosen colimit and limit constructions.

\begin{definition} \label{def:z3}
 Let $Z(3)$ be the simplicial space defined by the colimit of the following diagram.
 \[
  \begin{tikzcd}[row sep=.5cm, column sep=0.5cm]
   F(1) & & F(1) & & F(1) \\
   & F(0) \arrow[ur, "1"'] \arrow[ul, "1"]& & F(0) \arrow[ur, "0"'] \arrow[ul, "0"]& 
  \end{tikzcd}
  \]
\end{definition}

In particular, by \cref{prop:colimit limit map}, we have 
\[\Map(Z(3), T) \cong T_1 {\times}^{1^*,1^*}_{T_0} T_1 {\times}^{0^*,0^*}_{T_0} T_1.\]
As a next step, we use this description to construct two maps into $\Map(Z(3),T)$.
 
For this next lemma, we denote the elements in 
\[(F(2) \coprod_{F(1)}^{12,01} F(2))_1 = \{00_0,01_0,11_0,12_0,22_0,02_0,12_1,02_1,22_1\},\]
where the subscript indicates whether the element comes from the left or right copy of $F(2)$.
 
\begin{lemma} \label{lemma:zthree in fthree}
  Let $T$ be a Segal space. There are natural maps
  \[ 02_0 + 12_0 + 02_1\colon Z(3) \to F(2) \coprod_{F(1)}^{12,01} F(2),\]
  \[ 012 + 123 \colon F(2) \coprod_{F(1)}^{12,01} F(2) \to F(3)\]
  which induce a chain of Kan fibrations
  \[ T_3 \overset{P_{32}}{\twoheadrightarrow}  T_2 \times_{T_1}^{12^*,01^*} T_2 \overset{P_{21}}{\twoheadrightarrow}  T_1 \times^{1^*,1^*}_{T_0} T_1 \times^{0^*,0^*}_{T_0} T_1.\]
\end{lemma}

 \begin{proof}
  By \cref{prop:colimit limit}, we need to construct an element in the set 
  {\scriptsize\[
  \Hom(F(1),F(2) \coprod_{F(1)}^{12,01} F(2)) \underset{\Hom(F(0),F(2) \coprod_{F(1)}^{12,01} F(2))}{\times} \Hom(F(1),F(2) \coprod_{F(1)}^{12,01} F(2)) \underset{\Hom(F(0),F(2) \coprod_{F(1)}^{12,01} F(2))}{\times} \Hom(F(1),F(2) \coprod_{F(1)}^{12,01} F(2)),
  \]}
  which by the Yoneda lemma (\cref{lemma:FDelta yoneda}) corresponds to an element in the set 
  \[ 
  (F(2) \coprod_{F(1)}^{12,01} F(2))_1 \times_{(F(2) \coprod_{F(1)}^{12,01} F(2))_0} (F(2) \coprod_{F(1)}^{12,01} F(2))_1 \times_{(F(2) \coprod_{F(1)}^{12,01} F(2))_0} (F(2) \coprod_{F(1)}^{12,01} F(2))_1.
  \]
  Such an element is given by the following evidently commutative diagram
  \[ 
  \begin{tikzcd}[row sep=0.5cm, column sep=0.5cm]
    & F(0) \arrow[dl, "1"'] \arrow[dr, "1"]& & F(0) \arrow[dl, "0"'] \arrow[dr, "0"]& \\
   F(1) \arrow[drr, "02_0"'] & & F(1) \arrow[d, "12_0"] & & F(1) \arrow[dll, "02_1"]\\
   & & F(2) \coprod_{F(1)}^{12,01} F(2) & & 
  \end{tikzcd}.
  \]
  With the same argument, again by \cref{prop:colimit limit}, we get the second map $F(2) \coprod_{F(1)}^{12,01} F(2) \to F(3)$ via the diagram 
  \[
  \begin{tikzcd}
    F(1) \arrow[r, "12"] \arrow[d, "01"'] & F(2) \arrow[d, "012"] \\ 
    F(2) \arrow[r, "123"'] & F(3)
  \end{tikzcd}.
  \]
  Now, by \cref{lemma:reedy}, these maps induce Kan fibrations
  \[ 
  \begin{tikzcd}[column sep=2cm]
   \Map(F(3),T) \arrow[r, twoheadrightarrow, "(012 + 123 )^*"] &  \Map(F(2) \coprod_{F(1)}^{12,01} F(2),T) \arrow[r, twoheadrightarrow, "(02_0 + 12_0 + 02_1)^*"] & \Map(Z(3),T), 
  \end{tikzcd},
  \]
  as $T$ is Reedy fibrant. Finally, by \cref{prop:colimit limit map}, this corresponds to the desired maps 
  \[ T_3 \overset{P_{32}}{\twoheadrightarrow}  T_2 \times_{T_1}^{12^*,01^*} T_2 \overset{P_{21}}{\twoheadrightarrow}  T_1 \times^{1^*,1^*}_{T_0} T_1 \times^{0^*,0^*}_{T_0} T_1.\]
\end{proof}

We want to make some observations about $P_{32}$. This requires defining a helper function.

\begin{lemma} \label{lemma:ctwo}
 Let $T$ be a Segal space. Then there is a map 
 \[ 01_0 + 12_0 + 12_1\colon F(1) \coprod^{1,0}_{F(0)} F(1) \coprod^{1,0}_{F(0)} F(1) \to F(2) \coprod^{12,01}_{F(1)} F(2) ,\]
 which induces a Kan fibration 
 \[C_{21}\colon T_2 \times^{12^*,01^*}_{T_1} T_2 \twoheadrightarrow T_1 \times^{1^*,0^*}_{T_0} T_1 \times^{1^*,0^*}_{T_0} T_1  \] 
\end{lemma}

\begin{proof}
  Following the same steps of the proof of \cref{lemma:zthree in fthree}, the existence of the map $01_0 + 12_0 + 12_1$ follows from the following commutative diagram and the property of colimits (\cref{prop:colimit limit}).
  \[ 
  \begin{tikzcd}[row sep=0.5cm, column sep=0.5cm]
    & F(0) \arrow[dl, "1"'] \arrow[dr, "0"]& & F(0) \arrow[dl, "1"'] \arrow[dr, "0"]& \\
   F(1) \arrow[drr, "01_0"'] & & F(1) \arrow[d, "12_0"] & & F(1) \arrow[dll, "12_1"]\\
   & & F(2) \coprod_{F(1)}^{12,01} F(2) & & 
  \end{tikzcd}.
  \]
  Then, by applying \cref{lemma:reedy} to $(01_0 + 12_0 + 12_1)^*$ and the limit computation of \cref{prop:colimit limit map}, we get the desired Kan fibration 
  \[C_{21}\colon T_2 \times^{12^*,01^*}_{T_1} T_2 \twoheadrightarrow T_1 \times^{1^*,0^*}_{T_0} T_1 \times^{1^*,0^*}_{T_0} T_1. \qedhere\]
\end{proof}

\begin{lemma} \label{lemma:ptwoone equivalence}
 Let $T$ be a Segal space. The map $C_{21}\colon T_2 \times^{12^*,01^*}_{T_1} T_2 \to T_1 \times^{1^*,0^*}_{T_0} T_1 \times^{1^*,0^*}_{T_0} T_1$ is isomorphic to the map $i_2^* \times i_2^* \colon T_2 \times^{12^*,01^*}_{T_1} T_2 \to (T_1 \times_{T_0} T_1) \times^{\pi_2,\pi_1}_{T_1}  (T_1 \times_{T_0} T_1)$ and is a trivial Kan fibration.
\end{lemma}

\begin{proof}
  First we establish the isomorphism. First of all we have the following commutative diagram
  \[
  \begin{tikzcd}
   F(1) \arrow[r, "12"] \arrow[d, "01"] & G(2) \arrow[d, "01 + 12"] \\
   G(2) \arrow[r, "12 + 23"] & G(3) 
  \end{tikzcd}, 
  \]
  which induces a map $(01 + 12) + (12 + 23)\colon G(2) \coprod^{12,01}_{F(1)} G(2) \to G(3)$. Via the same argument as in \cref{lemma:spine colimit}, this map is an isomorphism. Next, we observe that the following diagram commutes
  \[ 
  \begin{tikzcd}
   \displaystyle G(2) \coprod^{12,01}_{F(1)} G(2) \arrow[dr, "i_2 + i_2"'] \arrow[rr, "\cong"', "(01 + 12) + (12 + 23)"] & & G(3) \arrow[dl, "01_0 + 01_1 + 12_1"] \\ 
   & \displaystyle F(2) \coprod_{F(1)}^{12,01} F(2)  
  \end{tikzcd}.
  \] 
   Finally, applying $\Map(-,T)$ to this commutative diagram and applying \cref{prop:colimit limit map}, we get the desired commutative diagram
 \[
 \begin{tikzcd}
  & T_2 \times_{T_1}^{12^*,01^*} T_2 \arrow[dl, "C_{21}"', twoheadrightarrow] \arrow[dr, "i_2^* \times i_2^*", twoheadrightarrow] & \\ 
  T_1 \times^{1^*,0^*}_{T_0} T_1 \times^{1^*,0^*}_{T_0} T_1 \arrow[rr, "\cong"] & & (T_1 \times^{1^*,0^*}_{T_0} T_1) \times^{\pi_2, \pi_1}_{T_1} (T_1 \times^{1^*,0^*}_{T_0} T_1)
 \end{tikzcd},
 \]
 where all morphisms are Kan fibrations,  by \cref{lemma:reedy}, as $T$ is Reedy fibrant.
 
 We now move on to the second part. By $2$-out-of-$3$ (\cref{rem:two out of three}), to show that $C_{21}$ is a trivial Kan fibration, it suffices to show that $i_2^* \times i_2^*$ is a homotopy equivalence. This map is the pullback of the following diagram 
 \[ 
 \begin{tikzcd}
  T_2 \arrow[r, twoheadrightarrow, "12^*"] \arrow[d, "\simeq", twoheadrightarrow] & T_1 \arrow[d, equal] & T_2 \arrow[l, twoheadrightarrow, "01^*"'] \arrow[d, "\simeq", twoheadrightarrow]\\
  T_1 \times_{T_0} T_1 \arrow[r, twoheadrightarrow, "\pi_2"] & T_1 & T_1 \times_{T_0} T_1 \arrow[l, twoheadrightarrow, "\pi_1"'] 
 \end{tikzcd},
 \]
 and so the result follows from \cref{lemma:homotopy pullback}, as all horizontal maps are Kan fibrations, and the vertical maps are homotopy equivalences. 
\end{proof}

\begin{lemma} \label{lemma:pthree equivalence}
 Let $T$ be a Segal space. The map $P_{32}\colon T_3 \to T_2 \times^{12^*,01^*}_{T_1} T_2$ is a trivial Kan fibration.
\end{lemma}

\begin{proof}
  The argument proceeds analogously to the proof of \cref{lemma:ptwoone equivalence}. By \cref{lemma:zthree in fthree}, $P_{32}$ is a Kan fibration. So, we only need to show it is a homotopy equivalence. The following diagram commutes (where we used \cref{lemma:spine colimit} for the case $n=3$)
 \[ 
 \begin{tikzcd}
  & G(3) \arrow[dl, "01_0 + 12_0 + 12_1"'] \arrow[dr, "01 + 12 + 23"] & \\
  F(2) \coprod_{F(1)}^{12,01} F(2) \arrow[rr, "012 + 123"'] & & F(3) 
 \end{tikzcd},
 \]
 which induces a commutative diagram
 \[
 \begin{tikzcd}
  T_3 \arrow[dr, twoheadrightarrow] \arrow[rr, "P_{32}", twoheadrightarrow] & & T_2 \times^{12^*,01^*}_{T_1} T_2  \arrow[dl, "i_2^* \times i_2^*", twoheadrightarrow] \\
  & T_1 \times^{1^*,0^*}_{T_0} T_1 \times^{1^*,0^*}_{T_0} T_1 & 
 \end{tikzcd}.
 \]
 By the Segal condition and \cref{lemma:ptwoone equivalence}, the diagonal maps are trivial Kan fibrations. Thus, by $2$-out-of-$3$  (\cref{rem:two out of three}), $P_{32}$ is a homotopy equivalence of Kan complexes.
\end{proof}

See \cref{intuition:homotopy equivalence tetrahedron} for some intuition about this equivalence.

\begin{lemma} \label{lemma:tetra to square}
  Let $T$ be a Segal space. There is a natural map
 \[ (11^* , 01^* , 00^*) \colon T_1 \to T_1 \times^{1^*,1^*}_{T_0} T_1 \times^{0^*,0^*}_{T_0} T_1.\]
\end{lemma}

\begin{proof}
  By \cref{prop:limit limit}, we need to construct an element in the set 
  \[\Hom(T_1,T_1) \times_{\Hom(T_1,T_0)} \Hom(T_1,T_1) \times_{\Hom(T_1,T_0)} \Hom(T_1,T_1).\] 
  Such an element is given by the following diagram
  \[ 
  \begin{tikzcd}[row sep=0.5cm, column sep=0.5cm]
      & & T_1 \arrow[dll, "11^*"'] \arrow[drr, "00^*"] \arrow[d, "01^*"] & \\
    T_1 \arrow[dr, "1^*"] & & T_1 \arrow[dl, "1^*"] \arrow[dr, "0^*"] & & T_1 \arrow[dl, "0^*"]\\
    & T_0  & & T_0 & 
  \end{tikzcd},
  \]
  where commutativity follows from the evident equalities $01 \circ 1 = 11 \circ 1, 01 \circ 0 = 00 \circ 0$.
\end{proof}

\begin{intuition}
  Intuitively, the map in \cref{lemma:tetra to square} takes a morphism $f\colon x \to y$ and sends it to the triple $(\id_y, f, \id_x)$, where the pullback specifies the first morphism is the identity on the codomain of $f$ and the third morphism is the identity on the domain of $f$.
\end{intuition}
 
We now have all the pieces in place to present another way to define homotopy equivalences.

 \begin{lemma} \label{lemma:homotopy equivalence tetrahedron}
  Let $T$ be a Segal space, let $x, y$ be objects in $T$, and let $f\colon x \to y$ be a morphism in $T$. Then $f$ is a homotopy equivalence if and only if the following diagram admits a lift $H_2$ or $H_3$
  \[
  \begin{tikzcd}[column sep= 2cm, row sep=0.5cm]
    & & T_3 \arrow[d, "P_{32}", twoheadrightarrow] \\
    & &  T_2 \times^{12^*,01^*}_{T_1} T_2 \arrow[d, "P_{21}", twoheadrightarrow] \\
    \Delta[0] \arrow[urr, dashed, "H_2", bend left = 10] \arrow[uurr, dashed, "H_3", bend left = 10] \arrow[r, "f"] & T_1 \arrow[r, "{(11^* , 01^* , 00^*)}"] & T_1 \times^{1^*,1^*}_{T_0} T_1 \times^{0^*,0^*}_{T_0} T_1
  \end{tikzcd}.
  \]
 \end{lemma}
 \begin{proof}
  By \cref{lemma:pthree equivalence}, the map $P_{32}$ is a trivial fibration of Kan complexes, so it suffices to show that $f$ is a homotopy equivalence if and only if the diagram admits a lift $H_2$. Let $f\colon x \to y$ be a homotopy equivalence. Denote its right inverse by $g$ and its left inverse by $h$. Let $\sigma_g, \sigma_h\colon F(2) \to T$ be the maps corresponding to the fact that $fg$ composes to $\id_y$ and $hf$ composes to $\id_x$. Then $H_2 \coloneq (\sigma_g, \sigma_h)$ is a vertex in $T_2 \times^{12^*,01^*}_{T_1} T_2$ that lifts $(\id_y ,f ,\id_x)$ in $T_1 \times^{1^*,1^*}_{T_0} T_1 \times^{0^*,0^*}_{T_0} T_1$. On the other hand, assume that $(\id_y ,f ,\id_x)$ lifts to $H_2 = (\sigma_1, \sigma_2)$ in $T_2 \times^{12^*,01^*}_{T_1} T_2$. Let us denote $01^*\sigma_1 = g$ and $12^*\sigma_2 = h$. Then $\sigma_1$ is a witness for $fg \simeq \id_y$ and $\sigma_2$ is a witness for $hf \simeq \id_x$, meaning $f$ is a homotopy equivalence. 
\end{proof}

 \begin{intuition} \label{intuition:homotopy equivalence tetrahedron}
  This characterization appears technical, but is diagrammatically intuitive. The element $(\id_y, f, \id_x)$ in $T_1 \times^{1^*,1^*}_{T_0} T_1 \times^{0^*,0^*}_{T_0} T_1$  can be represented by the diagram:
    \[
    \begin{tikzcd}[row sep=0.5cm, column sep=0.5cm]
     & x & \\
     & & \\
     y \arrow[dr, "\id_y"'] & & x \arrow[dl, "f"] \arrow[uul, "\id_x"'] \\
     & y 
    \end{tikzcd}
   \]
   A lift to an element in $T_2 \times_{T_1} T_2$ is a diagram of the following form.
    \[
    \begin{tikzcd}[row sep=0.6cm, column sep=1cm]
     & x & \\
     & & \\
     y \arrow[dr, "\id_y"'] \arrow[rr, ""{name=U, above}, "g" near start, dashed] & & x \arrow[dl, "f"] \arrow[uul, "\id_x"'] \\
     & y \arrow[uuu, dashed, "h" near end]
      \arrow[start anchor={[xshift=10ex, yshift=10ex]}, to=U, phantom, "\sigma_2"]
     \arrow[start anchor={[xshift=-4ex, yshift=2ex]}, to=U, phantom, "\sigma_1"]
    \end{tikzcd}
    \]
  The data in this diagram is exactly that of two morphisms $g,h\colon y \to x$ and two $2$-cells $\sigma_1, \sigma_2$ that give us the right and left inverses. 

  On the other hand, a lift to an element in $T_3$ is a diagram of the following form.
    \[
    \begin{tikzcd}[row sep=0.6cm, column sep=1cm]
     & x & \\
     & & \\
     y \arrow[uur, dashed, "g"] \arrow[dr, "\id_y"'] \arrow[rr, ""{name=U, above}, "g" near start, dashed] & & x \arrow[dl, "f"] \arrow[uul, "\id_x"'] \\
     & y \arrow[uuu, dashed, "h" near end]
      \arrow[start anchor={[xshift=-10ex, yshift=10ex]}, to=U, phantom, "g"]
      \arrow[start anchor={[xshift=10ex, yshift=10ex]}, to=U, phantom, "\sigma_2"]
     \arrow[start anchor={[xshift=-4ex, yshift=2ex]}, to=U, phantom, "\sigma_1"]
    \end{tikzcd}
   \]
   The data in this diagram is precisely the previous data along with the data that the two morphisms $g,h\colon y \to x$ are equivalent (which is witnessed by the front $2$-cell with boundaries $\id_y,g,h$).

   Hence, the fact that the existence of lifts to $T_3$ and $T_2 \times^{12^*,01^*}_{T_1} T_2$ is equivalent (\cref{lemma:pthree equivalence}), corresponds to the fact that if a morphism has a left inverse and a right inverse, then those inverses are homotopic, as we observed in \cref{lemma:inverses homotopic}.
   
 \end{intuition}

This new method gives us an easy way of showing that the notion of a homotopy equivalence is homotopy invariant.
\begin{lemma} \label{lemma:hoequiv arrows}
  Let $\gamma\colon \Delta[1] \to T_1$ be a path from $\gamma(0) = f$ to $\gamma(1) = f'$ in $T_1$. If $f'$ is a homotopy equivalence, then $f$ is also a homotopy equivalence. 
\end{lemma}

\begin{proof}
  By \cref{lemma:homotopy equivalence tetrahedron}, $f'$ admits a lift $H'\colon \Delta[0] \to T_3$. Now, we have the following commutative diagram
  \begin{center}
    \begin{tikzcd}[row sep=0.7in, column sep=0.7in]
      \Delta[0] \arrow[d, " 1 "] \arrow[rr, "H'"] &   & T_3 \arrow[d, "P_{21} \circ P_{32}", twoheadrightarrow] \\
      \Delta[1] \arrow[urr, dashed, "\tilde{\gamma}"'] \arrow[r, "\gamma"'] & T_1 \arrow[r, "{(11^* , 01^* , 00^*)}"] & T_1 \times^{1^*,1^*}_{T_0} T_1 \times^{0^*,0^*}_{T_0} T_1
    \end{tikzcd}.
  \end{center}
  This diagram lifts because the right-hand map is a fibration. Thus $\tilde{\gamma}(0)$ is the lift of the given map $(11^*f,f,00^*f)$, which, again by \cref{lemma:homotopy equivalence tetrahedron}, proves $f$ is a homotopy equivalence. 
\end{proof}
 
 We now want to proceed to show that there are various spaces of equivalences, which are all homotopy equivalent to each other and give us a more intrinsic way to understand homotopy equivalences. We start with the most straightforward one.

\begin{definition}[Space of homotopy equivalences] \label{def:hoequiv space}
  Let $T$ be a Segal space. Let $T_{\hoequiv}$ be the sub-simplicial set of $T_1$ generated by the elements in $T_{10}$ that are homotopy equivalences (\cref{ex:sub simplicial set}). We call $T_{\hoequiv}$ the \emph{space of homotopy equivalences}.
\end{definition}
 
We now immediately have the following observation about $T_{\hoequiv}$.

\begin{lemma} \label{lemma:hoequiv inclusion kan fibration}
  The map $i\colon T_{\hoequiv} \hookrightarrow T_1$ is an inclusion of path-components (\cref{def:inclusion path components}) and a Kan fibration.
\end{lemma}

\begin{proof}
  Let $\sigma$ be in $T_{1n}$. As the vertices $0^*\sigma, ... n^*\sigma$ in $T_{10}$ are all connected by edges in $T_{1n}$, if one of them is a homotopy equivalence, then all of them are homotopy equivalences, by \cref{lemma:hoequiv arrows}. This means that $\sigma$ is in $(T_{\hoequiv})_n$ if and only if $0^*\sigma$ is in $T_{\hoequiv}$, which is equivalent to saying that $\tau_0\sigma = [0^*\sigma]$ lies in $\pi_0(T_{\hoequiv})$. The argument implies that the following is a pullback diagram
  \[
  \begin{tikzcd}
  T_{\hoequiv} \arrow[r, "\tau_0"] \arrow[d, "i", twoheadrightarrow] \arrow[dr, phantom, "\ulcorner", very near start] & \pi_0(T_{\hoequiv}) \arrow[d, "\pi_0(i)", twoheadrightarrow] \\
  T_1  \arrow[r, "\tau_0"] & \pi_0(T_1)
  \end{tikzcd},
  \]  
  meaning that $T_{\hoequiv} \to T_1$ is an inclusion of path-components (\cref{def:inclusion path components}). Finally, by \cref{prop:inclusion path component kan fibration}, it is also a Kan fibration.
\end{proof}

We now use $P_{32}, P_{21}$ to pursue alternative characterizations of the space $T_{\hoequiv}$ and the map $i\colon T_{\hoequiv} \to T_1$, which will give us more insight into homotopy equivalences in Segal spaces.

\begin{definition}[Space of homotopy equivalences with (compatible) inverse] \label{def:hoeqchoice}
  Let $T$ be a Segal space. We define the \emph{space of homotopy equivalences with compatible inverse} $T_{\hoeqcomp}$ and the \emph{space of homotopy equivalences with inverses} $T_{\hoeqchoice}$ via the following pullback diagrams
  \[
  \begin{tikzcd}[column sep=1.5cm]
    T_{\hoeqcomp} \arrow[d, "P_{\forgcomp}"', twoheadrightarrow] \arrow[r] \arrow[dr, phantom, "\ulcorner", very near start] & T_3 \arrow[d, "P_{32}", twoheadrightarrow] \\
    T_{\hoeqchoice} \arrow[d, "P_{\forgchoice}"', twoheadrightarrow] \arrow[r] \arrow[dr, phantom, "\ulcorner", very near start] & T_2 \times^{12^*,01^*}_{T_1} T_2 \arrow[d, "P_{21}", twoheadrightarrow] \\
    T_1 \arrow[r, "{(11^* , 01^* , 00^*)}"] & T_1 \times^{1^*,1^*}_{T_0} T_1 \times^{0^*,0^*}_{T_0} T_1
   \end{tikzcd}
  \]
\end{definition}

\begin{remark} \label{rem:kan fibrations hoeqchoice}
  The maps $P_{\forgchoice}, P_{\forgcomp}$ are Kan fibrations because they are pullbacks of the fibrations $P_{21}, P_{32}$ (\cref{lemma:zthree in fthree}).
\end{remark}

We now have the following immediate observation about $P_{\forgcomp}$.

\begin{lemma} \label{lemma:hoeqcomp hoeqchoice equivalence}
  The map $P_{\forgcomp}\colon T_{\hoeqcomp} \to T_{\hoeqchoice}$ is a trivial fibration of Kan complexes.
\end{lemma}

\begin{proof}
  By \cref{lemma:pthree equivalence}, the map $P_{32}\colon T_3 \to T_2 \times^{12^*,01^*}_{T_1} T_2$ is a trivial Kan fibration, and so the result follows from \cref{lemma:trivial fibration pullback}.
\end{proof}

\begin{remark}
 In light of the explanations in \cref{intuition:homotopy equivalence tetrahedron}, we can understand a point in the space $T_{\hoeqcomp}$ as a homotopy equivalence together with compatible inverses. On the other hand, a point in $T_{\hoeqchoice}$ is a homotopy equivalence together with not necessarily compatible inverses, and the map $P_{\forgcomp}$ is a forgetful map that forgets the compatibility. The fact that $P_{\forgcomp}$ is a trivial fibration equivalently states that for a morphism $(f,g,h, \sigma_g, \sigma_h)\colon \Delta[0] \to T_\hoeqchoice$, not only is there a lift to $T_{\hoeqcomp}$, which corresponds to a compatibility of the inverses (as proven in \cref{lemma:inverses homotopic}), but that there is a homotopically unique way to choose such compatibility data (consistent with the philosophy of \cref{sec:contractibility}).
\end{remark}

We can even simplify the pullback diagram to characterize $T_{\hoeqcomp},T_{\hoeqchoice}$ in a simpler manner.
 
 \begin{lemma} \label{lemma:hoeqchoice second}
 In the commutative diagram below 
 \[
   \begin{tikzcd}[row sep=0.5cm, column sep=1.2cm]
    T_{\hoeqcomp} \arrow[d] \arrow[r, "P_{\forgcomp}", twoheadrightarrow] & T_{\hoeqchoice} \arrow[d] \arrow[r, "P_{\forgchoice}", twoheadrightarrow] \arrow[dr, phantom, "\ulcorner", very near start] & T_1 \arrow[d, "{(11^* , 01^* , 00^*)}"] \arrow[r, "{(1^* , 0^*)}", twoheadrightarrow] & T_0 \times T_0 \arrow[d, "00^* \times 00^*"]\\
    T_3 \arrow[r, "P_{32}", twoheadrightarrow] &  T_2 \times^{12^*,01^*}_{T_1} T_2 \arrow[r, "P_{21}", twoheadrightarrow] & T_1 \times^{1^*,1^*}_{T_0} T_1 \times^{0^*,0^*}_{T_0} T_1 \arrow[r, "{(\pi_1,\pi_3)}", twoheadrightarrow] &  T_1 \times T_1 
   \end{tikzcd}
  \]
  all squares and rectangles are pullback squares.
 \end{lemma}

\begin{proof}
 The middle square and left-hand rectangle are pullbacks by definition. So, we can apply \cref{lemma:pullback pasting} to deduce that the left-hand square is a pullback. Now, again by \cref{lemma:pullback pasting}, all remaining squares and rectangles will be pullbacks as soon as we show that the right-hand square is a pullback. 
 
 However, this follows from \cref{prop:colimit limit map} and the fact that the following is a pushout diagram (\cref{def:pushout square})
 \[
\begin{tikzcd}
 F(1) \coprod F(1) \arrow[r, "01_0 + 01_2"] \arrow[d, "00_0 + 00_1"'] & F(1) \coprod_{F(0)}^{1,1} F(1) \coprod_{F(0)}^{0,0} F(1)\arrow[d, "11 + 01 + 00"] \\
  F(0) \coprod F(0) \arrow[r, "1 + 0"] &  F(1) \arrow[ul, phantom, "\ulcorner", very near start]
\end{tikzcd}.
\]
Here we freely used the Yoneda lemma (\cref{lemma:FDelta yoneda}) to characterize the maps. The subscripts indicate which copy of $F(1)$ in $ F(1) \coprod_{F(0)} F(1) \coprod_{F(0)} F(1)$ we are mapping into.
\end{proof}

We now want to show that both $T_{\hoeqcomp}, T_{\hoeqchoice}$ are equivalent to $T_{\hoequiv}$. For that we need some additional definitions.

\begin{definition}[Space of left and right inverses] \label{def:linv rinv}
  Let $T$ be a Segal space, let $x,y$ be two objects, and let $f\colon x \to y$ be a morphism. We define the \emph{space of left inverses} $\LInv(f)$ and the \emph{space of right inverses} $\RInv(f)$ via the following pullback diagrams
\[
\begin{tikzcd}
\LInv(f) \arrow[r] \arrow[d, twoheadrightarrow] \arrow[dr, phantom, "\ulcorner", very near start] & T_2 \arrow[d, "{(02^*,01^*)}", twoheadrightarrow] \\
\Delta[0] \arrow[r, "{(\id_x , f)}"'] & T_1 \times^{0^*,0^*}_{T_0} T_1
\end{tikzcd}
\begin{tikzcd}
\RInv(f) \arrow[r] \arrow[d, twoheadrightarrow] \arrow[dr, phantom, "\ulcorner", very near start] & T_2 \arrow[d, "{(02^*,12^*)}", twoheadrightarrow] \\
\Delta[0] \arrow[r, "{(\id_y , f)}"'] & T_1 \times^{1^*,1^*}_{T_0} T_1
\end{tikzcd}
\]
\end{definition}

\begin{intuition}
 A point in $\LInv(f)$ corresponds to a morphism $g\colon y \to x$ together with a $2$-cell from $gf$ to $\id_x$, and a point in $\RInv(f)$ corresponds to a morphism $h\colon y \to x$ together with a $2$-cell from $fh$ to $\id_y$. Thus, $\LInv(f)$ and $\RInv(f)$ can be thought of as the spaces of left and right inverses of $f$ respectively.
\end{intuition}

\begin{lemma} \label{lemma:linv rinv pullback}
  Let $T$ be a Segal space, let $x,y$ be two objects, and let $f\colon x \to y$ be a morphism. Then the following is a pullback diagram
\[
\begin{tikzcd}
   \RInv(f) \times \LInv(f) \arrow[r] \arrow[d, twoheadrightarrow] & T_{\hoeqchoice} \arrow[d, "P_{\forgchoice}", twoheadrightarrow] \\
  \Delta[0] \arrow[r, "f"] & T_1
\end{tikzcd}
\]
\end{lemma}

\begin{proof}
  We can extend the diagram to the following diagram
 \[
\begin{tikzcd}
   \RInv(f) \times \LInv(f) \arrow[r] \arrow[d, twoheadrightarrow] & T_{\hoeqchoice} \arrow[d, "P_{\forgchoice}", twoheadrightarrow] \arrow[r] & T_2 \times^{12^*,01^*}_{T_1} T_2 \arrow[d, twoheadrightarrow]\\
  \Delta[0] \arrow[r, "f"] & T_1 \arrow[r]  & T_1 \times^{1^*,1^*}_{T_0} T_1 \times^{0^*,0^*}_{T_0} T_1
\end{tikzcd},
\]
where the right-hand square is a pullback by definition of $T_{\hoeqchoice}$. Hence, by \cref{lemma:pullback pasting}, it suffices to show that the outer rectangle is a pullback. We now have the following diagram
 \[
  \begin{tikzcd}
    \Delta[0] \arrow[r, "{(\id_y,f)}"] \arrow[d, equal] 
      & T_1 \times^{1^*,1^*}_{T_0} T_1 \arrow[d, "\pi_2"] 
      & T_2 \arrow[d, "12^*"] \arrow[l, twoheadrightarrow, "{(02^*,12^*)}"'] 
      &[-.7cm] \RInv(f) \arrow[d] \\
    \Delta[0] \arrow[r, "\{f\}"] 
      & T_1 
      & T_1 \arrow[l, equal] \arrow[r, rightsquigarrow] 
      & \Delta[0] \\
    \Delta[0] \arrow[r, "{(\id_x,f)}"] \arrow[u, equal] 
      & T_1 \times^{0^*,0^*}_{T_0} T_1 \arrow[u, "\pi_2"'] \arrow[d, rightsquigarrow] 
      & T_2 \arrow[u, "01^*"'] \arrow[l, twoheadrightarrow, "{(02^*,01^*)}"'] 
      & \LInv(f) \arrow[u] \\ 
    \Delta[0] \arrow[r, "{(\id_y,f,\id_x)}"] 
      & T_1 \times^{1^*,1^*}_{T_0} T_1 \times^{0^*,0^*}_{T_0} T_1 
      & T_2 \times^{12^*,01^*}_{T_1} T_2 \arrow[l, twoheadrightarrow, "{P_{21}}"']
  \end{tikzcd}.
\]
 We can directly see that the pullback of the columns gives us the bottom row, which is the pullback we want to compute. On the other hand, by \cref{def:linv rinv}, the pullbacks of the rows give us the right column, whose pullback is the product $\RInv(f) \times \LInv(f)$ (\cref{ex:product}). Finally, by pullback commutativity (\cref{lemma:pullback commutativity}), these two coincide, giving us the desired result. 
\end{proof}
 
\begin{intuition}
  Considering the detailed explanation given in \cref{intuition:homotopy equivalence tetrahedron}, a point in the pullback $T_{\hoeqchoice} \times^{P_{\forgchoice},\{f\}}_{T_1} \Delta[0]$ corresponds to a choice of left and right inverses for the chosen morphism $f\colon x \to y$. \cref{lemma:linv rinv pullback} is telling us that these choices of data are fully independent not just at the level of vertices, but for the whole simplicial structure, which is consistent with our intuition regarding left and right inverses. 
\end{intuition}

\begin{lemma} \label{lemma:linv rinv local}
  Let $T$ be a Segal space, let $x,y$ be two objects, and let $f\colon x \to y$ be a morphism. Then the following diagrams are pullback squares
\[
\begin{tikzcd}
\LInv(f) \arrow[r] \arrow[d, twoheadrightarrow] \arrow[dr, phantom, "\ulcorner", very near start] & \map_T(x,y,x) \arrow[d, "{(02^*,01^*)}", twoheadrightarrow] \\
\Delta[0] \arrow[r, "{(\id_x , f)}"'] & \map_T(x,x) \times \map_T(x,y)
\end{tikzcd}
\begin{tikzcd}
\RInv(f) \arrow[r] \arrow[d, twoheadrightarrow] \arrow[dr, phantom, "\ulcorner", very near start] & \map_T(y,x,y) \arrow[d, "{(02^*,12^*)}", twoheadrightarrow] \\
\Delta[0] \arrow[r, "{(\id_y , f)}"'] & \map_T(y,y) \times \map_T(x,y)
\end{tikzcd}
\]
\end{lemma}

\begin{proof}
  We only show the case for $\LInv(f)$ and the other case is analogous. We have the following diagram
\[
\begin{tikzcd}
\LInv(f) \arrow[r] \arrow[d, twoheadrightarrow]  & \map_T(x,y,x) \arrow[d, "{(02^*,01^*)}", twoheadrightarrow] \arrow[r] & T_2 \arrow[d, "{(02^*,01^*)}", twoheadrightarrow]  \\
\Delta[0] \arrow[r, "{(\id_x , f)}"'] & \map_T(x,x) \times \map_T(x,y) \arrow[r] \arrow[d, twoheadrightarrow] & T_1 \times_{T_0}^{0^*,0^*} T_1 \arrow[d, twoheadrightarrow, "{(0^*\pi_1,1^*\pi_2,1^*\pi_1)}"] \\
 & \Delta[0] \arrow[r, "{(x,y,x)}"] & T_0 \times T_0 \times T_0
\end{tikzcd}.
\]
We make the following observations about this diagram. First the right-hand rectangle is a pullback by definition of $\map_T(x,y,x)$. We now move on to show the bottom square is a pullback. We have the following diagram
\[ 
  \begin{tikzcd}
    \Delta[0] \arrow[r, "{(x,y)}"] \arrow[d, equal] & T_0 \times T_0 \arrow[d, "\pi_1"] &[1cm] T_1 \arrow[d, "0^*"]  \arrow[l, twoheadrightarrow, "{(0^*,1^*)}"'] &[-.7cm] \map_T(x,y) \arrow[d] \\
    \Delta[0] \arrow[r, "\{x\}"] & T_0 & T_0 \arrow[l, equal] \arrow[r, rightsquigarrow] & \Delta[0]\\
    \Delta[0] \arrow[r, "{(x,x)}"] \arrow[u, equal] &T_0 \times T_0 \arrow[u, "\pi_2"'] \arrow[d, rightsquigarrow] & T_1 \arrow[u, "0^*"'] \arrow[l, twoheadrightarrow, "{(0^*,1^*)}"'] & \map_T(x,x) \arrow[u] \\ 
    \Delta[0] \arrow[r, "{(x,y,x)}"] & T_0 \times T_0 \times T_0 & T_1 \times^{0^*,0^*}_{T_0} T_1 \arrow[l, twoheadrightarrow, "{(0^*\pi_1,1^*\pi_2,1^*\pi_1)}"'] \\
  \end{tikzcd}.
 \]
 We can directly see that the pullback of the columns gives us the bottom row, which is the pullback of the bottom right square. On the other hand, by \cref{def:linv rinv}, the pullbacks of the rows give us the right column, whose pullback is the product $\map_T(x,x) \times \map_T(x,y)$ (\cref{ex:product}). Thus, by pullback commutativity (\cref{lemma:pullback commutativity}), these two coincide, proving the bottom square is indeed a pullback square. Hence, by \cref{lemma:pullback pasting}, the top right square is also a pullback.
 
 Finally, by \cref{def:linv rinv}, the top rectangle is also a pullback, and so, again by \cref{lemma:pullback pasting}, the left-hand square is also a pullback. Hence we are done.
\end{proof}

\begin{definition} \label{def:map xyz f g}
  Let $T$ be a Segal space, let $x,y,z$ be three objects, and let $f\colon x \to y, g\colon y \to z$ be two morphisms. We define $\map_T(x,y,z)_{f/},\map_T(x,y,z)_{/g}$ as the pullback of the following diagram
 \[
  \begin{tikzcd}
  \map_T(x,y,z)_{f/} \arrow[r] \arrow[d, twoheadrightarrow] \arrow[dr, phantom, "\ulcorner", very near start] & \map_T(x,y,z) \arrow[d, "01^*", twoheadrightarrow] \\
  \Delta[0] \arrow[r, "{f}"'] & \map_T(x,y)
  \end{tikzcd}
  \begin{tikzcd}
  \map_T(x,y,z)_{/g} \arrow[r] \arrow[d, twoheadrightarrow] \arrow[dr, phantom, "\ulcorner", very near start] & \map_T(x,y,z) \arrow[d, "12^*", twoheadrightarrow] \\
  \Delta[0] \arrow[r, "g"'] & \map_T(y,z)
  \end{tikzcd}.
 \]
\end{definition}

We have the following result about these spaces.

\begin{lemma} \label{lemma:done kan fibration}
  Let $T$ be a Segal space, let $x,y,z$ be three objects, and let $f\colon x \to y, g\colon y \to z$ be two morphisms. Then $02^*\colon\map_T(x,y,z)_{f/} \to \map_T(x,z)$, $02^*\colon\map_T(x,y,z)_{/g} \to \map_T(x,z)$ are Kan fibrations.
\end{lemma}

\begin{proof}
  We prove the case for $\map_T(x,y,z)_{f/}$ and the other case is analogous. We have the following commutative diagram 
  \[
  \begin{tikzcd}
   \map_T(x,y,z)_{f/} \arrow[r] \arrow[d, "02^*", twoheadrightarrow] \arrow[dr, phantom, "\ulcorner", very near start] &  \map_T(x,y,z) \arrow[d, "{(02^*,01^*)}", twoheadrightarrow] \arrow[dd, bend left = 45, "01^*", shift left=1.5cm]\\
   \map_T(x,z) \arrow[r] \arrow[d, twoheadrightarrow] & \map_T(x,z) \times \map_T(x,y) \arrow[d, twoheadrightarrow, "\pi_2"] \\
   \Delta[0] \arrow[r, "{f}"'] & \map_T(x,y)
  \end{tikzcd}.
  \]
  The bottom square is a pullback by direct computation, and the rectangle is a pullback square, by \cref{def:map xyz f g}, and so the top square is also a pullback square by \cref{lemma:pullback pasting}. As the right-hand vertical map is a Kan fibration, so is the left-hand vertical map (\cref{lemma:kan fibration pullback}), as desired.
\end{proof}

\begin{lemma} \label{lemma:map xyz f g}
 Let $T$ be a Segal space, let $x,y,z$ be three objects, and let $f\colon x \to y, g\colon y \to z$ be two morphisms. Then the maps $12^*\colon\map_T(x,y,z)_{f/} \to \map_T(y,z)$, $01^*\colon\map_T(x,y,z)_{/g} \to \map_T(x,y)$ are equivalences.
\end{lemma}

\begin{proof}
  We prove the case for $\map_T(x,y,z)_{f/}$ and the other case is analogous. We have the following diagram
  \[
  \begin{tikzcd}
    \map_T(x,y,z) \arrow[r, twoheadrightarrow] \arrow[d, "\simeq", "{(01^*,12^*)}"', twoheadrightarrow] & \map_T(x,y) \arrow[d, equal] & \Delta[0] \arrow[l, "\{f\}"'] \arrow[d, equal] \\
    \map_T(x,y) \times \map_T(y,z) \arrow[r, "\pi_1", twoheadrightarrow] & \map_T(x,y) & \Delta[0] \arrow[l, "\{f\}"']
  \end{tikzcd},
  \]
  where the vertical map is an equivalence by \cref{prop:generalized mapping space mapping space}, and the left horizontal maps are Kan fibrations. Hence, by \cref{lemma:homotopy pullback,rem:homotopy pullback}, the induced map on pullbacks, namely $\map_T(x,y,z)_{f/} \to \map_T(y,z)$ is also an equivalence, as desired.
\end{proof}

\begin{lemma} \label{lemma:linv rinv fiber}
 Let $T$ be a Segal space, let $x,y$ be two objects, and let $f\colon x \to y$ be a morphism. Then the following are pullback diagrams.
\[
\begin{tikzcd}
\LInv(f) \arrow[r] \arrow[d, twoheadrightarrow] \arrow[dr, phantom, "\ulcorner", very near start] & \map_T(x,y,x)_{f/} \arrow[d, "02^*", twoheadrightarrow] \\
\Delta[0] \arrow[r, "{\id_x}"'] & \map_T(x,x)
\end{tikzcd}
\begin{tikzcd}
\RInv(f) \arrow[r] \arrow[d, twoheadrightarrow] \arrow[dr, phantom, "\ulcorner", very near start] & \map_T(y,x,y)_{/f} \arrow[d, "02^*", twoheadrightarrow] \\
\Delta[0] \arrow[r, "{\id_y}"'] & \map_T(y,y)
\end{tikzcd}
\]
\end{lemma}

\begin{proof}
 We prove the case for $\LInv(f)$ and the other case is analogous. We now have the following diagram
 \[ 
  \begin{tikzcd}
   \LInv(f) \arrow[r] \arrow[d, twoheadrightarrow] & \map_T(x,y,x)_{f/} \arrow[r] \arrow[d, twoheadrightarrow] & \map_T(x,y,x) \arrow[d, "{(02^*,01^*)}", twoheadrightarrow] \\
   \Delta[0] \arrow[rr, "{(\{\id_x\},\{f\})}"', bend right=10] \arrow[r, "\{\id_x\}"] &  \map_T(x,x) \arrow[r, "\id \times \{f\}"] & \map_T(x,x) \times \map_T(x,y)
  \end{tikzcd}.
 \]
 By \cref{lemma:linv rinv local} the rectangle is a pullback square, and by \cref{def:map xyz f g} the right square is also a pullback square. Thus, by \cref{lemma:pullback pasting}, the left square is also a pullback square, giving us the desired result.
\end{proof}

We now have the following key observations about $\LInv(f),\RInv(f)$, when $f$ is a homotopy equivalence.

\begin{proposition} \label{prop:linv rinv contractible}
  Let $T$ be a Segal space, let $x,y$ be two objects, and let $f\colon x \to y$ be a homotopy equivalence. Then $\LInv(f)$ and $\RInv(f)$ are contractible Kan complexes.
\end{proposition}

\begin{proof}
 It suffices to consider the case of $\LInv(f)$, and the case for $\RInv(f)$ is analogous. First we observe that the composition $\map_T(x,y,x)_{f/} \to \map_T(y,x) \xrightarrow{f^*} \map_T(x,x)$ is homotopic to the map $02^*\colon \map_T(x,y,x)_{f/} \to \map_T(x,x)$, as both correspond to pre-composition with $f$. However, by \cref{lemma:map xyz f g,prop:pre post composition homotopy equivalence}, the first map is a composition of equivalences and hence an equivalence. This means $02^*\colon \map_T(x,y,x)_{f/} \to \map_T(x,x)$ is also an equivalence. As it is also a Kan fibration (\cref{lemma:done kan fibration}), it is in fact a trivial Kan fibration.
 
 By \cref{prop:trivial fib}, it hence follows that the fiber over $\{\id_x\}$, namely $\LInv(f)$ (\cref{lemma:linv rinv fiber}), is contractible.
\end{proof}

\begin{intuition}
  It is a very classical result in category theory that if a morphism is an isomorphism, then the choice of (right or left) inverse is unique. \cref{prop:linv rinv contractible} is the correct homotopical generalization thereof, proving left and right inverses of equivalences are homotopically unique (\cref{def:homotopicallyunique}).
\end{intuition}

\begin{remark}
  Notice the condition that $f$ is a homotopy equivalence is crucial in the above result. Indeed, this result even fails in classical category theory, such as the category of sets, where every injection admits many left inverses and every surjection admits many right inverses. 
\end{remark}

We now have the following elegant result.

 \begin{theorem} \label{thm:hoeqchoice is hoequiv}
   Let $T$ be a Segal space. The map $T_{\hoeqchoice} \to T_1$ factors through $T_{\hoequiv} \hookrightarrow T_1$ and the resulting map $T_{\hoeqchoice} \to T_{\hoequiv}$ is a homotopy equivalence. 
 \end{theorem}

\begin{proof}
  We first show that the map $T_{\hoeqchoice} \to T_1$ factors through $T_{\hoequiv}$. By definition a point in $T_{\hoeqchoice}$ is a tuple $(f,g,h, \sigma_g,\sigma_h)$, where $\sigma_g,\sigma_h$ witness $g,h$ being right and left inverses of $f$. The map $T_{\hoeqchoice} \to T_1$ takes this tuple to $f$. Since $g,h$ are inverses of $f$, $f$ is a homotopy equivalence and so the map factors through $T_{\hoequiv}$. 
  
  Now, by \cref{rem:kan fibrations hoeqchoice}, $T_{\hoeqchoice} \to T_1$ is a Kan fibration. As $T_{\hoequiv} \to T_1$ is an inclusion of path components, the restricted map $T_{\hoeqchoice} \to T_{\hoequiv}$ remains a Kan fibration (\cref{lemma:kan fibration factorization}). By \cref{prop:trivial fib}, to prove that $T_{\hoeqchoice} \to T_{\hoequiv}$ is a homotopy equivalence, it suffices to show that for every $f$ in $T_{\hoequiv}$, the fiber over $f$ is contractible. By \cref{lemma:linv rinv pullback}, this fiber is equivalent to $\RInv(f) \times \LInv(f)$, which is contractible by \cref{prop:linv rinv contractible}, combined with the fact that $f$ is a homotopy equivalence and the product of contractible Kan complexes is contractible.
\end{proof}

\begin{intuition} \label{intuition:hoeqchoice is hoequiv}
  By definition, a point in $T_{\hoequiv}$ is a map $f\colon x \to y$ in $T$ that has the \emph{property} of being a homotopy equivalence. On the other hand a point in $T_{\hoeqchoice}$ is a map $f\colon x \to y$ in $T$ and choices of inverses $g,h\colon y \to x$ along with homotopies that witness the inverse properties. The map $P_{\forgchoice}\colon T_{\hoeqchoice} \to T_{\hoequiv}$ forgets the specific chosen inverse and only keeps the map. The proof above basically says that if we know a map is a homotopy equivalence, there is only a contractible way of choosing inverses. Hence, following the philosophy of \cref{sec:contractibility}, the choice of inverse is homotopically unique.  

  This observation also relates back to \cref{prop:homotopy equiv iff iso in homotopy category}, where we showed that we can recover being a homotopy equivalence at the level of homotopy categories. The homotopy category cannot see all these homotopies, as it does not remember the mapping space. \cref{thm:hoeqchoice is hoequiv} confirms that this does not cause any loss of information, as we only need the mere existence of a choice of inverse.
\end{intuition}
 
Let us end with the local analogue of the space of homotopy equivalences.
\begin{definition} \label{def:hoequiv space local}
 Let $T$ be a Segal space. For two objects $x, y$  in $T$ we define the \emph{space of homotopy equivalences}, $\hoequiv_{T}(x,y)$, as the following pullback
 \[
 \begin{tikzcd}
  \hoequiv_{T}(x,y) \arrow[r] \arrow[d, twoheadrightarrow] \arrow[dr, phantom, "\ulcorner", very near start] & T_{\hoequiv} \arrow[d, "{(0^* , 1^*)}", twoheadrightarrow] \\
  \Delta[0] \arrow[r, "{(x , y)}"'] & T_0 \times T_0
 \end{tikzcd}
 \]
\end{definition}

We have the following result regarding this definition.

\begin{lemma}
  Let $T$ be a Segal space and let $x,y$ be two objects in $T$. Then the space $\hoequiv_{T}(x,y)$ is isomorphic to the sub-simplicial set of $\map_{T}(x,y)$ generated by the homotopy equivalences (\cref{ex:sub simplicial set}).
\end{lemma}

\begin{proof}
 By \cref{lemma:hoequiv inclusion kan fibration}, the map $T_{\hoequiv} \hookrightarrow T_1$ is an inclusion of path-components, and the map $\hoequiv_{T}(x,y) \to \map_{T}(x,y)$ is given as the pullback of the following diagram
 \[
  \begin{tikzcd}
   T_{\hoequiv} \arrow[d, hookrightarrow] \arrow[r, "{(0^* , 1^*)}", twoheadrightarrow] & T_0 \times T_0 \arrow[d, equal] & \Delta[0] \arrow[l, "{(x,y)}"'] \arrow[d, equal] \\ 
   T_1 \arrow[r, "{(0^* , 1^*)}", twoheadrightarrow]  & T_0 \times T_0 & \Delta[0] \arrow[l, "{(x,y)}"']
  \end{tikzcd}.
 \] 
 This means it again satisfies the pullback condition defining inclusion of path-components (\cref{def:inclusion path components}). 
\end{proof}

\section{Examples of Segal Spaces} \label{sec:example Segal spaces}
We now turn to some examples of Segal spaces. We focus on the case of categories and spaces. In both cases we consider how our categorical notions manifest in those particular examples.
 
\begin{example}[Nerve as a Segal space] \label{ex:nerve}
  Let $\cC$ be a category. Then $\inccat(N\cC)$ is a discrete simplicial space, and hence Reedy fibrant (\cref{prop:inccat reedy fibrant}). Moreover, it satisfies the Segal condition (\cref{prop:nerve Segal}). Thus $\inccat(N\cC)$ is a Segal space.
  
  Fortunately, the definitions we are used to from category theory perfectly match up with the ones for Segal space theory. In particular, an object in the Segal space $\inccat(N\cC)$ is just an object in the category $\cC$. The same is true for morphisms.

  However, as $\inccat(N\cC)_1$ is just a set, the mapping space is actually just a set as well. This in particular implies that composition is well-defined not just up to homotopy. In fact for any collection of objects $x_0, ... , x_n$ in $\inccat(N\cC)$, the space $\map_{\inccat(N\cC)}(x_0,...,x_n)$ is in bijection with $\map_{\inccat(N\cC)}(x_0,x_1) \times ... \times \map_{\inccat(N\cC)}(x_{n-1},x_n)$. Thus the pullback 
  \[
    \begin{tikzcd}
     \Comp_{\inccat(N\cC)}(f,g) \arrow[r] \arrow[d, twoheadrightarrow, "\cong"'] \arrow[dr, phantom, "\ulcorner", very near start] & \map_{\inccat(N\cC)}(x_0,x_1,x_2) \arrow[d, "{(01^*, 12^*)}", twoheadrightarrow] \\
     \Delta[0] \arrow[r, "{(f , g)}"'] & \map_{\inccat(N\cC)}(x_0,x_1) \times \map_{\inccat(N\cC)}(x_1,x_2)
    \end{tikzcd}
  \]
  is not just contractible, but actually just a point.
   
  In addition to all of these, as $\map_{\inccat(N\cC)}(x_0,x_1)$ is just a set, two maps are homotopic if and only if they are equal to each other. This implies that a map is a homotopy equivalence if and only if it is an isomorphism. In particular, the homotopy category of the Segal space, $\Ho(\inccat (N \cC))$ is exactly $\cC$, as the two categories have the same objects and the same morphisms.
 \end{example}

 As expected in the case of an ordinary category, the corresponding Segal space has all the category theory we desire, but has no non-trivial homotopical information. We now see an example in the other direction, where we have a Segal space that has rich homotopical information, but no categorical information.

 \begin{example}[Spaces as Segal spaces] \label{ex:spaces are segal spaces}
  Let $K$ be a space (Kan complex). Our first guess for the corresponding Segal space might be to take the simplicial space $\incspace(K)$. While it does satisfy the Segal condition, it is not Reedy fibrant (\cref{ex:ione not reedy fibrant}). Fortunately, there is an equivalent simplicial space that is Reedy fibrant. Namely, let $\iota K = \Map(\Delta[\bullet],K)$ be the simplicial space defined as 
  \[
  \begin{tikzcd}[column sep=1.2cm]
   K \arrow[r, shorten >=1ex,shorten <=1ex]
   & {\Map(\Delta[1],K)}
   \arrow[l, shift left=1.2, "t"] \arrow[l, shift right=1.2, "s"'] 
   \arrow[r, shift right, shorten >=1ex,shorten <=1ex ] \arrow[r, shift left, shorten >=1ex,shorten <=1ex] 
   & {\Map(\Delta[2],K)} 
   \arrow[l] \arrow[l, shift left=2] \arrow[l, shift right=2] 
   \arrow[r, shorten >=1ex,shorten <=1ex] \arrow[r, shift left=2, shorten >=1ex,shorten <=1ex] 
   \arrow[r, shift right=2, shorten >=1ex,shorten <=1ex]
   & \cdots 
   \arrow[l, shift right=1] \arrow[l, shift left=1] \arrow[l, shift right=3] \arrow[l, shift left=3] 
 \end{tikzcd}
  \]
  where the boundary maps are induced by the maps between the simplices. Alternatively we can see that $\iota K_{k,l} = \Hom_{\sSet}(\Delta[k] \times \Delta[l],K)$. 

  This means that Reedy fibrancy of $\iota K$ corresponds to the fact that the following diagram admits a lift for every $k \geq 0,l \geq 1$ and $0 \leq i \leq l$
  \[
  \begin{tikzcd}
   \displaystyle \partial \Delta[k] \times \Delta[l] \coprod_{\partial \Delta[k] \times \Lambda[l]_i} \Delta[k] \times \Lambda[l]_i \arrow[r] \arrow[d, hook] & K \\
   \Delta[k] \times \Delta[l] \arrow[ur, dashed] & 
  \end{tikzcd}
  \]
  Analogously, applying \cref{rem:segal space lifting}, the Segal condition for $\iota K$ corresponds to the fact that for $k \geq 2$, $l \geq 0$ the following diagram admits a lift
  \[
   \begin{tikzcd}
    \displaystyle \Sp[k] \times \Delta[l] \coprod_{\Sp[k] \times \partial \Delta[l]} \Delta[k] \times \partial \Delta[l] \arrow[r] \arrow[d, hook] & K \\
    \Delta[k] \times \Delta[l] \arrow[ur, dashed] & 
   \end{tikzcd}.
  \]
  These are both conceptually straightforward, yet combinatorially tedious arguments, and so we refer the reader to \cite[Corollary 4.3, Corollary 4.6]{goerssjardine2009simplicialhomotopytheory}.
  
  What does the category theory of a Segal space look like in this case? An object in this Segal space $\iota K$ is a point in $K$. A morphism is now a point in the space $\Map(\Delta[1],K)$ and so is just a path $\Delta[1] \to K$. Composition of morphisms corresponds to concatenation of paths in the space. As concatenation of paths is not unique and only well-defined up to homotopy, here we do need the contractibility condition for composition (\cref{prop:composition contractible}).

  More generally (using \cref{def:path space}), for two given points $x,y$, 
  \[\map_{\iota K}(x,y) = \Delta[0] \times_{(K \times K)} \Map(\Delta[1],K)  = \Path_K(x,y).\]
  So, two morphisms are homotopic if they are homotopic as paths in $K$. In particular, every morphism is a homotopy equivalence, as every path $\Delta[1] \to K$ admits an inverse path (explicitly constructed in \cref{lemma:homotopy equivalence relation} to show homotopy is symmetric). Thus $\iota K$ is an example of a Segal space groupoid (\cref{def:segal space groupoid}).
  
  Notice that the homotopy category of this Segal space is the category which has objects the points $x$ in $K$ and has morphisms homotopy classes of paths in $K$, $\pi_0(\Path_K(x,y))$. This category is commonly called the \emph{fundamental groupoid} of $K$ and is denoted by $\Pi_1(K)$ \cite[Section 2.5]{tomdieck2008algebraictopology}.
 \end{example}

 This last case is a non-example.

 \begin{example}[$G(n)$ is not a Segal space] \label{ex:gn not a segal space}
  For every $n \geq 2$, the simplicial space $G(n)$ (\cref{not:FDelta}) is Reedy fibrant (by \cref{prop:inccat reedy fibrant}), but not a Segal space. Here we explicitly consider the case $n=2$, and the general case is analogous. We can directly see that $G(2)$ is the following simplicial space:
  \[
   \begin{tikzcd}[column sep=1.2cm]
    \{0,1,2 \} \arrow[r, shorten >=1ex,shorten <=1ex]
    & \{00,01,11,12,22 \} 
    \arrow[l, shift left=1.2, "d_1"] \arrow[l, shift right=1.2, "d_0"'] 
    \arrow[r, shift right, shorten >=1ex,shorten <=1ex ] \arrow[r, shift left, shorten >=1ex,shorten <=1ex] 
    & \{ 000,001,011,111,112,122,222 \} 
    \arrow[l] \arrow[l, shift left=2, "d_2"] \arrow[l, shift right=2, "d_0"'] 
    \arrow[r, shorten >=1ex,shorten <=1ex] \arrow[r, shift left=2, shorten >=1ex,shorten <=1ex] 
    \arrow[r, shift right=2, shorten >=1ex,shorten <=1ex]
    & \cdots 
    \arrow[l, shift right=1] \arrow[l, shift left=1] \arrow[l, shift right=3] \arrow[l, shift left=3] 
   \end{tikzcd}
  \]
  where the numbers indicate how the simplicial maps act. Recall that $d_i$ drops the $i$th digit. Thus, we have
  \[G(2)_1 \underset{G(2)_0}{\times} G(2)_1 = \{ (00,00),(00,01),(01,11),(01,12),(11,11),(11,12),(12,22),(22,22)\}\]
  So, clearly $G(2)_2$ is not equivalent to $G(2)_1 \times_{G(2)_0} G(2)_1$ as $G(2)_2$ has $7$ elements and the other has $8$ elements. Concretely, $G(2)_1 \times_{G(2)_0} G(2)_1$ has the element $(01,12)$ which wants to be composed to a $012$ in $G(2)_2$, which is the element in $F(2)_2$ that is missing in $G(2)_2$.
 \end{example}
 
 \section{Complete Segal Spaces} \label{section:css}
 Up until this point we defined Segal spaces and developed their category theory. The comparison with the case of classical categories (\cref{cor:nerve equivalence}) suggests that Segal spaces are a good model for higher categories. However, in this section we will see that in the higher categorical setting we need one additional property, namely completeness. This work is fundamentally based on \cite[Section 3, Section 6]{rezk2001css}.
 
\subsection{Why are Segal Spaces not Enough?} \label{subsec:why css}
By definition, a Segal space has both a category theory and a homotopy theory. However, there is no condition that guarantees that these two theories are compatible with each other. This causes several problems, which we will illustrate in this section. For this first example, let us recall the category $I(1)$ (\cref{ex:ione}). Then, following \cref{ex:nerve}, $\inccat N I(1)$ is a Segal space.

\begin{notation}
  We denote the Segal space $\inccat N I(1)$ by $E(1)$.
\end{notation}

 \begin{example}
  Following \cref{ex:i equivalence}, the functor $0\colon[0] \to I(1)$ is an equivalence of categories. However, the map of Segal spaces $\inccat N [0]  = F(0)\to E(1)$ is not an equivalence of Segal spaces. Indeed, by \cref{rem:equivalence to homotopy equivalence}, if it were the case, then
  \[\{0\}\colon (\inccat N[0])_0 \to (E(1))_0 = \{0,1\}\]
  would be a homotopy equivalence of spaces, which is evidently not the case.
 \end{example}

 \begin{intuition}
  What essentially happened here is that the category theory has an underlying homotopy theory of groupoids ($I(1)$ is a groupoid), which is completely ignored and thus missed by the Segal space.
 \end{intuition}

 \begin{example}
  Let us go back to $E(1)$ once more. It is a discrete Segal space, with two objects denoted $0,1$. Moreover, it has two non-identity morphisms, $01, 10$, which are inverses of each other. Thus the two objects are equivalent to each other in the sense that there is a homotopy equivalence between them. However, they are not equivalent in the space $E(1)_0 = \{0,1\}$, as there is no path between them.
 \end{example}

 \begin{intuition}
  Here we see a clear mismatch between homotopy theory and category theory. Categorically the two points are equivalent, but homotopically they are not. 
 \end{intuition}

 For this next example, we recall equivalent characterizations of equivalence of categories (\cref{thm:equivalence of categories}).
 \begin{example}
  We would expect that a map of Segal spaces $f\colon T \to U$ is an equivalence of Segal spaces if it is fully faithful and essentially surjective. However, this already fails for the simple example $\inccat N [0] \to \inccat N I(1)$.
 \end{example}
 
 \begin{intuition}
  As in the previous example, the problem is that $x$ and $y$ are equivalent in the Segal space $E(1)$, but as points in the space $E(1)_0$ they are not homotopic.
 \end{intuition}
  
 \begin{example} \label{ex:segal space groupoid homotopy hypothesis}
  In \cref{def:segal space groupoid} we defined a Segal groupoid as a Segal space in which every morphism is a homotopy equivalence. In \cref{ex:spaces are segal spaces} we discussed how every space $K$ gives us a Segal groupoid $\iota K$. However, the opposite is not true, meaning it is not the case that for every Segal groupoid $W$, there exists a space $K$ and an equivalence $W\simeq\iota K$. Indeed, the Segal space $E(1)$ is a Segal groupoid, but not equivalent to a space.
 \end{example}
 
\begin{intuition} \label{int:homotopy hypothesis}
 This example is contrary to our understanding of higher category theory. Intuitively, a higher category has homotopical data and categorical data. However, in a groupoid every morphism is invertible, which means it does not contain any non-trivial categorical data. Therefore, our notion of groupoid should really correspond to just a space. The idea we just explained is commonly called the \emph{homotopy hypothesis} and is one of the guiding ideas in the realm of higher category theory, going back to ideas of Grothendieck in the 80s (republished recently in \cite{grothendieck2022pursuingstacks}). 
\end{intuition}

Seeing these examples, we realize that we need to impose one additional condition to make sure the homotopy and category theory suitably coincide.

 \subsection{Defining Complete Segal Spaces} \label{subsec:css}
 In \cref{subsec:why css}, we saw that Segal spaces fall short of our expectations for a good model of higher categories. In particular, we saw that the homotopy theory and category theory of a Segal space are completely disconnected from each other, and we need to impose some compatibility condition. This will result in complete Segal spaces, which is the focus of this section.

\begin{lemma} \label{lemma:identity is hoequiv}
 Let $T$ be a Segal space and let $x$ be an object in $T$. Then the identity morphism $\id_x$ is a homotopy equivalence.
\end{lemma}

\begin{proof}
 Following \cref{prop:associativity and identity of segal space}, $\id_x \circ \id_x \simeq \id_x$, which proves that $\id_x$ is the homotopy inverse of itself.
\end{proof}

\begin{remark} \label{rem:identity is hoequiv}
  \cref{lemma:identity is hoequiv} implies that the image of $00^*\colon T_0 \to T_1$ is actually contained in $T_{\hoequiv}$, meaning there is a factorization $T_0 \xrightarrow{ \ 00^* \ } T_{\hoequiv} \hookrightarrow T_1$. 
\end{remark}

This map $00^*\colon T_0 \to T_{\hoequiv}$ is generally not an equivalence, as the following example shows.
  
\begin{example} \label{ex:eone not complete}
  In the Segal space $E(1)$, the space of objects is the set $\{ 0,1 \}$, but the space of morphisms is the set with four elements $\{ 00, 01, 10, 11 \}$, all of which are equivalences. Thus the map 
  \[00^*\colon \{0,1\} \to \{00, 01, 10, 11\} \]
  is not an equivalence, as it is not surjective.
\end{example}

We saw in several examples of \cref{subsec:why css} that the lack of paths between equivalent objects is the source of all the problems. We therefore introduce the next definition to address this issue.
 
 \begin{definition}[Complete Segal space] \label{def:css}
  A \emph{complete Segal space} is a Segal space $W$ for which the map 
  \[00^*\colon W_0 \to W_{\hoequiv}\] 
  defined in \cref{rem:identity is hoequiv} is an equivalence.
 \end{definition}
 
We now want to characterize the completeness condition in several equivalent ways. Before we can do that, we need to make several technical observations and constructions.

 \begin{remark} \label{rem:lift}
  Given a Kan fibration $p\colon Y \twoheadrightarrow X$ and a Kan complex $K$, every commutative square of the form 
  \[ 
  \begin{tikzcd}
   K \arrow[r] \arrow[d, "00^*"] & Y \arrow[d, "p"] \\
   {\Map(\Delta[1],K)} \arrow[r] \arrow[ur, dashed] & X
   \end{tikzcd}
  \]
  admits a lift. For a proof see \cite[First paragraph of Section I.9, Theorem 11.3]{goerssjardine2009simplicialhomotopytheory}.
 \end{remark}

 \begin{construction} \label{const:hoequiv path}
   Let $W$ be a Segal space. Then following \cref{rem:lift}, the following diagram admits a lift
   \[
   \begin{tikzcd}
    W_0 \arrow[r] \arrow[d] & W_{\hoequiv} \arrow[d, twoheadrightarrow] \\ 
    {\Map(\Delta[1],W_0)} \arrow[r] \arrow[ur, dashed, "P" description] & W_0 \times W_0
   \end{tikzcd},
   \]
   that we denote $P\colon \Map(\Delta[1],W_0) \to W_{\hoequiv}$. 
  \end{construction}
 
\begin{intuition}
  Intuitively, this map takes a path $\gamma\colon\Delta[1] \to W_0$ from $x$ to $y$ and produces a homotopy equivalence from $x$ to $y$.
\end{intuition}

\begin{definition} \label{def:hoequiv path}
  Let $W$ be a Segal space and let $x,y$ be two objects in $W$. Let 
  \[P_{x,y}\colon \Path_{W_0}(x,y) \to \hoequiv_W(x,y)\]
  be given by the pullback of the following diagram
  \[
  \begin{tikzcd}
  \Delta[0] \arrow[r,"{(x,y)}"] \arrow[d, equal] & W_0 \times W_0 \arrow[d, equal] & \Map(\Delta[1],W_0) \arrow[l, twoheadrightarrow] \arrow[d, "P"]\\
  \Delta[0] \arrow[r,"{(x,y)}"] & W_0 \times W_0 & W_{\hoequiv} \arrow[l, twoheadrightarrow]
  \end{tikzcd}
  \]
\end{definition}

\begin{theorem}[{\cite[Theorem 6.2]{rezk2001css}}] \label{thm:hoequiv is Eone}
  Let $W$ be a Segal space. Then there is a homotopy equivalence $W_{\hoequiv} \simeq \Map(E(1),W)$.
\end{theorem}

\begin{proposition} \label{prop:completeness conditions}
  Let $W$ be a Segal space. The following are equivalent.
  \begin{enumerate}
   \item $W$ is a complete Segal space.
   \item Given the commutative diagram below,
   \[
   \begin{tikzcd}[column sep=1.5cm]
    W_0 \arrow[r, "0000^*"] \arrow[d, "00^*"'] & W_3 \arrow[d, "{(02^*,12^*,13^*)}", twoheadrightarrow] \\
    W_1 \arrow[r, "{(11^*,01^*,00^*)}"] & W_1 \times^{1^*,1^*}_{W_0} W_1 \times^{0^*,0^*}_{W_0} W_1
   \end{tikzcd},
   \]
   the induced map $W_0 \to W_{\hoeqcomp}$ (\cref{def:hoeqchoice}) is a homotopy equivalence.
   \item Given the commutative diagram below,
   \[
   \begin{tikzcd}[column sep=1.5cm]
    W_0 \arrow[r, "{(000^*,000^*)}"] \arrow[d, "{00^*}"'] & W_2 \times^{12^*,01^*}_{W_1} W_2 \arrow[d, "P_{21}", twoheadrightarrow] \\
    W_1 \arrow[r, "{(11^*,01^*,00^*)}"] & W_1 \times^{1^*,1^*}_{W_0} W_1 \times^{0^*,0^*}_{W_0} W_1
   \end{tikzcd},
   \]
   the induced map $W_0 \to W_{\hoeqchoice}$ (\cref{def:hoeqchoice}) is a homotopy equivalence.
   \item The map of spaces 
   \[0^*\colon \Map(E(1),W) \to \Map(F(0),W),\]
   induced by the inclusion $0\colon F(0) \to E(1)$, is a weak equivalence.
   \item For any two objects $x,y$ the natural map (\cref{def:hoequiv path})
   \[P_{x,y}\colon \Path_{W_0}(x,y) \to \hoequiv_W(x,y)\]
   is an equivalence of spaces.
  \end{enumerate}
 \end{proposition}
 
 \begin{proof}
  We consider various equivalent conditions.
  
  \emph{$(1 \Longleftrightarrow 2 \Longleftrightarrow 3)$} 
  We have the following diagram
  \[
  \begin{tikzcd}[column sep=1.5cm]
    & W_0 \arrow[dl] \arrow[d] \arrow[dr] & \\
    W_{\hoeqcomp} \arrow[r, "P_{\forgcomp}"', twoheadrightarrow, "\simeq"] & W_{\hoeqchoice} \arrow[r, "P_{\forgchoice}"', twoheadrightarrow, "\simeq"] & W_{\hoequiv}
  \end{tikzcd}
  \]
  By \cref{thm:hoeqchoice is hoequiv,lemma:hoeqcomp hoeqchoice equivalence}, the horizontal maps are trivial Kan fibrations. Thus, the equivalence of $(1)$, $(2)$ and $(3)$ follows from applying $2$-out-of-$3$ (\cref{rem:two out of three}).
  
  \emph{$(1 \Longleftrightarrow 4)$}
  By \cref{thm:hoequiv is Eone}, there is an equivalence $\Map(E(1),W) \simeq W_{\hoequiv}$. Thus, by $2$-out-of-$3$ (\cref{rem:two out of three}), the map $\Map(E(1),W) \to \Map(F(0),W)$ is an equivalence if and only if $W_{\hoequiv} \to W_0$ is an equivalence, which is exactly the completeness condition.
  
  \emph{$(1 \Longleftrightarrow 5)$}
  By \cref{const:hoequiv path}, the map $W_0 \to W_{\hoequiv}$ factors as $W_0 \to \Map(\Delta[1],W_0) \to W_{\hoequiv}$. Thus, by $2$-out-of-$3$ (\cref{rem:two out of three}), the map $W_0 \to W_{\hoequiv}$ is an equivalence if and only if the map $P\colon\Map(\Delta[1],W_0) \to W_{\hoequiv}$ is an equivalence. 

  Now, again by \cref{const:hoequiv path}, we have the following diagram 
  \[
  \begin{tikzcd}
    {\Map(\Delta[1],W_0)} \arrow[dr, twoheadrightarrow] \arrow[rr] & & W_{\hoequiv} \arrow[dl, twoheadrightarrow] \\
    & W_0 \times W_0 
  \end{tikzcd}.
  \]
  Here, the left diagonal map is a Kan fibration, by \cref{prop:kan fibration mapping space}, and the right diagonal map is a Kan fibration by \cref{lemma:kan fibration compose}, as it is the composition of the Kan fibration $W_{\hoequiv} \to W_1$ (\cref{lemma:hoequiv inclusion kan fibration}) and $W_1 \to W_0 \times W_0$ (\cref{ex:reedy fibrant one}). Thus, by \cref{lemma:homotopy pullback preimage}, the top map is a homotopy equivalence if and only if, for every pair of points $(x,y)$, the map between fibers is an equivalence. However, by \cref{def:hoequiv path}, this is precisely the map 
  \[ P_{x,y}\colon \Path_{W_0}(x,y) \to \hoequiv_W(x,y).\]
 \end{proof}

 \begin{intuition} \label{int:css}
  The completeness condition exactly addresses the problems we raised in \cref{subsec:why css}. By imposing the additional condition that every equivalence in the Segal space $W$ can be represented by a path in $W_0$, we are ensuring that the homotopy theory of the space $W_0$ corresponds to equivalences in $W_1$. 
\end{intuition}
 
We can now see how the completeness condition addresses the problems we raised in \cref{subsec:why css}. First we have the following result, the proof of which can be found in \cite{rezk2001css}.

\begin{theorem}[{\cite[Proposition 7.6]{rezk2001css}}]
  Let $F\colon W \to V$ be a functor of complete Segal spaces. The following are equivalent:
  \begin{enumerate}
   \item $F$ is an equivalence.
   \item $F$ is fully faithful and essentially surjective.
  \end{enumerate}
\end{theorem}

We now have a similar result regarding the homotopy hypothesis (\cref{int:homotopy hypothesis}).

\begin{definition}[Complete Segal space groupoid] \label{def:complete segal space groupoid}
  A \emph{complete Segal space groupoid} is a complete Segal space $W$ such that every morphism in $W$ is a homotopy equivalence.
\end{definition}

In \cref{ex:segal space groupoid homotopy hypothesis}, we saw that not every Segal space groupoid is equivalent to a Segal space $\iota K$, for some space $K$. We now want to see how the completeness condition rectifies this problem.

\begin{example} \label{ex:space complete}
  Let $K$ be a space. Then the Segal space $\iota K$ is a complete Segal space groupoid. Indeed, we already saw in \cref{ex:spaces are segal spaces} that $\iota K$ is a Segal space groupoid. Moreover, we saw in the aforementioned example that $\map_{\iota K}(x,y) \simeq \Path_{K}(x,y)$ for any $x,y$ in $K$. Hence in this case the map $P_{x,y}\colon \Path_{K}(x,y) \to \hoequiv_{\iota K}(x,y)$ is the identity, and so $\iota K$ is complete, by \cref{prop:completeness conditions}.
\end{example}

On the other hand, we have the following major result, which confirms the homotopy hypothesis we discussed in \cref{int:homotopy hypothesis}.
 
 \begin{proposition}
  A complete Segal space $W$ is a complete Segal space groupoid if and only if $W$ is equivalent to $\iota W_0$. Hence complete Segal spaces satisfy the homotopy hypothesis.
 \end{proposition}

 \begin{proof}
  One direction follows from \cref{ex:space complete}. Now, let us assume that $W$ is a complete Segal space groupoid. We want to show that $W$ is equivalent to $\iota W_0$. By assumption $W_{\hoequiv} = W_1$. Thus, by the completeness condition, $00^*\colon W_0 \to W_1$ is an equivalence. By the simplicial identities $0^*00^* = 1^*00^* = \id_{W_0}$. So, by $2$-out-of-$3$ (\cref{rem:two out of three}), the maps $0^*,1^*\colon W_1 \to W_0$ are also equivalences, and in fact trivial fibrations (\cref{ex:reedy fibrant one}), as $W$ is Reedy fibrant.

  We now have the following pullback diagram
  \[
  \begin{tikzcd}
   W_1 \times_{W_0} W_1 \arrow[r, "\pi_1", twoheadrightarrow, "\simeq"'] \arrow[d, "\pi_2", twoheadrightarrow, "\simeq"'] \arrow[dr, phantom, "\ulcorner", very near start] & W_1 \arrow[d, "1^*", twoheadrightarrow, "\simeq"'] \\
   W_1 \arrow[r, twoheadrightarrow, "0^*", "\simeq"'] & W_0
  \end{tikzcd},
  \]
  and so the map $W_1 \times_{W_0} W_1 \xrightarrow{\pi_1} W_1 \xrightarrow{0^*}W_0$, or equivalently $0^*\colon \Map(G(2),W) \to W_0$, is a trivial fibration (\cref{lemma:trivial fibration pullback}). Repeating this process inductively, we can conclude that the map $0^*\colon\Map(G(n),W) \to W_0$ is a trivial Kan fibration.

  Finally, by the Segal condition, the composition $W_n \xrightarrow{ i_n^*} \Map(G(n),W) \xrightarrow{0^*} W_0$, which composes to $0^*\colon W_n \to W_0$, is a trivial Kan fibration. However, by \cref{ex:map kan complex equivalence}, we also have the equivalence $W_0 \simeq \Map(\Delta[n],W_0)$. Combining these, we get the desired equivalence $W_n \simeq W_0 \simeq \Map(\Delta[n],W_0) = (\iota W_0)_n$.
 \end{proof}

 \subsection{Rezk Nerves as Complete Segal Spaces} \label{subsec:nerves as css}
 Our final goal is to discuss how we can build a complete Segal space out of a category. We have seen already that $\inccat N \cC$ is a Segal space (\cref{ex:nerve}), but is not necessarily complete (\cref{ex:eone not complete}). Before we can rectify the situation let us better comprehend the problem at hand. 
 
\begin{remark} \label{rem:nerve not complete}
 We saw in \cref{ex:nerve} that objects and morphisms in $\inccat N\cC$ are the objects and morphisms of $\cC$ and a morphism in the Segal space $\inccat N\cC$ is a homotopy equivalence if and only if it is an isomorphism in $\cC$. Hence, for two objects $x,y$ in $\cC$, $\hoequiv_{\inccat N\cC}(x,y) \cong \Hom_{\cC^\simeq}(x,y)$. Here we used the fact that the morphisms in the underlying groupoid $\cC^\simeq$ are precisely the isomorphisms in $\cC$ (\cref{prop:core}).

 Given this explanation and the content of \cref{prop:completeness conditions}, the completeness condition means $\Path_{(\inccat N\cC)_0}(x,y) \simeq \Hom_{\cC^\simeq}(x,y)$. However, $(\inccat N\cC)_0 = \Obj_\cC$ is just the set of objects of $\cC$, and so 
 \[
  \Path_{(\inccat N\cC)_0}(x,y)= \begin{cases} \emptyset & x \neq y \\ \{ \id_x \} & x=y \end{cases}.
 \]
\end{remark}
 
 What we need is a suitably adapted version of the nerve construction that not only adds objects and compositions in the horizontal directions, but also isomorphisms in the vertical direction. This leads to the following construction, which is due to Rezk \cite{rezk2001css}.

 For this definition, we use the evident inclusion $I\colon \DD \to \Cat$ (\cref{def:delta}), the opposite category $(-)^{\op}$ (\cref{lemma:opposite functor}), the functor category (\cref{prop:functor category functorial}), the core of a category (\cref{prop:core}), and the nerve of a category (\cref{def:nerve}).

 \begin{definition}[Rezk nerve] \label{def:rezk nerve}
  Let $\cC$ be a category. We define the \emph{Rezk nerve} of $\cC$ as the simplicial space $\sN \cC$ given by the following composition 
  \[ \DD^{\op} \xrightarrow{ I^{\op} } \Cat^{\op} \xrightarrow{\Fun(-,\cC)} \Cat \xrightarrow{ (-)^\simeq} \Cat \xrightarrow{N} \sSet.\] 
 \end{definition}

 \begin{remark} \label{rem:classification diagram}
  Unlike most other definitions in this section, the Rezk nerve is originally called the \emph{classification diagram} in \cite{rezk2001css}. Given the non-descriptive nature of this name, and the recent prevalence of the terminology ``Rezk nerve'' in the literature, we have opted to use the modern name in this text.
 \end{remark}

 Before we proceed we want to understand this definition more explicitly.

 \begin{remark} \label{rem:rezk nerve explicit}
  By construction, $\sN \cC_n = N(\Fun([n],\cC)^\simeq)$, meaning $\sN \cC_{nl} = N(\Fun([n],\cC)^\simeq)_l$.
\end{remark}

 Let us see an even more explicit way to characterize the Rezk nerve of a category.

\begin{intuition}
  Following \cref{rem:rezk nerve explicit}, $\sN\cC_0 = N(\Fun([0],\cC)^\simeq) = N(\cC^\simeq)$, where we used \cref{ex:functor category from point}. Thus two objects are connected via a path in $\sN\cC_0$ precisely when the objects are isomorphic in $\cC$. This is already an improvement over what we saw in \cref{rem:nerve not complete}.
\end{intuition}

 \begin{lemma} \label{lemma:rezk nerve as functors}
  Let $\cC$ be a category. Then we have the following bijection
    \[\sN(\cC)_{nl} \cong \Hom_{\Cat}([n] \times I(l), \cC).\]
 \end{lemma}

 \begin{proof}
  By \cref{rem:rezk nerve explicit}, $\sN(\cC)_{nl} = N(\Fun([n],\cC)^\simeq)_l \cong \Hom_{\Cat}([l], \Fun([n],\cC)^\simeq)$. However, every morphism in the core is actually an isomorphism, so this canonically coincides with functors  $I(l) \to \Fun([n],\cC)^\simeq$ via the lift (\cref{rem:n to in lift}). Finally, by \cref{prop:functor category as currying}, this is naturally in bijection with the set of functors $[n] \times I(l) \to \cC$.
 \end{proof}
  
  \begin{intuition} \label{int:classification diagram}
   With \cref{lemma:rezk nerve as functors} at hand, we can depict $\sN\cC$ very explicitly. Here, we use the following notation.
   \begin{enumerate}
    \item $\bullet$: objects
    \item $\longrightarrow$: morphisms
    \item $ \xrightarrow{ \ \ \sim \ \ }$: isomorphisms
   \end{enumerate}
   With this notation the Rezk nerve has the form of the following simplicial space:
 \[
 \nervediagram
 \]
 As we can see, the vertical direction (the homotopical direction) focuses on the isomorphisms and so the homotopical data of a category, whereas the horizontal direction (the categorical direction) focuses on arbitrary morphisms and its category theory.
 \end{intuition}

\begin{remark}
  The Rezk nerve is a special case of a more general construction, known as the \emph{relative nerve} or \emph{classification diagram}. See \cite[Section 3]{rezk2001css} for more details.
\end{remark}

Having improved our definition of a nerve, we now proceed to show it is a complete Segal space. The first step, Reedy fibrancy, is somewhat tedious and so we provide a proper reference.

\begin{lemma}[{\cite[Lemma 3.9]{rezk2001css}}]
 Let $\cC$ be a category. Then $\sN(\cC)$ is Reedy fibrant.
\end{lemma}

\begin{lemma} \label{lemma:rezk nerve segal}
  Let $\cC$ be a category. Then the Rezk nerve $\sN(\cC)$ is a Segal space.
\end{lemma}

\begin{proof}
 Let us fix $n \geq 2, l \geq 0$. Following \cref{lemma:rezk nerve as functors}, it suffices to prove the map 
 \[
  \begin{tikzcd}
   \Hom_{\Cat}([n] \times I(l),\cC) \arrow[d] \\ \Hom_{\Cat}([1] \times I(l),\cC) \times_{\Hom_{\Cat}([0] \times I(l),\cC)} ... \times_{\Hom_{\Cat}([0] \times I(l),\cC)} \Hom_{\Cat}([1] \times I(l),\cC)
  \end{tikzcd} 
 \]
 is a bijection. By \cref{prop:functor category as currying}, this map coincides with the map 
 {\small
 \[
 \begin{tikzcd}
  \Hom_{\Cat}([n], \Fun (I(l),\cC)) \arrow[d] \\ \Hom_{\Cat}([1], \Fun(I(l),\cC)) \times_{\Hom_{\Cat}([0], \Fun(I(l),\cC))} ... \times_{\Hom_{\Cat}([0], \Fun(I(l),\cC))} \Hom_{\Cat}([1], \Fun(I(l),\cC))
 \end{tikzcd} 
 \]}
 However, this map is a bijection if $N(\Fun(I(l),\cC))$ satisfies the Segal condition, which is true by \cref{prop:nerve Segal}.
\end{proof}

\begin{remark} \label{rem:objects morphisms rezk nerve}
  Following \cref{lemma:rezk nerve as functors}, an object in the Segal space $\sN(\cC)$ is an element in $\sN(\cC)_{00} = \Fun([0] \times I(0),\cC) \cong \Obj_{\cC}$. Similarly, a morphism in $\sN(\cC)$ is an element in $\sN(\cC)_{10} = \Fun([1] \times I(0),\cC) \cong \Mor_{\cC}$.
\end{remark}

More generally, we have the following.

\begin{lemma} \label{lemma:homotopic morphisms are equal}
  Let $\cC$ be a category, and let $x,y$ be two objects in the Segal space $\sN(\cC)$. Then the mapping space $\map_{\sN(\cC)}(x,y)$ is isomorphic to the constant space on the set $\Hom_{\cC}(x,y)$.
\end{lemma}

\begin{proof}
  By \cref{rem:rezk nerve explicit,def:mapping space Segal}, the mapping space $\map_{\sN(\cC)}(x,y)$ is the pullback of the following diagram
  \[
  \begin{tikzcd}
    \map_{\sN(\cC)}(x,y) \arrow[d, twoheadrightarrow] \arrow[r] & N(\Fun([1]{,}\cC)^\simeq) \arrow[d, "{(0^*,1^*)}", twoheadrightarrow] \\
    \Delta[0] \arrow[r, "{(x,y)}"] & N\cC^\simeq \times N\cC^\simeq
  \end{tikzcd}.
  \]
  We already saw in \cref{rem:objects morphisms rezk nerve} that $\map_{\sN(\cC)}(x,y)_0 \cong \Hom_{\cC}(x,y)$. Moreover, we have 
  \[\map_{\sN(\cC)}(x,y)_l = \{F\colon [1] \times I(l) \to \cC \mid F|_{\{0\} \times I(l)} = \{x\}, F|_{\{1\} \times I(l)} = \{y\}\}, \]
  where $\{x\}, \{y\}$ are the constant functors $I(l) \to \cC$ with value $x$, resp. $y$. A functor $F\colon [1] \times I(l) \to \cC$ satisfying this condition necessarily factors through $[1] \times I(0)$, which means it is determined by a morphism in $\cC$ from $x$ to $y$. Thus we have $\map_{\sN(\cC)}(x,y)_l \cong \Hom_{\cC}(x,y)$ for every $l$, which means $\map_{\sN(\cC)}(x,y)$ is the constant space on the set $\Hom_{\cC}(x,y)$. 
\end{proof}

\begin{intuition}
  An element in $\map_{\sN(\cC)}(x,y)_1$ is a functor $F\colon [1] \times I(1) \to \cC$ of the following form 
  \[ 
  \begin{tikzcd}
    x \arrow[r, "f"] \arrow[d, "\id_x"'] & y \arrow[d, "\id_y"] \\
    x \arrow[r, "g"'] & y
  \end{tikzcd},
  \]
  which implies $\id_y \circ f = g \circ \id_x$, which means $f=g$. 
\end{intuition}

We can use a similar argument to prove the following analogous result.

\begin{lemma} \label{lemma:homotopic morphisms are equal iso}
 Let $\cC$ be a groupoid. Then for two objects $x,y$ in the Kan complex $N\cC$, $\Path_{N\cC}(x,y)$ is isomorphic to the constant space on the set $\Hom_{\cC}(x,y)$.   
\end{lemma}

\begin{proposition}
  Let $\cC$ be a category. Then $\Ho(\sN(\cC))$ is isomorphic to $\cC$.
\end{proposition}

\begin{proof}
  By \cref{constr:homotopy category of segal space}, the objects of $\Ho(\sN(\cC))$ are the objects of $\sN(\cC)$, which, following \cref{rem:objects morphisms rezk nerve}, are the objects of $\cC$. Now, the morphisms in $\Ho(\sN(\cC))$ are homotopy classes of morphisms in $\sN(\cC)$. However, by \cref{lemma:homotopic morphisms are equal}, those are precisely the morphisms in $\cC$. Thus $\Ho(\sN(\cC))$ is isomorphic to $\cC$.
\end{proof}

\begin{lemma} \label{lemma:equivalences in rezk nerve}
  Let $\cC$ be a category. Then a morphism in $\sN(\cC)$ is a homotopy equivalence if and only if it is an isomorphism in $\cC$.
\end{lemma}

\begin{proof}
 By \cref{prop:homotopy equiv iff iso in homotopy category}, a morphism in $\sN(\cC)$ is a homotopy equivalence if and only if it is an isomorphism in $\Ho(\sN(\cC)) = \cC$. Hence we are done.
\end{proof}

\begin{theorem}[Rezk nerve is a complete Segal space]
  Let $\cC$ be a category. Then $\sN(\cC)$ is a complete Segal space.
\end{theorem}

\begin{proof}
  We have already shown it is a Segal space (\cref{lemma:rezk nerve segal}), so it suffices to check the completeness condition.   Let $x,y$ be two objects in $\sN(\cC)$. By \cref{lemma:homotopic morphisms are equal iso}, we have a bijection of constant spaces
  \[\Path_{N\cC^\simeq}(x,y) \cong \Hom_{\cC^\simeq}(x,y).\]
  On the other hand, in \cref{lemma:homotopic morphisms are equal}, we showed that we have the equality of constant simplicial sets
  \[\map_{\sN\cC}(x,y) = \Hom_{\cC}(x,y). \]
  Moreover, by \cref{lemma:equivalences in rezk nerve}, homotopy equivalences in $\sN(\cC)$ are precisely the isomorphisms in $\cC$. Hence, we similarly have an equality of constant simplicial sets 
   \[\hoequiv_{\sN\cC}(x,y) = \Hom_{\cC^\simeq}(x,y),\]
   where $\Hom_{\cC^\simeq}(x,y)$ is the constant simplicial set on the set of isomorphisms from $x$ to $y$ in $\cC$. Thus, following \cref{prop:completeness conditions}, the completeness condition translates to 
   \[P_{x,y}\colon\Hom_{\cC^\simeq}(x,y) \to \Hom_{\cC^\simeq}(x,y)\]
    being a bijection of constant spaces. Note that, under this correspondence, the map \(P_{x,y}\) sends an isomorphism from $x$ to $y$ to the same isomorphism regarded as a homotopy equivalence in $\sN(\cC)$. Hence $P_{x,y}$ is the identity map on \(\Hom_{\cC^\simeq}(x,y)\), which is indeed a bijection.
\end{proof}

\section{Where to go from here?} \label{sec:next}
What have we learned until now? 
\begin{enumerate}[leftmargin=*]
  \item We defined complete Segal spaces as a model of $\infty$-categories.
  \item We showed that they come with many categorical notions (objects, morphisms, composition).
  \item We showed they also come with many homotopical notions (homotopies, equivalences).
  \item We also saw that with an adjusted notion of uniqueness, namely homotopical uniqueness or contractibility, we get many analogues of uniqueness results in category theory, such as the uniqueness of compositions or inverses.
  \item We finally saw that both categories and spaces give us examples of complete Segal spaces.
\end{enumerate}

The next natural step is to properly develop category theory in this setting (limits, adjunctions, Kan extensions, ...), and see some examples. The following is not an exhaustive list of all references, but a selection based on reasonable next steps that an interested reader might want to pursue.

\subsection{Historical Approaches}
As with any area of mathematics, $\infty$-category theory has slowly evolved over time, and there are various historical developments that might not be a starting point, but are worth keeping in mind.
\begin{enumerate}[leftmargin=*]
    \item \textbf{A Simplicially-Enriched Approach to Higher Categories:} As we observed in the definitions of \cref{sec:twopaths} and the proofs of \cref{sec:segal spaces,sec:homotopy equiv in Segal Spaces}, the existence of mapping spaces plays a constructive role in the study of higher categories. This already plays a prominent role in the work of Dwyer and Kan \cite{dwyerkan1980simplocalization,dwyerkan1980calculatingsimplocalizations}.  
    \item \textbf{A Model Categorical Approach to Higher Categories:} On several occasions throughout this note, we observed that various notions of fibrations and lifting conditions played an important role in various higher categorical definitions. The idea of studying categories along with a class of ``fibrations'' that is characterized via certain lifting conditions goes back at least to Quillen, who introduced the notion of a \emph{model category} \cite{quillen1967modelcats}. It is hence unsurprising that many original approaches to higher categories involved model categorical techniques, some of which we encountered throughout this note. This has resulted in ``models of $\infty$-categories''. This includes the model of Segal categories, due to Hirschowitz and Simpson \cite{hirschowitzsimpson1998segalcat}, complete Segal spaces, due to Rezk \cite{rezk2001css}, Kan-enriched categories, due to Bergner \cite{bergner2007bergnermodelcat}, quasi-categories, due to Joyal \cite{joyal2008notes,joyal2008theory}, relative categories, due to Barwick and Kan \cite{barwickkan2012relativecategory}, among others.
    \item \textbf{Comparison with Other Models:} In \cref{thm:sing equivalence}, we saw that Kan complexes and CW-complexes give us equivalent homotopy categories, suggesting we should think of them as equivalent ways to characterize the homotopy theory of spaces. In a similar vein, with the development of models of $\infty$-categories, figures such as Joyal, Tierney, Bergner, Barwick, and Kan have proven that all of these aforementioned models do give us equivalent homotopy categories \cite{joyaltierney2007qcatvssegal,bergner2007threemodels,bergner2010survey,barwickkan2012relativecategory,bergner2018book}.
    \item \textbf{An Axiomatic Approach to Higher Categories:} Given the variety of models, there have been some historic efforts to introduce setups that capture $\infty$-categories axiomatically. This was successfully pursued by Toën \cite{toen2005unicity} and later generalized by Barwick and Schommer-Pries \cite{barwickschommerpries2021unicity}.
\end{enumerate}

\subsection{From Research to Education}
In recent years we have seen the rise of expositional sources for $\infty$-category theory and research sources that contain an expositional component. These have established themselves as a suitable modern resource for learning $\infty$-category theory, and can be the next step for an interested reader after this note.
\begin{enumerate}[leftmargin=*]
  \item \textbf{A Simplicial Approach to Higher Categories:} The simplicial approach to $\infty$-category theory can be seen as a natural generalization of the perspective developed here. It assumes a functioning theory of simplicial objects (primarily simplicial sets, and sometimes simplicial spaces), and then develops $\infty$-categories as a special class of these simplicial objects (quasi-categories as simplicial sets, complete Segal spaces as simplicial spaces). There is now a wide range of introductory and advanced material centering this perspective. Introductory notes include the work of Groth \cite{groth2010highercat} and Antol{\'i}n Camarena \cite{antolincamarena2016whirlwind}, and more recently Ravenel \cite{ravenel2023inftycategory}. More advanced introductions via simplicial sets include the notes by Rezk \cite{rezk2022qcats}, Haugseng \cite{haugseng2017qcats}, and Cisinski \cite{cisinski2019highercategories}, and a textbook by Land \cite{land2021infinitycategories}. More advanced sources based on a theory of simplicial sets include the collected work of Lurie \cite{lurie2009htt,lurie2017ha,lurie2026kerodon}, whose books combine expository and research aspects.
  \item \textbf{A $2$-Categorical Approach to Higher Categories:} One major historical advancement in classical category theory was the realization that many important aspects of categories (such as limits and the Yoneda lemma, that we did discuss in \cref{sec:ct}, or adjunction, that we skipped) can be studied abstractly using a suitable $2$-category known as a \emph{cosmos} \cite{street1974elementarycosmoi,street1980cosmoiinternalcategories}. An even bigger advancement is that the same perspective can be used to study $\infty$-categories via an $\infty$-categorical generalization of cosmoi, known as \emph{$\infty$-cosmoi}. This has been the insight of Riehl and Verity, who developed a whole $\infty$-categorical machinery avoiding any discussion of particular models \cite{riehlverity2017inftycosmos,riehlverity2022elements}.
  \item \textbf{A Synthetic Approach to Higher Categories:} More recently we are witnessing efforts to develop $\infty$-category theory based on some axiomatic principles, known as \emph{synthetic $\infty$-category theory}. This includes a research program \cite{cisinskicnossennguyenwalde2025higher}, as well as educational notes \cite{haugseng2025qcats}. 
\end{enumerate}

\subsection{Research via Complete Segal Spaces}
As discussed in \cref{subsec:why}, complete Segal spaces still have relevance in various aspects of cutting-edge research, including but not limited to the following areas.
\begin{enumerate}[leftmargin=*]
  \item \textbf{A Theory of Fibrations:} As we saw in \cref{ex:repcontra}, $\Hom(-,x)\colon \cC^{\op} \to \Set$ is a functor. On the other hand, for a complete Segal space $W$ a morphism $f\colon X \to Y$ induces a morphism of mapping spaces $f^*\colon\map_W(Y,Z) \to \map_W(X,Z)$ (\cref{lemma:pre post composition}), but there is no way to make it functorial, as it involves choosing an infinite tower of coherent homotopies. The standard solution employs a variety of fibrations, such as left or right fibrations. Such fibrations have been developed in many models of $\infty$-categories (in particular in the context of quasi-categories in \cite{lurie2009htt}), however, the complete Segal space approach closest matches the intuition from classical categories and is most suitable for generalizations to other contexts. See \cite[Section 1, Subsection 2.1]{rasekh2023left} for a more detailed review.
  \item \textbf{A Type Theoretic Approach to Higher Categories:} As discussed in \cref{rem:hott}, homotopy type theory is a foundation of mathematics that matches our homotopical intuition, and hence can be used to study homotopy theory. It has been an insight of Riehl and Shulman that regular homotopy type theory \cite{hottbook2013} can be generalized to a new type theory, known as \emph{simplicial homotopy type theory}. This setup admits a notion of \emph{Rezk type}, which should be thought of as an analogue to complete Segal spaces in type theory \cite{riehlshulman2017rezktypes}. This has been further advanced and extended by a variety of figures, including Buchholtz, Gratzer and Weinberger \cite{buchholtzweinberger2023cartesianfibhott,weinberger2024twosidedcartfib,weinberger2024internalsums,weinberger2022exttype,weinberger2024chev,gratzerweinbergerebuchholtz2024diruniv,gratzerweinbergerebuchholtz2025yonedaemb,gratzerweinbergerebuchholtz2026infinitycatsimp,bardomianomartinez2025limitscolimits}. Beyond being of theoretical interest, this approach has also been used to formalize $\infty$-category theory via proof assistants \cite{kudasovriehlweinberger2024rzkyoneda}.
  \item \textbf{An Internal Approach to Higher Categories:} Internal category theory studies categories that are internal to a fixed ambient category. The easiest examples are group objects, such as topological groups, which are group objects internal to the category of topological spaces. Unsurprisingly, various applications of $\infty$-category theory benefit from an internal perspective. It has long been established that the correct way to internalize $\infty$-categories is via complete Segal objects, which are analogues of complete Segal spaces internal to a fixed ambient $\infty$-category. After some early steps in this direction \cite{rasekh2018cso,rasekh2022cartesian}, this line of research has recently seen a lot of activity, both in theory, primarily due to Martini and Wolf \cite{martiniwolf2024colimits,martini2021yonedainternal,martini2022cartfibinternal,martiniwolf2022presentability,martiniwolf2023internalhighertopos}, and in applications \cite{wolf2022proetale,stenzel2025comprehension}.     
  \item \textbf{Completeness vs.~Univalence:} One central aspect of homotopy type theory is the \emph{univalence axiom}, which relates equalities and equivalences of types. From the early days of homotopy type theory, it has been expected that the univalence axiom is intimately related to the completeness condition. In recent work, using internal $\infty$-categories (in the sense of the previous item, i.e.~complete Segal objects) this question has been made precise and addressed in a variety of contexts \cite{rasekh2021univalence,stenzel2023univalence}.
  \item \textbf{Examples:} There are a variety of $\infty$-categories of interest that are constructed as complete Segal spaces. Examples include the $\infty$-category of whole-grain Petri nets \cite{kock2023wholegrainpetri} and the $\infty$-category of bordisms \cite{lurie2009cobordism,calaquescheimbauer2019cobordism}. Beyond that, a variety of $\infty$-categories can be described very explicitly using the complete Segal space perspective, such as twisted arrow $\infty$-categories \cite{mukherjeerasekh2022twisted} or the $\infty$-category of $\infty$-categories \cite{rasekh2024model}.
\end{enumerate}
    
\appendix

\section{Reviewing Category Theory} \label{sec:ct}
This appendix reviews some basic concepts from classical category theory. This is by no means a thorough review of all relevant topics. Rather it is designed to cover the basics that are used throughout this text. Here we do not provide detailed references for each statement and refer the reader to \cite[Sections 1--3]{maclane1998categories} or \cite[Chapters 1--3]{riehl2016context} for a more detailed introduction.

\subsection{Defining Categories}
We commence with the definition of a category. 
\begin{definition}
A \emph{category} $\cC$ consists of the following data:	
\begin{itemize}
    \item A class of objects $\Obj_{\cC}$.
    \item For every pair of objects $x,y$ in $\Obj_{\cC}$, a set of morphisms $\Hom_{\cC}(x,y)$.
    \item For every triple of objects $x,y,z$ in $\Obj_{\cC}$, a composition map
    \[
    - \circ - \colon \Hom_{\cC}(x,y) \times \Hom_{\cC}(y,z) \to \Hom_{\cC}(x,z)
    \]
    \item An identity morphism $\id_x$ in $\Hom_{\cC}(x,x)$ for every object $x$ in $\Obj_{\cC}$.
    \item Composition is unital: for every $f$ in $\Hom_{\cC}(x,y)$, we have $\id_y \circ f = f$ and $f \circ \id_x = f$.
		\item Composition is associative: for every $f$ in $\Hom_{\cC}(x,y)$, $g$ in $\Hom_{\cC}(y,z)$, $h$ in $\Hom_{\cC}(z,w)$, we have $h \circ (g \circ f) = (h \circ g) \circ f$.
\end{itemize}
\end{definition}

If $f$ is in $\Hom_{\cC}(x,y)$, we often write $f\colon x \to y$ and call $x$ the \emph{domain} and $y$ the \emph{codomain} of $f$. Let us observe some important examples. 

\begin{example} \label{ex:set}
  Let $\Set$ be the category defined as follows:
  \begin{itemize}
    \item $\Obj_{\Set}$ are sets.
    \item For sets $x,y$, $\Hom_{\Set}(x,y)$ are functions from $x$ to $y$.
    \item Composition and identity are defined as usual for functions.
  \end{itemize}
\end{example}

\begin{example} \label{ex:top}
  Let $\Top$ be the category defined as follows:
  \begin{itemize}
    \item $\Obj_{\Top}$ are topological spaces.
    \item For topological spaces $x,y$, $\Hom_{\Top}(x,y)$ are continuous functions from $x$ to $y$.
    \item Composition and identity are defined as usual for continuous functions.
  \end{itemize}
\end{example}

\begin{example} \label{ex:pointed top}
 We can analogously define the categories $\Set_*$ of pointed sets and the category $\Top_*$ of pointed topological spaces, where the objects are pairs $(X,x)$ with $X$ a set (resp.~topological space) and $x$ in $X$ a distinguished point, and morphisms are (continuous) functions that preserve the distinguished point.
\end{example}

\begin{example} \label{ex:ncat}
  Let $[n]$ be the category defined as follows:
  \begin{itemize}
    \item $\Obj_{[n]} = \{ 0, 1, \ldots, n \}$.
    \item For $i,j$ in $\Obj_{[n]}$, $\Hom_{[n]}(i,j)$ has a single element if $i \leq j$ and is empty otherwise.
    \item Composition and identity are defined uniquely by the above specification.
    \item Graphically, we can represent this category as follows:
  \end{itemize}
  \[
  \begin{tikzcd}
  0 \arrow[r] & 1 \arrow[r] & 2 \arrow[r] & \cdots \arrow[r] & n
  \end{tikzcd}
  \]
\end{example}

\begin{definition}
  Let $\cC$ be a category. A morphism $f\colon x \to y$ in $\cC$ is called an \emph{isomorphism} if there exists a morphism $g\colon y \to x$ in $\cC$ such that $g \circ f = \id_x$ and $f \circ g = \id_y$.
\end{definition}

\begin{example} \label{ex:iso in set}
  An isomorphism in $\Set$ (\cref{ex:set}) is a bijective function. An isomorphism in $\Top$ (\cref{ex:top}) is a homeomorphism. On the other hand, the category $[n]$ has no non-identity isomorphisms.
\end{example}

\begin{example} \label{ex:ione}
	Let $I(1)$ be the category with the following specification:
	\begin{itemize}
		\item $\Obj_{I(1)} = \{ 0, 1 \}$.
		\item $\Hom_{I(1)}(0,0) = \{ \id_0 \}$, $\Hom_{I(1)}(1,1) = \{ \id_1	\}$, $\Hom_{I(1)}(0,1) = \{ 01 \}$, 
    \item[] and	$\Hom_{I(1)}(1,0) = \{ 10 \}$.
		\item The composition is given by the identities and $10 \circ 01 = \id_0$ and $01 \circ 10 = \id_1$.
		\item Graphically, we can represent this category as follows:
	\end{itemize}
	\[
	\begin{tikzcd}
	0 \arrow[r, shift left, "01"] & 1 \arrow[l,	shift left, "10"]
	\end{tikzcd}		
\]
 We call this category the \emph{free isomorphism}.
\end{example}

We can generalize this example as follows.

\begin{example}
For $n \geq 0$, let $I(n)$ be the category with the following specification:
\begin{itemize}
  \item $\Obj_{I(n)} = \{ 0, 1, \ldots, n \}$.
  \item For $i,j$ in $\Obj_{I(n)}$, $\Hom_{I(n)}(i,j)$ has a single element.
  \item For all $i,j,k$ in $\Obj_{I(n)}$, composition is given by the unique morphisms.
\end{itemize}
Every morphism in $I(n)$ is an isomorphism, as the inverse of the unique morphism from $i$ to $j$ is the unique morphism from $j$ to $i$. 
\end{example}

\begin{definition}
  A \emph{groupoid} is a category in which every morphism is an isomorphism.
\end{definition}

\begin{example} \label{ex:groupoids}
  Following \cref{ex:iso in set}, the categories $\Set$, $\Top$, $[n]$ are not groupoids (if $n > 0$), whereas the $I(n)$ are groupoids.
\end{example}

If a category is not a groupoid, then we would like to be able to obtain one.

\begin{proposition} \label{prop:core}
  Let $\cC$ be a category. Then there is a groupoid $\cC^{\simeq}$, called the \emph{core} or \emph{maximal subgroupoid} of $\cC$, which has the same objects as $\cC$ and as morphisms the isomorphisms in $\cC$. Moreover, the core $(-)^{\simeq}$ defines a functor $(-)^{\simeq}\colon \Cat \to \Cat$.
\end{proposition}

Of course, we need to be able to relate different categories to each other. For that we have the notion of functors.

\begin{definition}
	Let $\cC,\cD$ be two categories. A \emph{functor} $F\colon \cC \to \cD$ consists of the following data:
	\begin{itemize}
		\item A function $F\colon \Obj_{\cC} \to \Obj_{\cD}$.
		\item For every pair of objects $x,y \in	\Obj_{\cC}$, a function $F\colon \Hom_{\cC}(x,y) \to \Hom_{\cD}(F(x),F(y))$.
		\item For every object $x$ in $\Obj_{\cC}$, we	have $F(\id_x) = \id_{F(x)}$.
		\item	For every triple of objects $x,y,z$ in $\Obj_{\cC}$, and every $f\colon x \to y, g\colon y \to z$, we have $F(g \circ f) = F(g) \circ F(f)$.	
	\end{itemize}
\end{definition}

Functors can be composed and we have identity functors. Thus, categories and functors themselves form a category.

\begin{example} \label{ex:cat}
 Let $\Cat$ be the category with the following specification:
	\begin{itemize}
		\item $\Obj_{\Cat}$ is the class of all categories.
		\item For $\cC,\cD$ in $\Obj_{\Cat}$, $\Hom_{\Cat}(\cC,\cD)$ is the set of all functors from $\cC$ to $\cD$.
		\item The composition is given by the composition of functors.
		\item The identity is given by the usual identity functor.
	\end{itemize}
\end{example}

Finally, we make the following observation about functors.
\begin{lemma} \label{lemma:functor as sequence of morphisms}
  Let $\cC$ be a category and let $n \geq 0$. A functor $F\colon [n]\to \cC$ is precisely a sequence of $n$ composable morphisms in $\cC$.
\end{lemma}

\begin{proof}
  By definition of $[n]$ (\cref{ex:ncat}), a functor $F\colon [n] \to \cC$ consists of the following data:
  \begin{itemize}
    \item For every $i \in \{ 0,1,\ldots,n \}$, an object $F(i)$ in $\cC$.
    \item For every pair of objects $i,j$ in $[n]$, a morphism $f_{ij} = F(i \to j)\colon F(i) \to F(j)$ in $\cC$ if $i \leq j$ and no morphism otherwise.
  \end{itemize}
  The identity condition implies $f_{ii}$ is the identity. Moreover, the composition condition implies if $i + 1 < j$, $f_{ij} = f_{i+1,j} \circ f_{i,i+1}$. Hence, the functor $F$ on morphisms is uniquely determined by $f_{i,i+1}\colon F(i) \to F(i+1)$, giving us the desired result.
\end{proof}

We can similarly get the following result for $I(n)$.

\begin{lemma} \label{lemma:functor as sequence of isomorphisms}
  Let $\cC$ be a category and let $n \geq 0$. A functor $F\colon I(n)\to \cC$ is precisely a sequence of $n$ composable isomorphisms in $\cC$.
\end{lemma}

\begin{remark} \label{rem:n to in lift}
  \cref{lemma:functor as sequence of morphisms,lemma:functor as sequence of isomorphisms} imply that for a functor $F\colon [n] \to \cC$ if the associated morphisms $f_i$ in the tuple $(f_i)_{1 \leq i \leq n}$ are isomorphisms in $\cC$, then $F$ lifts to a functor $\hat{F}\colon I(n) \to \cC$. Here the additional morphisms in $I(n)$ are sent to the inverses. 
  So, in particular, if $\cC$ is a groupoid, there is a bijection between functors $[n] \to \cC$ and functors $I(n) \to \cC$, and composable morphisms of length $n$ in $\cC$ are precisely composable isomorphisms of length $n$ in $\cC$.
\end{remark}

\subsection{Functor Categories}
One important aspect of category theory is its ability to compare functors via natural transformations.

\begin{definition}
	Let $\cC,\cD$ be two categories and let $F,G\colon \cC \to \cD$ be two functors. A \emph{natural transformation} $\eta\colon F \to G$ consists of the following data:
	\begin{itemize}
		\item For every object $x$ in $\Obj_{\cC}$, a morphism $\eta_x\colon F(x) \to G(x)$ in $\cD$.
		\item For every morphism $f\colon x \to y$ in $\cC$, the following diagram commutes:
		\[
		\begin{tikzcd}
		F(x) \arrow[r, "F(f)"] \arrow[d, "\eta_x"'] & F(y) \arrow[d, "\eta_y"] \\
		G(x) \arrow[r, "G(f)"'] & G(y)
		\end{tikzcd}
		\]
	\end{itemize}
  It is a \emph{natural isomorphism} if each $\eta_x$ is an isomorphism in $\cD$.
\end{definition}

Using natural transformations, functors themselves form a category.

\begin{definition} \label{def:functor category}
	Let $\cC,\cD$ be two categories. We define the \emph{functor category} $\Fun(\cC,\cD)$ with the following specification:
	\begin{itemize}
		\item $\Obj_{\Fun(\cC,\cD)}$ are functors from $\cC$ to $\cD$.
		\item For $F,G$ in $\Obj_{\Fun(\cC,\cD)}$, we have $\Hom_{\Fun(\cC,\cD)}(F,G) = \Nat(F,G)$, the set of natural transformations from $F$ to $G$.
	\end{itemize}
\end{definition}

\begin{example} \label{ex:functor category from point}
  The functor category $\Fun([0],\cC)$ is isomorphic to $\cC$. Indeed, a functor $F\colon [0] \to \cC$ is just a choice of object (\cref{lemma:functor as sequence of morphisms}) and a natural transformation between two objects is just a morphism between the corresponding objects in $\cC$. Hence, $\Fun([0],\cC) \cong \cC$.
\end{example}

The functor category is functorial, in the following sense.

\begin{proposition} \label{prop:functor category functorial}
  Let $\cE$ be a category. Then the assignment $\cD \mapsto \Fun(\cD,\cE)$ extends to a functor $\Fun(-,\cE)\colon \Cat^{\op} \to \Cat$, which sends a functor $F\colon \cC \to \cD$ to the functor $F^*\colon \Fun(\cD,\cE) \to \Fun(\cC,\cE)$, given by precomposition with $F$.
\end{proposition}

Functor categories can be iterated. For that we need product categories.

\begin{definition} \label{def:product category}
  Let $\cC,\cD$ be two categories. The \emph{product category} $\cC \times \cD$ is defined as follows:
\begin{itemize}
    \item $\Obj_{\cC \times \cD} = \Obj_{\cC} \times \Obj_{\cD}$.
    \item For objects $(x_1,y_1),(x_2,y_2)$ in $\Obj_{\cC \times \cD}$, $\Hom_{\cC \times \cD}((x_1,y_1),(x_2,y_2)) = \Hom_{\cC}(x_1,x_2) \times \Hom_{\cD}(y_1,y_2)$.
    \item The composition and identity are defined component-wise.
\end{itemize}
\end{definition}

\begin{proposition} \label{prop:functor category as currying}
  Let $\cC,\cD,\cE$ be three categories. Then there is an isomorphism of categories
  \[
  \Fun(\cC \times \cD, \cE) \cong \Fun(\cC, \Fun(\cD,\cE)).
  \]
\end{proposition}

\begin{proof}
 Given a functor $F\colon \cC \to \Fun(\cD,\cE)$, we define a functor $\hat{F}\colon \cC \times \cD \to \cE$ as follows:
\begin{itemize}[leftmargin=*]
  \item For every object $(x,y)$ in $\Obj_{\cC \times \cD}$, we have $\hat{F}(x,y) = F(x)(y)$.
  \item For every morphism $(g,f)\colon (x,y) \to (x',y')$ in $\cC \times \cD$, $\hat{F}(g,f)$ is the composite
  \[
    F(x)(y) \xrightarrow{F(x)(f)} F(x)(y') \xrightarrow{F(g)_{y'}} F(x')(y').
  \]
\end{itemize}
Similarly, for every natural transformation $\alpha\colon F \to G$ between functors $F,G\colon \cC \to \Fun(\cD,\cE)$, we define a natural transformation $\hat{\alpha}\colon \hat{F} \to \hat{G}$ as follows. For every object $(x,y)$ in $\Obj_{\cC \times \cD}$, we have $\hat{\alpha}_{(x,y)} = \alpha_x(y)$. This gives us a functor $\Fun(\cC, \Fun(\cD,\cE)) \to \Fun(\cC \times \cD, \cE)$.

It is immediate from the construction that this assignment is injective and surjective on objects and morphisms, hence is an isomorphism of categories.
\end{proof}

\subsection{The Yoneda Lemma}
We will now move on to a functor category of interest, which requires opposite categories. 

\begin{definition}
	Let $\cC$ be a category. The \emph{opposite category} $\cC^{\op}$ is defined as follows:
	\begin{itemize}
		\item $\Obj_{\cC^{\op}} = \Obj_{\cC}$.
		\item For objects $x,y$ in $\cC^{\op}$, $\Hom_{\cC^{\op}}(x,y) = \Hom_{\cC}(y,x)$.
		\item The composition and identity are defined as in $\cC$, but with the direction of morphisms reversed.
	\end{itemize}
\end{definition}

\begin{example} \label{ex:repcontra}
  Let $\cC$ be a category and let $x$ be an object in $\cC$. The \emph{representable functor} 
  \[\Hom_{\cC}(-,x)\colon \cC^{\op} \to \Set\] 
  is defined as follows:
  \begin{itemize}
    \item For every object $y$ in $\Obj_{\cC}$, we have $\Hom_{\cC}(y,x)$.
    \item For every morphism $f\colon y \to z$ in $\cC$, the map 
    \[
    f^*\colon \Hom_{\cC}(z,x) \to \Hom_{\cC}(y,x), \quad g \mapsto g \circ f.
    \]
  \end{itemize}
\end{example}

\begin{example}
 Similar to \cref{ex:repcontra}, we can define the \emph{corepresentable functor} 
 \[\Hom_{\cC}(x,-)\colon \cC \to \Set,\] 
 which maps an object $y$ to the hom-set $\Hom_{\cC}(x,y)$ and a morphism $f\colon y \to z$ in $\cC$ to the post-composition map $f_* \colon \Hom_{\cC}(x,y) \to \Hom_{\cC}(x,z)$, given by $g \mapsto f \circ g$.
\end{example}

In fact the opposite category is a functor.

\begin{lemma} \label{lemma:opposite functor}
  The opposite category construction extends to a functor $(-)^{\op}\colon \Cat \to \Cat$.
\end{lemma}

\begin{proof}
  For a given functor $F\colon \cC \to \cD$, we define $F^{\op}\colon \cC^{\op} \to \cD^{\op}$ as $F$ on objects and morphisms (where the domain and codomain have been exchanged).
\end{proof}

\begin{theorem}[Yoneda Lemma] \label{thm:yoneda}
	Let $\cC$ be a category, let $x$ be an object in $\cC$, and let $F\colon \cC^{\op} \to \Set$ be a functor. Then the map
	\[
 \Nat(\Hom_{\cC}(-,x),F) \to F(x), \quad \alpha \mapsto \alpha_x(\id_x)
	\]
  is a bijection of sets.
\end{theorem}

\begin{definition} \label{def:yoneda embedding}
  Let $\cC$ be a category. The \emph{Yoneda embedding} is the functor denoted $\Yon\colon\cC \to \Fun(\cC^{\op},\Set)$ and defined as follows:
  \begin{itemize}
   \item On objects: for every object $c$ in $\Obj_{\cC}$, $\Yon(c) = \Hom_{\cC}(-,c)$.
   \item On morphisms: for every morphism $f\colon x \to y$ in $\cC$, we define the natural transformation $\Yon(f)\colon \Hom_{\cC}(-,x) \to \Hom_{\cC}(-,y)$ to be the natural transformation given by postcomposition with $f$.
  \end{itemize}
\end{definition}

We in particular have the following consequence of the Yoneda lemma.

\begin{corollary} \label{cor:iso iff natural iso}
  Let $f \colon x \to y$ be a morphism in a category $\cC$. Then $f$ is an isomorphism if and only if the natural transformation $\Yon(f)\colon \Hom_{\cC}(-,x) \to \Hom_{\cC}(-,y)$ is a natural isomorphism. 
\end{corollary}

We can use the Yoneda embedding to construct functors into functor categories. 

\begin{construction} \label{constr:restricted yoneda}
 Let $G\colon \cC \to \cD$ be a functor. The restricted Yoneda functor of $G$ is the composition
 \[\hat{G}\colon \cD \xrightarrow{\Yon} \Fun(\cD^{\op},\Set) \xrightarrow{G^*} \Fun(\cC^{\op},\Set).\]
\end{construction}

\subsection{Equivalences of Categories}
We can use natural transformations to understand when two categories are equivalent.

\begin{definition} \label{def:equivalence of categories}
	Let $\cC,\cD$ be two categories. A functor $F\colon \cC \to \cD$ is an \emph{equivalence of categories} if there exists a functor $G\colon \cD \to \cC$ such that $G \circ F$ is naturally isomorphic to $\id_{\cC}$ and $F \circ G$ is naturally isomorphic to $\id_{\cD}$.
\end{definition}

\begin{definition} \label{def:essentially surjective}
  A functor $F\colon \cC \to \cD$ is called \emph{essentially surjective} if for every object $d$ in $\Obj_{\cD}$, there exists an object $c$ in $\Obj_{\cC}$ such that $F(c)$ is isomorphic to $d$ in $\cD$.
\end{definition}

\begin{definition} \label{def:fully faithful}
  A functor $F\colon \cC \to \cD$ is called \emph{fully faithful} if for every pair of objects $X,Y$ in $\Obj_{\cC}$, the map 
  \[
  F_{X,Y}\colon \Hom_{\cC}(X,Y) \to \Hom_{\cD}(F(X),F(Y))
  \]
  is a bijection of sets.
\end{definition}

\begin{theorem} \label{thm:equivalence of categories}
	Let $\cC,\cD$ be two categories and let $F\colon \cC \to \cD$ be a functor. The following are equivalent:
	\begin{itemize}
		\item $F$ is an equivalence of categories.
		\item $F$ is fully faithful and essentially surjective. 
	\end{itemize}
\end{theorem}

\begin{example} \label{ex:i equivalence}
  Let $n \geq 1$. The functor $0\colon [0] \to I(n)$ is an equivalence of categories, as it is fully faithful and $\{\id_0\} = \Hom_{[0]}(0,0) = \Hom_{I(n)}(0,0) = \{\id_0\}$, and essentially surjective, as every object in $I(n)$ is isomorphic to $0$.
\end{example}

\subsection{Limits} \label{subsec:limits}
We end this section with a discussion of colimits and limits. The theory of limits is one of the most fundamental aspects of category theory, and merits a detailed discussion, where we suggest \cite[Section V]{maclane1998categories}, \cite[Chapter 3]{riehl2016context} for a thorough introduction to this subject. Here we will exclusively review specific examples relevant to this text.

\begin{definition} \label{def:colimit}
  Let 
  \[A_1 \xleftarrow{f_{1,2}} A_{1,2} \xrightarrow{g_{1,2}} A_2 \xleftarrow{f_{2,3}} A_{2,3} ... \xleftarrow{f_{n-1,n}} A_{n-1,n}\xrightarrow{g_{n-1,n}}  A_n \]
   be a diagram of sets. The \emph{colimit} of this diagram is the set 
  \[A_1 \coprod_{A_{1,2}}^{f_{1,2},g_{1,2}} ... \coprod_{A_{n-1,n}}^{f_{n-1,n},g_{n-1,n}} A_n = (\coprod_{1 \leq i \leq n}A_i) / \sim,\]
  where $\coprod$ is the disjoint union of sets and  $\sim$ is the equivalence relation generated by $f_{i,i+1}(a) \sim g_{i,i+1}(a)$ for every $a \in A_{i,i+1}$.
\end{definition}

\begin{definition} \label{def:limit}
  Let 
  \[A_1 \xrightarrow{f_{1,2}} A_{1,2} \xleftarrow{g_{1,2}} A_2 \xrightarrow{f_{2,3}} A_{2,3} ... \xrightarrow{f_{n-1,n}} A_{n-1,n}\xleftarrow{g_{n-1,n}}  A_n \]
   be a diagram of sets. The \emph{limit} of this diagram is the set 
   {\small
  \[A_1 \times_{A_{1,2}}^{f_{1,2},g_{1,2}} ... \times_{A_{n-1,n}}^{f_{n-1,n},g_{n-1,n}} A_n = \{(a_i)_{1 \leq i \leq n} \in \prod_{1 \leq i \leq n}A_i \mid f_{i,i+1}(a_i) = g_{i,i+1}(a_{i+1}) \text{ for all } 1 \leq i \leq n-1 \}\]
   }
  where $\prod$ is the Cartesian product of sets.
\end{definition}

There are several special cases of interest.

\begin{example} \label{ex:pushout}
 If $n = 2$, then the colimit $A_1 \coprod_{A_{1,2}}^{f_{1,2},g_{1,2}} A_2$ is called the \emph{pushout} of the diagram $A_1 \xleftarrow{f_{1,2}} A_{1,2} \xrightarrow{g_{1,2}} A_2$.  
\end{example}

\begin{example} \label{ex:pullback}
  If $n = 2$, then the limit $A_1 \times_{A_{1,2}}^{f_{1,2},g_{1,2}} A_2$ is called the \emph{pullback} of the diagram $A_1 \xrightarrow{f_{1,2}} A_{1,2} \xleftarrow{g_{1,2}} A_2$.
\end{example}

\begin{example} \label{ex:coproduct}
 If all the sets $A_{i,i+1}$ in the diagram of \cref{def:colimit} are empty, then the colimit is just the disjoint union of the sets $A_i$ and is also called the \emph{coproduct}.
\end{example}

\begin{example} \label{ex:product}
 If all the sets $A_{i,i+1}$ in the diagram of \cref{def:limit} have one element, then the limit is just the Cartesian product of the sets $A_i$ and is also called the \emph{product}.
\end{example}

We can now generalize these notions from sets to functor categories.

\begin{definition} \label{def:colimit functor}
  Let $\cC$ be a category and let 
  \[F_1 \xleftarrow{\alpha_{1,2}} F_{1,2} \xrightarrow{\beta_{1,2}} F_{2} \xleftarrow{\alpha_{2,3}}  ... \xrightarrow{\beta_{n-1,n}} F_n\] 
  be a diagram in $\Fun(\cC^{\op},\Set)$. Then the \emph{colimit} of this diagram is the functor 
  \[F_1 \coprod_{F_{1,2}}^{\alpha_{1,2},\beta_{1,2}} ... \coprod_{F_{n-1, n}}^{\alpha_{n-1, n},\beta_{n-1, n}} F_n\colon \cC^{\op} \to \Set\] 
  defined as follows:
  \[(F_1 \coprod_{F_{1,2}}^{\alpha_{1,2},\beta_{1,2}} ... \coprod_{F_{n-1, n}}^{\alpha_{n-1, n},\beta_{n-1, n}} F_n)(c) = F_1(c) \coprod_{F_{1,2}(c)}^{\alpha_{1,2}(c),\beta_{1,2}(c)} ... \coprod_{F_{n-1, n}(c)}^{\alpha_{n-1, n}(c),\beta_{n-1, n}(c)} F_n(c),\]
  meaning the value of the colimit functor at an object $c$ is the colimit of the diagram obtained by evaluating the original diagram at $c$.
\end{definition}

\begin{definition} \label{def:limit functor}
  Let $\cC$ be a category and let 
  \[F_1 \xrightarrow{\alpha_{1,2}} F_{1,2} \xleftarrow{\beta_{1,2}} F_{2} \xrightarrow{\alpha_{2,3}}  ... \xleftarrow{\beta_{n-1,n}} F_n\] 
  be a diagram in $\Fun(\cC^{\op},\Set)$. Then the \emph{limit} of this diagram is the functor 
  \[F_1 \times_{F_{1,2}}^{\alpha_{1,2},\beta_{1,2}} ... \times_{F_{n-1, n}}^{\alpha_{n-1, n},\beta_{n-1, n}} F_n\colon \cC^{\op} \to \Set\] 
  defined as follows:
  \[(F_1 \times_{F_{1,2}}^{\alpha_{1,2},\beta_{1,2}} ... \times_{F_{n-1, n}}^{\alpha_{n-1, n},\beta_{n-1, n}} F_n)(c) = F_1(c) \times_{F_{1,2}(c)}^{\alpha_{1,2}(c),\beta_{1,2}(c)} ... \times_{F_{n-1, n}(c)}^{\alpha_{n-1, n}(c),\beta_{n-1, n}(c)} F_n(c),\]
  meaning the value of the limit functor at an object $c$ is the limit of the diagram obtained by evaluating the original diagram at $c$.
\end{definition}

\begin{remark}
  Analogous to \cref{ex:pushout,ex:pullback,ex:coproduct,ex:product}, we similarly have pushouts, pullbacks, coproducts, and products in $\Fun(\cC^{\op},\Set)$ as special cases of the above definitions.
\end{remark}

We can use pushouts and pullbacks to define a notion of a pullback or pushout square in $\Fun(\cC^{\op},\Set)$.
\begin{definition} \label{def:pullback square}
  If we are given a commutative square in $\Fun(\cC^{\op},\Set)$
  \[
  \begin{tikzcd}
    F \arrow[r, "\alpha"] \arrow[d, "\beta"'] \arrow[dr, phantom, "\ulcorner", very near start] & G \arrow[d, "\gamma"] \\
    H \arrow[r, "\delta"'] & K
  \end{tikzcd}.
  \]
  Then, we say it is a \emph{pullback square} if $F \cong G \times_K H$, and denote the square via $\ulcorner$ in the top left corner.
\end{definition}

\begin{definition} \label{def:pushout square}
  If we are given a commutative square in $\Fun(\cC^{\op},\Set)$
  \[
  \begin{tikzcd}
    F \arrow[r, "\alpha"] \arrow[d, "\beta"'] & G \arrow[d, "\gamma"] \\
    H \arrow[r, "\delta"'] & K  \arrow[ul, phantom, "\ulcorner", very near start]
  \end{tikzcd}.
  \]
  Then, we say it is a \emph{pushout square} if $K \cong G \coprod_F H$, and denote the square via $\ulcorner$ in the bottom right corner.
\end{definition}

\begin{remark} \label{rem:products functor}
  Products in $\Fun(\cC,\Set)$ give us a functor $X \times - \colon \Fun(\cC,\Set) \to \Fun(\cC,\Set)$ for every $X$ in $\Fun(\cC,\Set)$, which maps an object $Y$ to $X \times Y$, and a morphism $f\colon Y \to Z$ to the morphism $\id_X \times f \colon X \times Y \to X \times Z$.
\end{remark}
\begin{notation}
 If it is clear from the context we will denote the colimit by $F_1 \coprod_{F_{1,2}} ... \coprod_{F_{n-1, n}} F_n$ and the limit by $F_1 \times_{F_{1,2}} ... \times_{F_{n-1, n}} F_n$.  
\end{notation}

The notions of colimit and limit in $\Fun(\cC^{\op},\Set)$ interact well with limits in $\Set$.

\begin{proposition} \label{prop:colimit limit}
  Let $\cC$ be a category, 
  \[F_1 \xleftarrow{\alpha_{1,2}} F_{1,2} \xrightarrow{\beta_{1,2}} F_{2} \xleftarrow{\alpha_{2,3}}  ... \xrightarrow{\beta_{n-1,n}} F_n\] 
  be a diagram in $\Fun(\cC^{\op},\Set)$, and let $C\colon \cC^{\op} \to \Set$ denote the colimit of this diagram. Moreover, for a functor  $P\colon \cC^{\op} \to \Set$ let the set $L$ denote the limit of the diagram of sets 
  \[\Hom(F_1,P) \xrightarrow{\alpha_{1,2}^*} ... \xleftarrow{\beta_{n-1,n}^*} \Hom(F_n,P).\] 
  Then $L \cong \Hom(C,P)$. 
\end{proposition}

\begin{proposition} \label{prop:limit limit}
  Let $\cC$ be a category,
  \[F_1 \xrightarrow{\alpha_{1,2}} F_{1,2} \xleftarrow{\beta_{1,2}} F_{2} \xrightarrow{\alpha_{2,3}}  ... \xleftarrow{\beta_{n-1,n}} F_n\] 
  be a diagram in $\Fun(\cC^{\op},\Set)$, and let $L\colon \cC^{\op} \to \Set$ denote the limit of this diagram. Moreover, for a functor  $P\colon \cC^{\op} \to \Set$ let $L'$ denote the limit of the diagram of sets 
  \[\Hom(P,F_1) \xrightarrow{(\alpha_{1,2})_*} ... \xleftarrow{(\beta_{n-1,n})_*} \Hom(P,F_n).\] 
  Then $L' \cong \Hom(P,L)$.  
\end{proposition}

\begin{remark} \label{rem:colimit limit}
  \cref{prop:colimit limit} is stating that a natural transformation $\gamma \colon F_1 \coprod_{F_{1,2}} ... \coprod_{F_{n-1,n}} F_n \to P$ is uniquely determined by a tuple of natural transformations $\gamma_{F_i}\colon F_i \to P$ such that $\gamma_{F_i} \circ \alpha_{i,i+1} = \gamma_{F_{i+1}} \circ \beta_{i,i+1}$ for all $i$, which is precisely an element in the set on the right-hand side. Similarly, \cref{prop:limit limit} is stating that a natural transformation $\gamma\colon P \to F_1 \times_{F_{1,2}} ... \times_{F_{n-1,n}} F_n $ is uniquely determined by a tuple of natural transformations $\gamma_{F_i}\colon P \to F_i$ such that $\alpha_{i,i+1} \circ \gamma_{F_i} = \beta_{i,i+1} \circ \gamma_{F_{i+1}}$ for all $i$.
\end{remark}

One benefit of developing a theory of limits is the ability to obtain new pullbacks and maps. Here are several examples of such results. The first is the functoriality of pullbacks.

\begin{lemma} \label{lemma:limit map}
 Let $\cC$ be a category and assume we have the following commutative diagram in the category $\Fun(\cC^{\op},\Set)$
  \[ 
  \begin{tikzcd}
   A_1 \arrow[r] \arrow[d, "f_1"'] & A_{1,2} \arrow[d, "f_{1,2}"'] & ... \arrow[r] \arrow[l] & A_{n-1,n} \arrow[d, "f_{n-1,n}"'] & A_n \arrow[d, "f_n"'] \arrow[l] \\ 
   B_1 \arrow[r] & B_{1,2} & ... \arrow[r] \arrow[l] & B_{n-1,n} & B_n \arrow[l]
  \end{tikzcd},
  \]
  Then there exists an induced morphism on limits 
  \[(f_1,...,f_n) \colon A_1 \times_{A_{1,2}} ... \times_{A_{n-1,n}} A_n \to B_1 \times_{B_{1,2}} ... \times_{B_{n-1,n}} B_n,\]
  sending a tuple $(a_1,...,a_n)$ to the tuple $(f_1(a_1),...,f_n(a_n))$. Moreover, if the $f_i$ and $f_{i,i+1}$ are isomorphisms for all $i$, then the induced morphism on limits is also an isomorphism.
\end{lemma}

The second one is known as ``pullback pasting'' and ``pullback cancellation''.

\begin{lemma} \label{lemma:pullback pasting}
 Let $\cC$ be a category and assume we have the following commutative diagram in the category $\Fun(\cC^{\op},\Set)$
  \[
  \begin{tikzcd}
  A \arrow[r] \arrow[d] & B \arrow[r, "f"] \arrow[d, "g"] & C \arrow[d, "h"]\\ 
  D \arrow[r, "k"] & E \arrow[r, "l"] & F
  \end{tikzcd}
  \]
  If the right-hand square is a pullback square, then the left-hand square is a pullback square if and only if the outer rectangle is a pullback square.
\end{lemma}

This last one is known as ``pullback commutativity'', which is a special case of commuting limits \cite[Section IX.8]{maclane1998categories}.

\begin{lemma} \label{lemma:pullback commutativity}
 Let $\cC$ be a category and assume we have the following commutative diagram in the category $\Fun(\cC^{\op},\Set)$
  \[
  \begin{tikzcd}
  A \arrow[r] \arrow[d] & B \arrow[d] & C \arrow[l] \arrow[d] \\ 
  D \arrow[r] & E & F \arrow[l] \\
  G \arrow[r] \arrow[u] & H \arrow[u] & I \arrow[l] \arrow[u] 
  \end{tikzcd}.
  \]
  Then $(A \times_D G) \times_{(B \times_E H)} (C \times_F I) \cong (A \times_B C) \times_{(D \times_E F)} (G \times_H I)$. This means first taking the pullbacks of the rows and then the pullback of the resulting column coincides with first taking the pullbacks of the columns and then the pullback of the resulting row.
\end{lemma}

\bibliographystyle{alpha}
\bibliography{main}

\end{document}